\documentclass{article}
\usepackage{graphicx,epstopdf}
\usepackage{amsmath,amssymb,amsthm}
\usepackage{mathtools}
\usepackage{dsfont}
\usepackage{subfig}
\usepackage{epstopdf}
\usepackage{geometry}
\usepackage{booktabs}
\usepackage{longtable}
\usepackage{graphicx}
\usepackage{epstopdf}
\usepackage{comment}
\usepackage{hyperref}
\allowdisplaybreaks[1]
\numberwithin{equation}{section}
\newtheorem{theorem}{Theorem}[section]
\newtheorem{lemma}{Lemma}[section]

\normalsize
\large

\begin{document}
\title{ANALYSIS OF A FEM-MCM DISCRETIZATION FOR THE 2D/3D STOCHASTIC CLOSED-LOOP GEOTHERMAL SYSTEM
\thanks{Supported in part by Key project of National Natural Science Foundation (No. 12431016), and NSF of China (No. 12471408 and No. 12001347), Shaanxi Provincial Joint Laboratory of Artificial Intelligence (No. 2022JC-SYS-05), National program for the introduction of high-end foreign experts (No. G2023041032L), Innovative team project of Shaanxi Provincial Department of Education (No. 21JP013 and No. 21JP019) and Shaanxi Province Natural Science basic research program key project (No. 2023-JC-ZD-02).}}
\author{Yi Qin\thanks{School of Mathematics and Data Science, Shaanxi University of science and technology, Xi'an, Shaanxi 710016, China.({\tt 4545@sust.edu.cn})}
, Xinyue Gao\thanks{School of Mathematics and Data Science, Shaanxi University of science and technology, Xi'an, Shaanxi 710016, China.({\tt 210911022@sust.edu.cn})}
, Lele Chen\thanks{School of Mathematics and Data Science, Shaanxi University of science and technology, Xi'an, Shaanxi 710016, China.({\tt 221711036@sust.edu.cn})}
, Liangliang Jiang\thanks{Department of Chemical and Petroleum Engineering, University of Calgary, Calgary, Alberta, Canada.({\tt jial@ucalgary.ca})}
, Zhangxin Chen\thanks{Department of Chemical and Petroleum Engineering, University of Calgary, Calgary, AB T2N 1N4, Canada.
({\tt zhachen@ucalgary.ca})}
\ and Jian Li\thanks{School of Mathematics and Data Science, Shaanxi University of science and technology, Xi'an, Shaanxi 710016, China.({\tt jianli@sust.edu.cn})}
}
\date{}
\maketitle
{\bf Abstract. This paper develops a new 2D/3D stochastic closed-loop geothermal system with a random hydraulic conductivity tensor. We use the finite element method (FEM) and the Monte Carlo method (MCM) to discrete physical and probability spaces, respectively. This FEM-MCM method is effective. The stability for velocity and temperature is rigorously proved. Compared with the deterministic closed-loop geothermal system, the same optimal error estimate for approximate velocity and temperature is obtained. Furthermore, a series of numerical experiments were carried out to show this method has better stability and accuracy results.}\par
\par
\section{Introduction}
~\par With the increasing demand for energy worldwide, it is important to accelerate the exploration of renewable energy technologies to alleviate energy tensions. Compared with other renewable energy sources, geothermal energy has great development potential in power generation, heating, and comprehensive utilization, etc. Choosing a reasonable way to efficiently develop and utilize geothermal resources is a hot topic of research in recent years around the world \cite{M. Asif-2007,M Sun-2012,J. W. Lund-1996,I. B. Fridleifsson-2001,Y. Shi-2018}.  The closed-loop heat exchanger is an essential means of extracting geothermal energy. Its basic working principle is to circulate working fluid down to the rock mass with high temperatures and then back to the surface through a continuous and closed-loop pipe \cite{C M Oldenburg-2016,B S Wu-2017,A Bensoussan-1973}.\par
In a deterministic closed-loop geothermal system, both the flow of fluid and Darcy's law are coupled with heat transfer in the pipe region and reservoirs, respectively. Also, clean geothermal energy has no direct material exchange in pipelines and reservoirs \cite{M A A Mahbub-2020,Y Qin-2022,W. Zhang-2022}. For deterministic partial differential equation models \cite{L. L. Cao-2022,J. Li-2019,J. Li-2022,J. Li-2023,L. L. Cao-2021,J. Li-2020,R. Li-2018}, such as model coefficients, forcing terms, domain geometry, boundary conditions, and initial conditions, etc., are assumed to be completely known. However, in real life, due to the existence of measurement noise, the uncertainty of these real data makes it impossible to use deterministic partial differential equations to model accurately. Therefore, the stochastic partial differential equation model has attracted the attention of many experts. For this type of model, there have been many effective numerical methods, such as stochastic Galerkin methods \cite{R. Tempone-2004,H. G. Matthies-2005,L. J. Roman-2006}, stochastic configuration methods \cite{L. Guo-2017,F. Nobile-2010}, sparse grid methods \cite{F. Bao-2014,Z. Morrow-2020,F. Nobile-2008,R. Tempone-2008}, perturbation methods \cite{C. L. Winter-2002,R. A. Todor-2005} and many other methods \cite{M. Baccouch-2020,J. Lyu-2020,A. Tambue-2019,D. Wang-2021,Z. Wang-2018,X. Wei-2021,D. Zhang-2020}. The MCM is one of the most classical approaches for solving stochastic partial differential equations. The main idea of this method is to apply numerical simulation for many times, and then obtain the desired solution by averaging all of them \cite{G. S. Fishman-1996}. In addition, to our knowledge, there are currently few efficient way of solving the stochastic coupled free/porous system in 2D/3D has been developed. Moreover, the existing methods have high computational expense in terms of random and 3D computing. There are still many computational challenges in modeling this kind of complex system.\par
In this paper, a 2D/3D stochastic closed-loop geothermal system with random hydraulic conductivity tensor is investigated. Although the stochasticity comes only from the uncertainty of the hydraulic conductivity tensor in the porous media region, the solutions obtained are stochastic in the whole region due to the coupling system. We use the FEM-MCM method to obtain an approximate solution for a 2D/3D stochastic closed-loop geothermal system. The method consists of repeated sampling of the hydraulic conduction tensor and solving the corresponding deterministic closed-loop geothermal system, which generates identically distributed approximations of the solution. In theory, stability and convergence of the velocity and temperature of the 2D/3D stochastic closed-loop geothermal system were proved. It is worth mentioning that we obtain the same optimal convergence order as the deterministic closed-loop geothermal system. Finally, several numerical experiments, including 2D and 3D cases are performed to show validation of the proposed model and algorithms.\par
This paper is organized as follows. In Section 2, we introduce the 2D/3D stochastic closed-loop geothermal system with a random hydraulic conductivity tensor and some basics. The decoupled algorithm for the proposed model is given in Section 3. Then, in Section 4, stability and error estimate are obtained for the 2D/3D stochastic closed-loop geothermal system. Finally, numerical experiments are provided for the 2D/3D stochastic closed-loop geothermal system in Section 5.
\section{Model Foundation}
~\par The stochastic closed-loop geothermal system is defined in bounded domain $\Omega'\subset R^{d}(d=2,3)$. As illustrated in Figure 1, the global domain
 $\Omega'$ consists of two subdomains, $\Omega_{f}$ and $\Omega_{p}$. Here $\Omega_{f} \bigcap \Omega_{p}=\emptyset$, $\bar{\Omega}_{f} \bigcap \ \bar{\Omega}_{p}=\uppercase\expandafter{\romannumeral1}$ and $\bar{\Omega}_{f} \bigcup \bar{\Omega}_{p}=\overline{\Omega'}$.
$\mathbf{n}_{f}$ and $\mathbf{n}_{p}$ as usual the unit outward normal directions on $\partial\Omega_{f}$ and $\partial\Omega_{p}$. The time frame is considered in $[0,T]$.
Let $(\Omega,\mathcal{F},\mathbb{P})$ be a complete probability space. Here $\Omega$ is the set of outcomes, $\mathcal{F}\in2^{\Omega}$ is the $\sigma$-algebra of events, and $\mathbb{P}:\mathcal{F}\rightarrow [0, 1]$ is a probability measure.

\begin{figure}[ht]
  \centering
  \includegraphics[width=11cm]{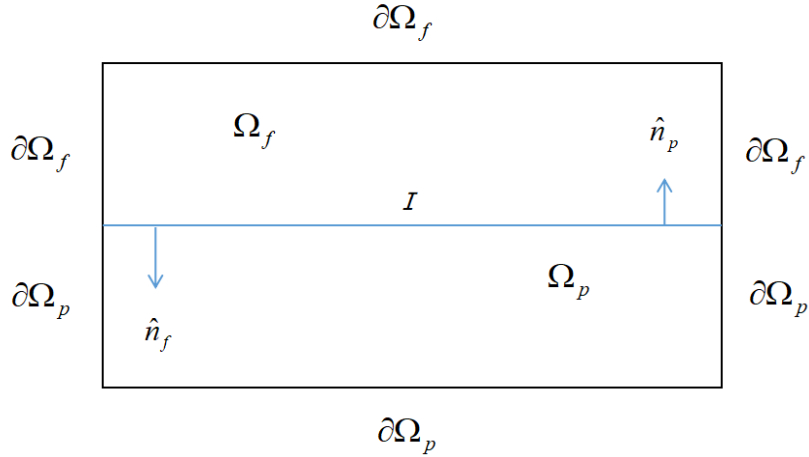}
  \caption{\small A sketch of the free fluid flow domain $\Omega_{f}$, porous media flow domain $\Omega_{p}$ , and the interface $\uppercase\expandafter{\romannumeral1}$.}\label{fig1}
\end{figure}

In $\Omega_{f}$, we assume that the fluid flow with heat transfer is governed by the Navier-Stokes equation and the heat equation:
\begin{equation}\label{f1}
\mathbf{u}_{f,t}-Pr\triangle\mathbf{u}_{f}+(\mathbf{u}_{f}\cdot\nabla)\mathbf{u}_{f}+\nabla p_{f}=PrRa\mathbf{\xi}\theta_{f}+\mathbf{f}_{f} \qquad  in\ \Omega_{f} \times (0,T] \times\Omega,\\
\end{equation}

\begin{equation}\label{f2}
\nabla\cdot\mathbf{u}_{f}=0 \qquad  in\ \Omega_{f} \times (0,T]\times\Omega,\\
\end{equation}

\begin{equation}\label{f3}
\mathbf{u}_{f}=0 \qquad  on\ \partial\Omega_{f}\backslash\uppercase\expandafter{\romannumeral1}\times (0,T],
\end{equation}

\begin{equation}\label{f4}
\mathbf{u}_{f}(0,x)=\mathbf{u}_{f}^{0}(x) \qquad  in\ \Omega_{f},
\end{equation}

\begin{equation}\label{f5}
\theta_{f,t}-k_{f}\triangle\theta_{f}+\mathbf{u}_{f}\cdot\nabla\theta_{f}=\mathbf{\Upsilon}_{f} \qquad in\ \Omega_{f} \times (0,T]\times\Omega,
\end{equation}

\begin{equation}\label{f6}
\theta_{f}=0 \qquad on\ \Gamma_{M}\times(0,T],
\end{equation}

\begin{equation}\label{f7}
\frac{\partial\theta_{f}}{\partial\mathbf{n}_{f}}=0 \qquad on\ \Gamma_{E}\times(0,T],
\end{equation}

\begin{equation}\label{f8}
\theta_{f}(0,x)=\theta_{f}^{0}(x) \qquad in\ \Omega_{f}.
\end{equation}
Here $\mathbf{u}_{f}=(\mathbf{u}_{1},\mathbf{u}_{2},...,\mathbf{u}_{d}), p_{f}$, and $\theta_{f}$ represent the free fluid flow region velocity vector field, pressure, and temperature, respectively. The direction of gravitational acceleration is denoted by the unit vector $\mathbf{\xi}=[0,1]^{T}$. $\mathbf{f}_{f}$ and $\mathbf{\Upsilon}_{f}$ are the external force terms. $Pr$, $Ra$, and $k_{f}$ are the Prandtl number, Rayleigh number, and thermal conductivity, respectively. In the pipe region boundaries, $\Gamma_{M}$ denotes the Dirichlet boundary conditions, and $\Gamma_{E}$ denotes the Neumann boundary condition. In addition, $\Gamma_{M} \bigcup \Gamma_{E}=\partial\Omega_{f}\backslash\uppercase\expandafter{\romannumeral1}$. For convenience, we define $\mathbf{u}_{f,t}=\frac{\partial\mathbf{u}_{f}}{\partial t}$ and $\theta_{f,t}=\frac{\partial\theta_{f}}{\partial t}$.

The porous media flow with heat transfer can be governed by the following Darcy's law coupled with a heat equation:
\begin{equation}\label{p1}
\frac{C_{a}K}{L^{2}}\mathbf{u}_{p,t}+Pr\mathbf{u}_{p}=-\frac{K}{L^{2}}\nabla\phi_{p}+\frac{PrRaK}{L^{2}}\mathbf{\xi}\theta_{p} \qquad in\ \Omega_{p}\times(0,T]\times\Omega,
\end{equation}

\begin{equation}\label{p2}
\nabla\cdot\mathbf{u}_{p}=0 \qquad in\ \Omega_{p}\times(0,T]\times\Omega,
\end{equation}

\begin{equation}\label{p3}
\mathbf{u}_{p}(0,x)=\mathbf{u}_{p}^{0}(x) \qquad in\ \Omega_{p},
\end{equation}

\begin{equation}\label{p4}
\mathbf{u}_{p}\cdot\mathbf{n}_{p}=0 \qquad on\ \partial\Omega_{p}\backslash\uppercase\expandafter{\romannumeral1},
\end{equation}

\begin{equation}\label{p5}
\theta_{p,t}-k_{p}\triangle\theta_{p}+\mathbf{u}_{p}\cdot\nabla\theta_{p}=\mathbf{\Upsilon}_{p} \qquad in \ \Omega_{p}\times(0,T]\times\Omega,
\end{equation}

\begin{equation}\label{p6}
\theta_{p}=0 \qquad on\ \Gamma_{N}\times(0,T],
\end{equation}

\begin{equation}\label{p7}
\frac{\partial\theta_{p}}{\partial\mathbf{n}_{p}}=0 \qquad on\ \Gamma_{Z}\times(0,T],
\end{equation}

\begin{equation}\label{p8}
\theta_{p}(0,x)=\theta_{p}^{0}(x) \qquad in\ \Omega_{p}.
\end{equation}
Here $\mathbf{u}_{p}=(\mathbf{u}_{p1},\mathbf{u}_{p2},...,\mathbf{u}_{pd}),\phi_{p}$, and $\theta_{p}$ denote the porous medium fluid flow region velocity vector field, pressure, and temperature,
respectively. $\mathbf{\Upsilon}_{p}$ is the source term. $L$, $K$, and $k_{p}$ are the characteristic length, hydraulic conductivity tensor, and thermal conductivity, respectively. $C_{a}$ is the dimensionless parameter that represents the acceleration coefficient. $K$ is assumed to be symmetric positive definite. $\Gamma_{N}$ and $\Gamma_{Z}$ denote the Dirichlet and the Neumann boundary conditions, respectively, in the porous media region boundaries where $\Gamma_{N}\bigcup\Gamma_{Z}=\partial\Omega_{p}$. In addition, we define $\mathbf{u}_{p,t}=\frac{\partial\mathbf{u}_{p}}{\partial t}$ and $\theta_{p,t}=\frac{\partial\theta_{p}}{\partial t}$.

We assume that the hydraulic conductivity tensor $K$ is a symmetric positive definite matrix, $K = diag(k, . . . , k)$ with $k > 0$, and it is uniformly bounded in $\Omega_{p}$: there are $k_{min} > 0$ and $k_{max} > 0$ such that
\begin{align}\label{K}
\begin{split}
&k_{min}\mid x\mid^{2} \leq Kx\cdot x \leq k_{max}\mid x\mid^{2} a.e. \  x \in \Omega_{p}.\nonumber
\end{split}
\end{align}

On the interface $\uppercase\expandafter{\romannumeral1}$, we impose the following four interface conditions \cite{P. Hansbo-2005, S. C. Hirata-2007}:

\begin{equation}\label{i1}
\theta_{f}=\theta_{p} \qquad on\ \uppercase\expandafter{\romannumeral1},
\end{equation}

\begin{equation}\label{i2}
\mathbf{n}_{f}\cdot k_{f}\nabla\theta_{f}=-\mathbf{n}_{p}\cdot k_{p}\nabla\theta_{p} \qquad on\ \uppercase\expandafter{\romannumeral1},
\end{equation}

\begin{equation}\label{i3}
\mathbf{u}_{f}\cdot\mathbf{n}_{f}=0, \quad \mathbf{u}_{f}\cdot\mathbf{\tau}=0 \qquad on\ \uppercase\expandafter{\romannumeral1},
\end{equation}

\begin{equation}\label{i4}
\mathbf{u}_{p}\cdot\mathbf{n}_{p}=0 \qquad on\ \uppercase\expandafter{\romannumeral1},
\end{equation}
where $\mathbf{\tau}$ denote the unit tangential vector along the interface.

Herein we consider computing an closed-loop geothermal system corresponding to $J$ different parameter $K_{j}$.
\begin{equation}\label{f11}
\mathbf{u}_{f,j,t}-Pr\triangle\mathbf{u}_{f,j}+(\mathbf{u}_{f,j}\cdot\nabla)\mathbf{u}_{f,j}+\nabla p_{f,j}=PrRa\mathbf{\xi}\theta_{f,j}+\mathbf{f}_{f,j},\ \nabla\cdot\mathbf{u}_{f,j}=0 \qquad  in\ \Omega_{f} \times (0,T],\\
\end{equation}

\begin{equation}\label{f22}
\mathbf{u}_{f,j}=0 \qquad  on\ \partial\Omega_{f}\backslash\uppercase\expandafter{\romannumeral1}\times (0,T],
\end{equation}

\begin{equation}\label{f33}
\theta_{f,j,t}-k_{f}\triangle\theta_{f,j}+\mathbf{u}_{f,j}\cdot\nabla\theta_{f,j}=\mathbf{\Upsilon}_{f,j} \qquad in\ \Omega_{f} \times (0,T],
\end{equation}

\begin{equation}\label{f44}
\theta_{f,j}=0 \quad on\ \Gamma_{M}\times(0,T],\quad and\quad \frac{\partial\theta_{f,j}}{\partial\mathbf{n}_{f}}=0 \quad on\ \Gamma_{E}\times(0,T].
\end{equation}

\begin{equation}\label{p11}
\frac{C_{a}K_{j}}{L^{2}}\mathbf{u}_{p,j,t}+Pr\mathbf{u}_{p,j}=-\frac{K_{j}}{L^{2}}\nabla\phi_{p,j}+\frac{PrRaK_{j}}{L^{2}}\mathbf{\xi}\theta_{p,j},\ \nabla\cdot\mathbf{u}_{p,j}=0 \qquad in\ \Omega_{p}\times(0,T],
\end{equation}

\begin{equation}\label{p22}
\mathbf{u}_{p,j}\cdot\mathbf{n}_{p}=0 \qquad on\ \partial\Omega_{p}\backslash\uppercase\expandafter{\romannumeral1},
\end{equation}

\begin{equation}\label{p33}
\theta_{p,j,t}-k_{p}\triangle\theta_{p,j}+\mathbf{u}_{p,j}\cdot\nabla\theta_{p,j}=\mathbf{\Upsilon}_{p,j} \qquad in \ \Omega_{p}\times(0,T],
\end{equation}

\begin{equation}\label{p44}
\theta_{p,j}=0 \quad on\ \Gamma_{N}\times(0,T],\quad and\quad \frac{\partial\theta_{p,j}}{\partial\mathbf{n}_{p}}=0 \quad on\ \Gamma_{Z}\times(0,T].
\end{equation}
We assume that there is uncertainty in the hydraulic conduction tensor $K$ and that $K_{j}(j=1,...,J)$ is one of the samples drawn from the probability distribution. $J$ is the number of total samples.

Denote by $(\cdot,\cdot)$ the inner product associated with $H^{k}(D)$, and the norm of $H^{k}(D)$ is denoted by $\parallel\cdot\parallel_{H^{k}(D)}$. The Sobolev space $H^{0}(D)$ coincides with the space of square integrable functions $L^{2}(D)$, in which case the norm and inner product are denoted by $\parallel\cdot\parallel_{L^{2}(D)}$ and $(\cdot,\cdot)$, respectively. And $D$ may be $\Omega_{f},\ \Omega_{p}$ or $\uppercase\expandafter{\romannumeral1}$.
Define the function spaces:
\begin{align*}
\begin{split}
&V_{f}:=\big[H_{0}^{1}(\Omega_{f})\big]^{d}:=\Big\{\mathbf{v}_{f}\in\big[H^{1}(\Omega_{f})\big]^{d}:\mathbf{v}_{f}=0 \quad on\ \partial\Omega_{f}\Big\},\\
&V_{p}:=H(div;\Omega_{p}):=\Big\{\mathbf{v}_{p}\in \big[L^{2}(\Omega_{p})\big]^{d},\nabla\cdot\mathbf{v}_{p}\in L^{2}(\Omega_{p}):\mathbf{v}_{p}\cdot\mathbf{n}_{p}=0 \quad on \ \partial\Omega_{p}\Big\},\\
&W_{f}:=H_{0}^{1}(\Omega_{f}):=\Big\{\varphi\in H^{1}(\Omega_{f}):\varphi=0 \quad on \ \Gamma_{M}\backslash\uppercase\expandafter{\romannumeral1}\Big\},\\
&W_{p}:=H_{0}^{1}(\Omega_{p}):=\Big\{\omega\in H^{1}(\Omega_{p}):\omega=0\quad on \ \Gamma_{N}\backslash\uppercase\expandafter{\romannumeral1}\Big\},\\
&U_{f}:=L_{0}^{2}(\Omega_{f}):=\bigg\{q\in L^{2}(\Omega_{f}):\int_{\Omega_{f}}q dx=0\bigg\},\\
&U_{p}:=L_{0}^{2}(\Omega_{p}):=\bigg\{\psi\in L^{2}(\Omega_{p}):\int_{\Omega_{p}}\psi dx=0\bigg\}.
\end{split}
\end{align*}

Set the solenoidal spaces
\begin{align*}
\begin{split}
&S_{f}:=\big\{\mathbf{v}_{f}\in V_{f}:\nabla\cdot\mathbf{v}_{f}=0\big\}, \quad S_{p}:=\big\{\mathbf{v}_{p}\in V_{p}:\nabla\cdot\mathbf{v}_{p}=0\big\}.
\end{split}
\end{align*}

We also define a product space
\begin{align*}
\begin{split}
&Y_{T}=W_{f}\times W_{p},
\end{split}
\end{align*}
and
\begin{align*}
\begin{split}
&Y_{\uppercase\expandafter{\romannumeral1}}:=\big\{(\varphi,\omega)\in(W_{f}\times W_{p}):\varphi\mid _{\uppercase\expandafter{\romannumeral1}}=\omega\mid _{\uppercase\expandafter{\romannumeral1}}\big\}\subset Y_{T}.
\end{split}
\end{align*}

Then we introduce the trilinear form
\begin{equation}
c(\mathbf{u}_{f},\mathbf{v}_{f},\mathbf{w})_{\Omega_{f}}
=\frac{1}{2}((\mathbf{u}_{f}\cdot\nabla)\mathbf{v}_{f},\mathbf{w})-\frac{1}{2}((\mathbf{u}_{f}\cdot\nabla)\mathbf{w},\mathbf{v}_{f})\quad \forall\  \mathbf{u}_{f},\mathbf{v}_{f},\mathbf{w}\in V_{f}.
\end{equation}

Similarly, we can define another two trilinear forms for any $(\varphi,\omega)\in Y_{\uppercase\expandafter{\romannumeral1}}$ and $(\theta_{f},\theta_{p})\in Y_{T}$:
\begin{align}
\begin{split}
&t_{f}(\mathbf{u}_{f},\theta_{f},\varphi)_{\Omega_{f}}=\frac{1}{2}((\mathbf{u}_{f}\cdot\nabla)\theta_{f},\varphi)-\frac{1}{2}((\mathbf{u}_{f}\cdot\nabla)\varphi,\theta_{f}) \quad \forall\ \mathbf{u}_{f}\in V_{f},\\
&t_{p}(\mathbf{u}_{p},\theta_{p},\omega)_{\Omega_{p}}=\frac{1}{2}((\mathbf{u}_{p}\cdot\nabla)\theta_{p},\omega)-\frac{1}{2}((\mathbf{u}_{p}\cdot\nabla)\omega,\theta_{p}) \quad \forall\ \mathbf{u}_{p}\in V_{p}.
\end{split}
\end{align}

In addition, if $\mathbf{u}_{f}\in S_{f}$ and $\mathbf{u}_{p}\in S_{p}$, then for any $(\mathbf{v}_{f},\mathbf{w})\in V_{f}$, $(\theta_{f},\theta_{p})\in Y_{T}$, and $(\varphi,\omega)\in Y_{\uppercase\expandafter{\romannumeral1}}$, we have \cite{J. M. Connors-2011}
\begin{align*}
\begin{split}
&c(\mathbf{u}_{f},\mathbf{v}_{f},\mathbf{w})_{\Omega_{f}}=((\mathbf{u}_{f}\cdot\nabla)\mathbf{v}_{f},\mathbf{w})_{\Omega_{f}},\quad
t_{f}(\mathbf{u}_{f},\theta_{f},\varphi)_{\Omega_{f}}=((\mathbf{u}_{f}\cdot\nabla)\theta_{f},\varphi)_{\Omega_{f}},\\
&t_{p}(\mathbf{u}_{p},\theta_{p},\omega)_{\Omega_{p}}=((\mathbf{u}_{p}\cdot\nabla)\theta_{p},\omega)_{\Omega_{p}}.
\end{split}
\end{align*}

\begin{lemma}\label{lemma2.1}
(\cite{R. Temam-1984}). When $D$ is of class $C^{l}$ and regular enough, the trilinear form $((u\cdot\nabla)v,w)_{D}$ ,satisfies the following continuity property:
\begin{align}
\begin{split}
&\mid((u\cdot\nabla)v,w)_{D}\mid \leq N\parallel u\parallel_{s_{1}}\parallel v\parallel_{s_{2}+1}\parallel w\parallel_{s_{3}},\qquad
\forall\ u\in H^{s_{1}}(D),\ v\in H^{s_{2}+1}(D), w\in H^{s_{3}}(D),
\end{split}
\end{align}
where $s_{1} + s_{2} + s_{3} \geq 1$, $s_{i} \neq 1, i = 1, 2, 3, 0 \leq s_{1} \leq l, 0 \leq s_{3} \leq l, 0 \leq s_{2} \leq l - 1$. If $s_{i} = 1$ holds for some or other $i$, then we need $s_{1} + s_{2} + s_{3} > 1$. $N$ is a positive constant depending only on domain $D$, where $D$ may be $\Omega_{f}$ or $\Omega_{p}$.
\end{lemma}

The weak formulations of the coupled system (\ref{f11}) - (\ref{p44}) as follows: find $(\mathbf{u}_{f,j},\mathbf{u}_{p,j}; p_{f,j},\phi_{p,j}; \theta_{f,j},\theta_{p,j})\in (V_{f}\times V_{p} \times U_{f}\times U_{p}\times Y_{T})$ such that for any $(\mathbf{v}_{f},\mathbf{v}_{p}; q,\psi; \varphi,\omega)\in (V_{f}\times V_{p} \times U_{f}\times U_{p}\times Y_{\uppercase\expandafter{\romannumeral1}})$,
\begin{align}\label{RUOXINGSHI1}
\begin{split}
&(\mathbf{u}_{f,j,t},\mathbf{v}_{f})_{\Omega_{f}}
+\frac{C_{a}K_{j}}{L^{2}}(\mathbf{u}_{p,j,t},\mathbf{v}_{p})_{\Omega_{p}}
-(p_{f,j},\nabla\cdot\mathbf{v}_{f})_{\Omega_{f}}+Pr(\nabla\mathbf{u}_{f,j},\nabla\mathbf{v}_{f})_{\Omega_{f}}\\
&+c(\mathbf{u}_{f,j},\mathbf{u}_{f,j},\mathbf{v}_{f})_{\Omega_{f}}+Pr(\mathbf{u}_{p,j},\mathbf{v}_{p})_{\Omega_{p}}
-\frac{K_{j}}{L^{2}}(\phi_{p,j},\nabla\cdot\mathbf{v}_{p})_{\Omega_{p}}\\
=&PrRa(\mathbf{\xi}\theta_{f,j},\mathbf{v}_{f})_{\Omega_{f}}+\frac{PrRa}{L^{2}}(\mathbf{\xi}K_{j}\theta_{p,j},\mathbf{v}_{p})_{\Omega_{p}}
+(\mathbf{f}_{f,j},\mathbf{v}_{f})_{\Omega_{f}},
\end{split}
\end{align}

\begin{align}\label{RUOXINGSHI2}
\begin{split}
&(q,\nabla\cdot\mathbf{u}_{f,j})_{\Omega_{f}}=0,
\end{split}
\end{align}

\begin{align}\label{RUOXINGSHI3}
\begin{split}
&\frac{K_{j}}{L^{2}}(\psi,\nabla\cdot\mathbf{u}_{p,j})_{\Omega_{p}}=0,
\end{split}
\end{align}

\begin{align}\label{RUOXINGSHI4}
\begin{split}
&(\theta_{f,j,t},\varphi)_{\Omega_{f}}
+(\theta_{p,j,t},\omega)_{\Omega_{p}}
+k_{f}(\nabla\theta_{f,j},\nabla\varphi)_{\Omega_{f}}+k_{p}(\nabla\theta_{p,j},\nabla\omega)_{\Omega_{p}}\\
&+t_{f}(\mathbf{u}_{f,j},\theta_{f,j},\varphi)_{\Omega_{f}}+t_{p}(\mathbf{u}_{p,j},\theta_{p,j},\omega)_{\Omega_{p}}
-k_{f}\int_{\uppercase\expandafter{\romannumeral1}}\mathbf{n}_{f}\cdot\nabla\theta_{f,j}(\varphi-\omega)dl\\
&+\frac{k_{f}\gamma}{h}\int_{\uppercase\expandafter{\romannumeral1}}(\theta_{f,j}-\theta_{p,j})(\varphi-\omega)dl
=(\mathbf{\Upsilon}_{f,j},\varphi)_{\Omega_{f}}+(\mathbf{\Upsilon}_{p,j},\omega)_{\Omega_{p}}.
\end{split}
\end{align}

The proof of the existence and uniqueness of the solution to (\ref{RUOXINGSHI1}) - (\ref{RUOXINGSHI4}) is similar to \cite{M A A Mahbub-2020} and will not be elaborated here.
\section{Numerical algorithms}
~\par To discretize the stochastic closed-loop geothermal system problem in space by the finite element method, we choose conforming velocity, pressure, temperature finite element spaces based on an edge to edge triangulation of the domain $D$ with maximum element diameter $h$:
\begin{align*}
\begin{split}
&V_{f}^{h}\subset V_{f},\quad U_{f}^{h}\subset U_{f},\quad W_{f}^{h}\subset W_{f},\\
&V_{p}^{h}\subset V_{p},\quad U_{p}^{h}\subset U_{p},\quad W_{p}^{h}\subset W_{p}.
\end{split}
\end{align*}

The computation procedure is as follows.

(1) Generate a number of independently, identically distributed samples for the random hydraulic conductivity $K_{j},j=1,...,J$;

(2) The fully discrete approximation of (\ref{f11}) - (\ref{p44}) is:

Find $(\mathbf{u}_{f,j}^{h},\mathbf{u}_{p,j}^{h};\ p_{f,j}^{h},\phi_{p,j}^{h};\ \theta_{f,j}^{h},\theta_{p,j}^{h})\in (V_{f}^{h},V_{p}^{h};\ U_{f}^{h},U_{p}^{h};\ W_{f}^{h},W_{p}^{h})$ satisfying $\forall (\mathbf{v}_{f}^{h},\mathbf{v}_{p}^{h};\ q^{h},$
$\psi^{h};\ \varphi^{h},\omega^{h})\in (V_{f}^{h},V_{p}^{h};\ U_{f}^{h},U_{p}^{h};\ W_{f}^{h},W_{p}^{h})$,

$\bullet$Step 1
\begin{align}\label{suanfa 1}
\begin{split}
&\Bigg(\frac{\mathbf{u}_{f,j}^{h,n+1}-\mathbf{u}_{f,j}^{h,n}}{\Delta t},\mathbf{v}_{f}^{h}\Bigg)_{\Omega_{f}}
-(p_{f,j}^{h,n+1},\nabla\cdot\mathbf{v}_{f}^{h})_{\Omega_{f}}
+Pr(\nabla \mathbf{u}_{f,j}^{h,n+1},\nabla\mathbf{v}_{f}^{h})_{\Omega_{f}}
+c(\mathbf{u}_{f,j}^{h,n},\mathbf{u}_{f,j}^{h,n+1},\mathbf{v}_{f}^{h})_{\Omega_{f}}\\
=&PrRa(\mathbf{\xi}\theta_{f,j}^{h,n},\mathbf{v}_{f}^{h})_{\Omega_{f}}
+(\mathbf{f}_{f,j}(t_{n+1}),\mathbf{v}_{f}^{h})_{\Omega_{f}},
\end{split}
\end{align}

\begin{align}\label{suanfa 2}
\begin{split}
&(q^{h},\nabla\cdot\mathbf{u}_{f,j}^{h,n+1})_{\Omega_{f}}=0,
\end{split}
\end{align}

$\bullet$Step 2
\begin{align}\label{suanfa 3}
\begin{split}
&\Bigg(\frac{\theta_{f,j}^{h,n+1}-\theta_{f,j}^{h,n}}{\Delta t},\varphi^{h}\Bigg)_{\Omega_{f}}
+k_{f}(\nabla\theta_{f,j}^{h,n+1},\nabla\varphi^{h})_{\Omega_{f}}
+t_{f}(\mathbf{u}_{f,j}^{h,n},\theta_{f,j}^{h,n+1},\varphi^{h})_{\Omega_{f}}\\
&-k_{f}\int_{\uppercase\expandafter{\romannumeral1}}\mathbf{n}_{f}\cdot\nabla\theta_{f,j}^{h,n}\varphi^{h}dl
+\frac{k_{f}\gamma}{h}\int_{\uppercase\expandafter{\romannumeral1}}(\theta_{f,j}^{h,n+1}-\theta_{p,j}^{h,n})\varphi^{h}dl
=(\mathbf{\Upsilon}_{f,j}(t_{n+1}),\varphi^{h})_{\Omega_{f}},
\end{split}
\end{align}

$\bullet$Step 3
\begin{align}\label{suanfa 4}
\begin{split}
&\frac{C_{a}}{L^{2}}\Bigg(K_{j}\frac{\mathbf{u}_{p,j}^{h,n+1}-\mathbf{u}_{p,j}^{h,n}}{\Delta t},\mathbf{v}_{p}^{h}\Bigg)_{\Omega_{p}}
+Pr(\mathbf{u}_{p,j}^{h,n+1},\mathbf{v}_{p}^{h})_{\Omega_{p}}-\frac{1}{L^{2}}(K_{j}\phi_{p,j}^{h,n+1},\nabla\cdot\mathbf{v}_{p}^{h})_{\Omega_{p}}\\
=&\frac{PrRa}{L^{2}}(\mathbf{\xi}K_{j}\theta_{p,j}^{h,n},\mathbf{v}_{p}^{h})_{\Omega_{p}},
\end{split}
\end{align}

\begin{equation}\label{suanfa 5}
\frac{1}{L^{2}}(K_{j}\psi^{h},\nabla\cdot\mathbf{u}_{p,j}^{h,n+1})_{\Omega_{p}}=0,
\end{equation}

$\bullet$Step 4
\begin{align}\label{suanfa 6}
\begin{split}
&\Bigg(\frac{\theta_{p,j}^{h,n+1}-\theta_{p,j}^{h,n}}{\triangle t},\omega^{h}\Bigg)_{\Omega_{p}}
+k_{p}(\nabla\theta_{p,j}^{h,n+1},\nabla\omega^{h})_{\Omega_{p}}
+t_{p}(\mathbf{u}_{p,j}^{h,n},\theta_{p,j}^{h,n+1},\omega^{h})_{\Omega_{p}}\\
&+k_{f}\int_{\uppercase\expandafter{\romannumeral1}}\mathbf{n}_{f}\cdot\nabla\theta_{f,j}^{h,n}\omega^{h}dl
-\frac{k_{f}\gamma}{h}\int_{\uppercase\expandafter{\romannumeral1}}(\theta_{f,j}^{h,n+1}-\theta_{p,j}^{h,n+1})\omega^{h}dl
=(\mathbf{\Upsilon}_{p,j}(t_{n+1}),\omega^{h})_{\Omega_{p}}.
\end{split}
\end{align}

(3) Output required statistical information such as the
expectation of $v$: $E[v]\approx \frac{1}{J}\sum_{j=1}^{J}v_{j}^{h}$, where $v=\vec{u}_{f}^{h},\vec{u}_{p}^{h},\theta_{f}^{h},\theta_{p}^{h}$.

Now, we will introduce some inequalities that are needed for analysis.
\begin{lemma}
(Discrete Grnowall's lemma) Let $C$ and $a_{k}, b_{k}, c_{k}, d_{k}$ for integer $k\geq 0$ be non-negative number such that \cite{S. Li-2009}
\begin{align}
\begin{split}
&a_{n}+\Delta t\sum\limits_{k=0}^{n}b_{k}\leq \Delta t\sum\limits_{k=0}^{n}d_{k}a_{k}+\Delta t\sum\limits_{k=0}^{n}c_{k}+C,\qquad \forall\ n\geq 1.
\end{split}
\end{align}
Suppose that $\Delta td_{k}\leq 1$ for all $k$, then we have
\begin{align}
\begin{split}
&a_{n}+\Delta t\sum\limits_{k=0}^{n}b_{k}\leq exp(\Delta t\sum\limits_{k=0}^{n}d_{k})(\Delta t\sum\limits_{k=0}^{n}c_{k}+C),\qquad \forall\ n\geq 1.
\end{split}
\end{align}
\end{lemma}

Furthermore, for the sake of theoretical analysis, we will review the following Poincar$\acute{e}$ and trace inequalities: there exist
constant $C_{f}, C_{t},C_{T}$ that only depend on the domain $\Omega_{f}$, and $C_{p}, \tilde{C_{t}},\tilde{C}_{T}$ that only depend on the domain $\Omega_{p}$, such that for all $\mathbf{v}_{f}\in V_{f}, \mathbf{v}_{p}\in V_{p}, \varphi\in W_{f}, \omega\in W_{p}$,
\begin{align}
\begin{split}
&\parallel\mathbf{v}_{f}\parallel_{L^{2}(\Omega_{f})}\leq C_{f}\parallel\nabla\mathbf{v}_{f}\parallel_{L^{2}(\Omega_{f})},\qquad \qquad\quad \parallel\mathbf{v}_{p}\parallel_{L^{2}(\Omega_{p})}\leq C_{p}\parallel\nabla\mathbf{v}_{p}\parallel_{L^{2}(\Omega_{p})},\\
&\parallel\varphi\parallel_{L^{2}(\Omega_{f})}\leq C_{t}\parallel\nabla\varphi\parallel_{L^{2}(\Omega_{f})},\qquad\qquad\qquad
\parallel\omega\parallel_{L^{2}(\Omega_{p})}\leq \tilde{C}_{t}\parallel\nabla\omega\parallel_{L^{2}(\Omega_{p})},\\
&\parallel\varphi\parallel_{L^{2}(\uppercase\expandafter{\romannumeral1})}\leq C_{T}^{\frac{1}{2}}\parallel\varphi\parallel^{\frac{1}{2}}_{L^{2}(\Omega_{f})}\parallel\nabla\varphi\parallel^{\frac{1}{2}}_{L^{2}(\Omega_{f})},\qquad    \parallel\omega\parallel_{L^{2}(\uppercase\expandafter{\romannumeral1})}\leq \tilde{C}_{T}^{\frac{1}{2}}\parallel\omega\parallel^{\frac{1}{2}}_{L^{2}(\Omega_{p})}\parallel\nabla\omega\parallel^{\frac{1}{2}}_{L^{2}(\Omega_{p})}.
\end{split}
\end{align}

The following inverse inequality will be used in the proofs, \cite{P.G. Ciarlet-1978, Z. Chen-2005}.
\begin{align*}
\begin{split}
&\parallel\nabla \theta_{f}^{h}\parallel_{L^{2}(\Omega_{f})}\leq C_{inv}h^{-1}\parallel\theta_{f}^{h}\parallel_{L^{2}(\Omega_{f})}\quad \forall\ \theta_{f}^{h}\in W_{f}^{h},\\
&\parallel\nabla \theta_{p}^{h}\parallel_{L^{2}(\Omega_{p})}\leq \tilde{C}_{inv}h^{-1}\parallel\theta_{p}^{h}\parallel_{L^{2}(\Omega_{p})}\quad \forall\ \theta_{p}^{h}\in W_{p}^{h},
\end{split}
\end{align*}
where $C_{inv}$ and $\tilde{C}_{inv}$ are constants.

We define the linear projection operators: for all $(\mathbf{v}_{f}^{h},\mathbf{v}_{p}^{h};q^{h},\psi^{h};\varphi^{h},\omega^{h})
\in(V_{f}^{h},V_{p}^{h};U_{f}^{h},U_{p}^{h};W_{f}^{h},W_{p}^{h})$ and $t\in (0,T]$, $\Pi_{f}^{h}: \mathbf{u}_{f}(t)\in V_{f}\longrightarrow \Pi_{f}^{h}\mathbf{u}_{f}(t)\in V_{f}^{h},\ \Pi_{p}^{h}: \mathbf{u}_{p}(t)\in V_{p}\longrightarrow \Pi_{p}^{h}\mathbf{u}_{p}(t)\in V_{p}^{h},\
\rho_{f}^{h}: p_{f}(t)\in U_{f}\longrightarrow \rho_{f}^{h}p_{f}(t)\in U_{f}^{h},\ \rho_{p}^{h}: \phi_{p}(t)\in U_{p}\longrightarrow \rho_{p}^{h}\phi_{p}(t)\in U_{p}^{h},\ Q_{f}^{h}: \theta_{f}(t)\in W_{f}\longrightarrow Q_{f}^{h}\theta_{f}(t)\in W_{f}^{h}$, and $\ Q_{p}^{h}: \theta_{p}(t)\in W_{p}\longrightarrow Q_{p}^{h}\theta_{p}(t)\in W_{p}^{h}$ satisfies

\begin{align}\label{touyin1}
\begin{split}
&Pr(\nabla \mathbf{u}_{f}(t),\nabla\mathbf{v}_{f}^{h})_{\Omega_{f}}+Pr(\mathbf{u}_{p}(t),\mathbf{v}_{p}^{h})_{\Omega_{p}}
+k_{f}(\nabla \theta_{f}(t),\nabla\varphi^{h})_{\Omega_{f}}+k_{p}(\nabla \theta_{p}(t),\nabla\omega^{h})_{\Omega_{p}}\\
&+(p_{f}(t),\nabla\cdot\mathbf{v}_{f}^{h})_{\Omega_{f}}+\frac{K}{L^{2}}(\phi_{p}(t),\nabla\cdot\mathbf{v}_{p}^{h})_{\Omega_{p}}
-k_{f}\int_{\uppercase\expandafter{\romannumeral1}}\mathbf{n}_{f}\cdot\nabla\theta_{f}(\varphi^{h}-\omega^{h})ds\\
=&Pr(\nabla \Pi_{f}^{h}\mathbf{u}_{f}(t),\nabla\mathbf{v}_{f}^{h})_{\Omega_{f}}+Pr(\Pi_{p}^{h}\mathbf{u}_{p}(t),\mathbf{v}_{p}^{h})_{\Omega_{p}}
+k_{f}(\nabla Q_{f}^{h}\theta_{f}(t),\nabla\varphi^{h})_{\Omega_{f}}+k_{p}(\nabla Q_{p}^{h}\theta_{p}(t),\nabla\omega^{h})_{\Omega_{p}}\\
&+(\rho_{f}^{h}p_{f}(t),\nabla\cdot\mathbf{v}_{f}^{h})_{\Omega_{f}}+\frac{K}{L^{2}}(\rho_{p}^{h}\phi_{p}(t),\nabla\cdot\mathbf{v}_{p}^{h})_{\Omega_{p}}
-k_{f}\int_{\uppercase\expandafter{\romannumeral1}}\mathbf{n}_{f}\cdot\nabla Q_{f}^{h}\theta_{f}(\varphi^{h}-\omega^{h})ds,\\
&\\
&(q^{h},\nabla\cdot\Pi_{f}^{h}\mathbf{u}_{f}(t))_{\Omega_{f}}=0,\\
&\\
&\frac{K}{L^{2}}(\psi^{h},\nabla\cdot\Pi_{f}^{h}\mathbf{u}_{p}(t))_{\Omega_{p}}=0.\\
\end{split}
\end{align}

Besides we assume that $(\mathbf{u}_{f}(t),\mathbf{u}_{p}(t); p_{f}(t),\phi_{p}(t); \theta_{f}(t),\theta_{p}(t))$ is smooth enough and the projection $(\Pi_{f}^{h},\Pi_{p}^{h}; $
$\rho_{f}^{h},\rho_{p}^{h}; Q_{f}^{h}, Q_{p}^{h})$ of $(\mathbf{u}_{f}(t),\mathbf{u}_{p}(t); p_{f}(t),\phi_{p}(t); \theta_{f}(t),\theta_{p}(t))$ satisfies following approximation properties:
\begin{align}\label{touyin2}
\begin{split}
&\parallel \mathbf{u}_{f}(t)-\Pi_{f}^{h}\mathbf{u}_{f}(t)\parallel_{L^{2}(\Omega_{f})}
+h\parallel \nabla(\mathbf{u}_{f}(t)-\Pi_{f}^{h}\mathbf{u}_{f}(t))\parallel_{L^{2}(\Omega_{f})}
\leq Ch^{k+1}\parallel\mathbf{u}_{f}(t)\parallel_{H^{k+1}(\Omega_{f})},\\
&\parallel \mathbf{u}_{p}(t)-\Pi_{p}^{h}\mathbf{u}_{p}(t)\parallel_{L^{2}(\Omega_{p})}
\leq Ch^{k+1}\parallel\mathbf{u}_{p}(t)\parallel_{H^{k+1}(\Omega_{p})},\\
&\parallel \theta_{f}(t)-Q_{f}^{h}\theta_{f}(t)\parallel_{L^{2}(\Omega_{f})}
+h\parallel \nabla(\theta_{f}(t)-Q_{f}^{h}\theta_{f}(t))\parallel_{L^{2}(\Omega_{f})}
\leq Ch^{k+1}\parallel\theta_{f}(t)\parallel_{H^{k+1}(\Omega_{f})},\\
&\parallel \theta_{p}(t)-Q_{p}^{h}\theta_{p}(t)\parallel_{L^{2}(\Omega_{p})}
+h\parallel \nabla(\theta_{p}(t)-Q_{p}^{h}\theta_{p}(t))\parallel_{L^{2}(\Omega_{p})}
\leq Ch^{k+1}\parallel\theta_{p}(t)\parallel_{H^{k+1}(\Omega_{p})},
\end{split}
\end{align}
where the generic constant $C>0$ is independent of the mesh size $h$.
\section{Theoretical analysis}
~\par In this section, we will analyze the stability and convergence of the stochastic closed-loop geothermal system. In addition, $C$ in the analysis below represents different positive constants in different cases.
\subsection{Stability analysis}
~\par The following theorem illustrates the stability of the decoupled finite element method for the stochastic closed-loop geothermal system.

\begin{theorem}\label{wending4-1}
Let $\{\mathbf{u}_{f}^{h,n},p_{f}^{h,n},\theta_{f}^{h,n};1\leq n\leq N\}$, $\{\mathbf{u}_{p}^{h,n},\phi_{p}^{h,n},\theta_{p}^{h,n};1\leq n\leq N\}$ be the solution of the problem (\ref{suanfa 1}) - (\ref{suanfa 6}), when penalty parameter satisfies the condition $\gamma>2C_{T}C_{inv}$, where $C_{T}, C_{inv}>0$, we have
\begin{align}
&\max\limits_{1\leq n\leq N}E\big[\parallel\mathbf{u}_{f}^{h,n}\parallel_{L^{2}(\Omega_{f})}^{2}\big]
+E\Big[\sum\limits_{n=1}^{N}\parallel\mathbf{u}_{f}^{h,n}-\mathbf{u}_{f}^{h,n-1}\parallel_{L^{2}(\Omega_{f})}^{2}\Big]
+E\Big[\Delta t\sum\limits_{n=1}^{N}\parallel\nabla\mathbf{u}_{f}^{h,n}\parallel_{L^{2}(\Omega_{f})}^{2}\Big]\nonumber\\
&+\max\limits_{1\leq n\leq N}E\big[\parallel\mathbf{u}_{p}^{h,n}\parallel_{L^{2}(\Omega_{p})}^{2}\big]
+E\Big[\sum\limits_{n=1}^{N}\parallel\mathbf{u}_{p}^{h,n}-\mathbf{u}_{p}^{h,n-1}\parallel_{L^{2}(\Omega_{p})}^{2}\Big]
+E\Big[\Delta t\sum\limits_{n=1}^{N}\parallel\mathbf{u}_{p}^{h,n}\parallel_{L^{2}(\Omega_{p})}^{2}\Big]\nonumber
\end{align}
\begin{align}
&+\max\limits_{1\leq n\leq N}E\big[\parallel\theta_{f}^{h,n}\parallel_{L^{2}(\Omega_{f})}^{2}\big]
+E\Big[\sum\limits_{n=1}^{N}\parallel\theta_{f}^{h,n}-\theta_{f}^{h,n-1}\parallel_{L^{2}(\Omega_{f})}^{2}\Big]
+E\Big[\Delta t\sum\limits_{n=1}^{N}\parallel\nabla\theta_{f}^{h,n}\parallel_{L^{2}(\Omega_{f})}^{2}\Big]\nonumber\\
&+\max\limits_{1\leq n\leq N}E\big[\parallel\theta_{p}^{h,n}\parallel_{L^{2}(\Omega_{p})}^{2}\big]
+E\Big[\sum\limits_{n=1}^{N}\parallel\theta_{p}^{h,n}-\theta_{p,j}^{h,n-1}\parallel_{L^{2}(\Omega_{p})}^{2}\Big]
+E\Big[\Delta t\sum\limits_{n=1}^{N}\parallel\nabla\theta_{p}^{h,n}\parallel_{L^{2}(\Omega_{p})}^{2}\Big]\nonumber\\
\leq& C\Bigg(E\Bigg[\sum\limits_{n=1}^{N}\parallel\mathbf{f}_{f}(t_{n})\parallel_{H^{-1}(\Omega_{f})}^{2}\Big]
+E\Big[\sum\limits_{n=1}^{N}\parallel\mathbf{\Upsilon}_{f}(t_{n})\parallel_{H^{-1}(\Omega_{f})}^{2}\Big]
+E\Big[\sum\limits_{n=1}^{N}\parallel\mathbf{\Upsilon}_{p}(t_{n})\parallel_{H^{-1}(\Omega_{p})}^{2}\Big]\nonumber\\
&+E\big[\parallel\mathbf{u}_{f}^{h,0}\parallel_{L^{2}(\Omega_{f})}^{2}\big]
+E\big[\parallel\mathbf{u}_{p}^{h,0}\parallel_{L^{2}(\Omega_{p})}^{2}\big]
+E\big[\parallel\theta_{f}^{h,0}\parallel_{L^{2}(\Omega_{f})}^{2}\big]
+E\big[\parallel\theta_{p}^{h,0}\parallel_{L^{2}(\Omega_{p})}^{2}\big]\Bigg),
\end{align}
\end{theorem}
where $E[\cdot]$ stands for the expectation. The above theorem is similar to the literature \cite{M A A Mahbub-2020}, and will not be proved in detail here for convenience.
\subsection{Convergence analysis}
~\par In this section, we derive error estimates for the finite element method for the stochastic closed-loop geothermal system. We assume that the true solutions of the closed-loop geothermal system model have the following regularity.
\begin{align*}
\begin{split}
&\mathbf{u}_{f}\in L^{\infty}(0,T;H^{k+1}(\Omega_{f})),\ \mathbf{u}_{f,t}\in L^{2}(0,T;H^{k+1}(\Omega_{f})),\
\mathbf{u}_{f,tt}\in L^{2}(0,T;L^{2}(\Omega_{f})),\\
&\theta_{f}\in L^{\infty}(0,T;H^{k+1}(\Omega_{f})),\ \theta_{f,t}\in L^{2}(0,T;H^{k+1}(\Omega_{f})),\
\theta_{f,tt}\in L^{2}(0,T;L^{2}(\Omega_{f})),\\
&\mathbf{u}_{p}\in L^{\infty}(0,T;H^{k+1}(\Omega_{p})),\ \mathbf{u}_{p,t}\in L^{2}(0,T;H^{k+1}(\Omega_{p})),\
\mathbf{u}_{p,tt}\in L^{2}(0,T;L^{2}(\Omega_{p})),\\
&\theta_{p}\in L^{\infty}(0,T;H^{k+1}(\Omega_{p})),\ \theta_{p,t}\in L^{2}(0,T;H^{k+1}(\Omega_{p})),\
\theta_{p,tt}\in L^{2}(0,T;L^{2}(\Omega_{p})).
\end{split}
\end{align*}

\begin{theorem}\label{wucha4-1}
Let $\{\mathbf{u}_{f}^{h,n},p_{f}^{h,n},\theta_{f}^{h,n};1\leq n\leq N\}$, $\{\mathbf{u}_{p}^{h,n},\phi_{p}^{h,n},\theta_{p}^{h,n};1\leq n\leq N\}$ be the solution of the problem (\ref{suanfa 1}) - (\ref{suanfa 6}), when penalty parameter satisfies the condition $\gamma>8C_{T}C_{inv}$, where $C_{T}, C_{inv}>0$, we have
\begin{align}\label{wucha}
\begin{split}
&\max\limits_{1\leq n\leq N}\bigg(E\big[\parallel \mathbf{u}_{f}(t_{n})-\mathbf{u}_{f}^{h,n}\parallel_{L^{2}(\Omega_{f})}^{2}\big]\bigg)^{\frac{1}{2}}
+\Bigg(E\Big[\sum\limits_{n=1}^{N}\parallel(\mathbf{u}_{f}(t_{n})-\mathbf{u}_{f}^{h,n})
-(\mathbf{u}_{f}(t_{n-1})-\mathbf{u}_{f}^{h,n-1})\parallel_{L^{2}(\Omega_{f})}^{2}\Big]\Bigg)^{\frac{1}{2}}\\
&+\max\limits_{1\leq n\leq N}\bigg(E\big[\parallel \mathbf{u}_{p}(t_{n})-\mathbf{u}_{p}^{h,n}\parallel_{L^{2}(\Omega_{p})}^{2}\big]\bigg)^{\frac{1}{2}}
+\Bigg(E\Big[\sum\limits_{n=1}^{N}\parallel(\mathbf{u}_{p}(t_{n})-\mathbf{u}_{p}^{h,n})
-(\mathbf{u}_{p}(t_{n-1})-\mathbf{u}_{p}^{h,n-1})\parallel_{L^{2}(\Omega_{p})}^{2}\Big]\Bigg)^{\frac{1}{2}}\\
&+\max\limits_{1\leq n\leq N}\bigg(E\big[\parallel \theta_{f}(t_{n})-\theta_{f}^{h,n}\parallel_{L^{2}(\Omega_{f})}^{2}\big]\bigg)^{\frac{1}{2}}
+\Bigg(E\Big[\sum\limits_{n=1}^{N}\parallel(\theta_{f}(t_{n})-\theta_{f}^{h,n})
-(\theta_{f}(t_{n-1})-\theta_{f}^{h,n-1})\parallel_{L^{2}(\Omega_{f})}^{2}\Big]\Bigg)^{\frac{1}{2}}\\
&+\max\limits_{1\leq n\leq N}\bigg(E\big[\parallel \theta_{p}(t_{n})-\theta_{p}^{h,n}\parallel_{L^{2}(\Omega_{p})}^{2}\big]\bigg)^{\frac{1}{2}}
+\Bigg(E\Big[\sum\limits_{n=1}^{N}\parallel(\theta_{p}(t_{n})-\theta_{p}^{h,n})
-(\theta_{p}(t_{n-1})-\theta_{p}^{h,n-1})\parallel_{L^{2}(\Omega_{p})}^{2}\Big]\Bigg)^{\frac{1}{2}}\\
&+\Bigg(E\Big[\Delta t\sum\limits_{n=1}^{N}\parallel\nabla(\mathbf{u}_{f}(t_{n})-\mathbf{u}_{f}^{h,n})\parallel_{L^{2}(\Omega_{f})}^{2}\Big]\Bigg)^{\frac{1}{2}}
+\Bigg(E\Big[\Delta t\sum\limits_{n=1}^{N}\parallel\mathbf{u}_{p}(t_{n})-\mathbf{u}_{p}^{h,n}\parallel_{L^{2}(\Omega_{p})}^{2}\Big]\Bigg)^{\frac{1}{2}}\\
&+\Bigg(E\Big[\Delta t\sum\limits_{n=1}^{N}\parallel\nabla(\theta_{f}(t_{n})-\theta_{f}^{h,n})\parallel_{L^{2}(\Omega_{f})}^{2}\Big]\Bigg)^{\frac{1}{2}}
+\Bigg(E\Big[\Delta t\sum\limits_{n=1}^{N}\parallel\nabla(\theta_{p}(t_{n})-\theta_{p}^{h,n})\parallel_{L^{2}(\Omega_{p})}^{2}\Big]\Bigg)^{\frac{1}{2}}\\
\leq&C(\Delta t+h^{k}).
\end{split}
\end{align}
\end{theorem}

$\mathbf{Proof.}$ First, the definition of the error functions are given below
\begin{align*}
\begin{split}
&e_{f,j}^{n}=a_{f,j}^{n}+b_{f,j}^{n},\  a_{f,j}^{n}:=\mathbf{u}_{f,j}(t_{n})-\Pi_{f}^{h}\mathbf{u}_{f,j}(t_{n}),\  b_{f,j}^{n}:=\Pi_{f}^{h}\mathbf{u}_{f,j}(t_{n})-\mathbf{u}_{f,j}^{h,n},\\
&\xi_{f,j}^{n}=c_{f,j}^{n}+d_{f,j}^{n},\  c_{f,j}^{n}:=p_{f,j}(t_{n})-\rho_{f}^{h}p_{f,j}(t_{n}),\  d_{f,j}^{n}:=\rho_{f}^{h}p_{f,j}(t_{n})-p_{f,j}^{h,n},\\
&\eta_{f,j}^{n}=\tau_{f,j}^{n}+\varepsilon_{f,j}^{n},\  \tau_{f,j}^{n}:=\theta_{f,j}(t_{n})-Q_{f}^{h}\theta_{f,j}(t_{n}),\  \varepsilon_{f,j}^{n}:=Q_{f}^{h}\theta_{f,j}(t_{n})-\theta_{f,j}^{h,n},\\
\end{split}
\end{align*}
\begin{align*}
\begin{split}
&e_{p,j}^{n}=a_{p,j}^{n}+b_{p,j}^{n},\  a_{p,j}^{n}:=\mathbf{u}_{p,j}(t_{n})-\Pi_{p}^{h}\mathbf{u}_{p,j}(t_{n}),\  b_{p,j}^{n}:=\Pi_{p}^{h}\mathbf{u}_{p,j}(t_{n})-\mathbf{u}_{p,j}^{h,n},\\
&\xi_{p,j}^{n}=c_{p,j}^{n}+d_{p,j}^{n},\  c_{p,j}^{n}:=\phi_{p,j}(t_{n})-\rho_{p}^{h}\phi_{p,j}(t_{n}),\  d_{p,j}^{n}:=\rho_{p}^{h}\phi_{p,j}(t_{n})-\phi_{p,j}^{h,n},\\
&\eta_{p,j}^{n}=\tau_{p,j}^{n}+\varepsilon_{p,j}^{n},\  \tau_{p,j}^{n}:=\theta_{p,j}(t_{n})-Q_{p}^{h}\theta_{p,j}(t_{n}),\  \varepsilon_{p,j}^{n}:=Q_{p}^{h}\theta_{p,j}(t_{n})-\theta_{p,j}^{h,n}.
\end{split}
\end{align*}
For $\forall (\mathbf{v}_{f}^{h},\ q^{h},\ \varphi^{h};\ \mathbf{v}_{p}^{h},\ \psi^{h},\ \omega^{h}) \in (V_{f}^{h},\ U_{f}^{h},\ W_{f}^{h};\
V_{p}^{h},\ U_{p}^{h},\ W_{p}^{h}),$ the true solution $(\mathbf{u}_{f,j},p_{f,j},\theta_{f,j};\ \mathbf{u}_{p,j},\phi_{p,j},\theta_{p,j})$ satisfies
\begin{align}\label{wucha1}
\begin{split}
&\Bigg(\frac{\mathbf{u}_{f,j}(t_{n+1})-\mathbf{u}_{f,j}(t_{n})}{\Delta t},\mathbf{v}_{f}^{h}\Bigg)_{\Omega_{f}}
+\frac{C_{a}}{L^{2}}\Bigg(K_{j}\frac{\mathbf{u}_{p,j}(t_{n+1})-\mathbf{u}_{p,j}(t_{n})}{\Delta t},\mathbf{v}_{p}^{h}\Bigg)_{\Omega_{p}}
-(p_{f,j}(t_{n+1}),\nabla\cdot\mathbf{v}_{f}^{h})_{\Omega_{f}}\\
&+Pr(\nabla\mathbf{u}_{f,j}(t_{n+1}),\nabla\mathbf{v}_{f}^{h})_{\Omega_{f}}
+c(\mathbf{u}_{f,j}(t_{n+1}),\mathbf{u}_{f,j}(t_{n+1}),\mathbf{v}_{f}^{h})_{\Omega_{f}}
+Pr(\mathbf{u}_{p,j}(t_{n+1}),\mathbf{v}_{p}^{h})_{\Omega_{p}}\\
&-\frac{K_{j}}{L^{2}}(\phi_{p,j}(t_{n+1}),\nabla\cdot\mathbf{v}_{p}^{h})_{\Omega_{p}}\\
=&\Bigg(\frac{\mathbf{u}_{f,j}(t_{n+1})-\mathbf{u}_{f,j}(t_{n})}{\Delta t}-\mathbf{u}_{f,j,t}(t_{n+1}),\mathbf{v}_{f}^{h}\Bigg)_{\Omega_{f}}
+\frac{C_{a}}{L^{2}}\Bigg(K_{j}\frac{\mathbf{u}_{p,j}(t_{n+1})-\mathbf{u}_{p,j}(t_{n})}{\Delta t}-\mathbf{u}_{p,j,t}(t_{n+1}),\mathbf{v}_{p}^{h}\Bigg)_{\Omega_{p}}\\
&+PrRa(\mathbf{\xi}\theta_{f,j}(t_{n+1}),\mathbf{v}_{f}^{h})_{\Omega_{f}}
+\frac{PrRa}{L^{2}}(K_{j}\mathbf{\xi}\theta_{p,j}(t_{n+1}),\mathbf{v}_{p}^{h})_{\Omega_{p}}
+(\mathbf{f}_{f,j}(t_{n+1}),\mathbf{v}_{f}^{h})_{\Omega_{f}},
\end{split}
\end{align}

\begin{align}\label{wucha25}
\begin{split}
&(q^{h},\nabla\cdot\mathbf{u}_{f,j}(t_{n+1}))_{\Omega_{f}}=0,
\end{split}
\end{align}

\begin{align}\label{wucha26}
\begin{split}
&\frac{K_{j}}{L^{2}}(\psi^{h},\nabla\cdot\mathbf{u}_{p,j}(t_{n+1}))_{\Omega_{p}}=0,
\end{split}
\end{align}

\begin{align}\label{wucha2}
\begin{split}
&\Bigg(\frac{\theta_{f,j}(t_{n+1})-\theta_{f,j}(t_{n})}{\Delta t},\varphi^{h}\Bigg)_{\Omega_{f}}
+\Bigg(\frac{\theta_{p,j}(t_{n+1})-\theta_{p,j}(t_{n})}{\Delta t},\omega^{h}\Bigg)_{\Omega_{p}}
+k_{f}(\nabla\theta_{f,j}(t_{n+1}),\nabla \varphi^{h})_{\Omega_{f}}\\
&+k_{p}(\nabla\theta_{p,j}(t_{n+1}),\nabla \omega^{h})_{\Omega_{p}}
+t_{f}(\mathbf{u}_{f,j}(t_{n+1}),\theta_{f,j}(t_{n+1}),\varphi^{h})_{\Omega_{f}}
+t_{p}(\mathbf{u}_{p,j}(t_{n+1}),\theta_{p,j}(t_{n+1}),\omega^{h})_{\Omega_{p}}\\
&-k_{f}\int_{\uppercase\expandafter{\romannumeral1}}\mathbf{n}_{f}\cdot\nabla\theta_{f,j}(t_{n+1})(\varphi^{h}-\omega^{h})dl
+\frac{k_{f}\gamma}{h}\int_{\uppercase\expandafter{\romannumeral1}}(\theta_{f,j}(t_{n+1})-\theta_{p,j}(t_{n+1}))(\varphi^{h}-\omega^{h})dl\\
=&\Bigg(\frac{\theta_{f,j}(t_{n+1})-\theta_{f,j}(t_{n})}{\Delta t}-\theta_{f,j,t}(t_{n+1}),\varphi^{h}\Bigg)_{\Omega_{f}}
+\Bigg(\frac{\theta_{p,j}(t_{n+1})-\theta_{p,j}(t_{n})}{\Delta t}-\theta_{p,j,t}(t_{n+1}),\omega^{h}\Bigg)_{\Omega_{p}}\\
&+(\mathbf{\Upsilon}_{f,j}(t_{n+1}),\varphi^{h})_{\Omega_{f}}
+(\mathbf{\Upsilon}_{p,j}(t_{n+1}),\omega^{h})_{\Omega_{p}}.
\end{split}
\end{align}
Subtracting (\ref{suanfa 1}) - (\ref{suanfa 6}) from (\ref{wucha1}) - (\ref{wucha2}) to get
\begin{align}\label{wucha3}
\begin{split}
&\Bigg(\frac{e_{f,j}^{n+1}-e_{f,j}^{n}}{\Delta t},\mathbf{v}_{f}^{h}\Bigg)_{\Omega_{f}}
+\frac{C_{a}}{L^{2}}\Bigg(K_{j}\frac{e_{p,j}^{n+1}-e_{p,j}^{n}}{\Delta t},\mathbf{v}_{p}^{h}\Bigg)_{\Omega_{p}}
-(\xi_{f,j}^{n+1},\nabla\cdot\mathbf{v}_{f}^{h})_{\Omega_{f}}
+Pr(\nabla e_{f,j}^{n+1},\nabla\mathbf{v}_{f}^{h})_{\Omega_{f}}\\
&+c(\mathbf{u}_{f,j}(t_{n+1}),\mathbf{u}_{f,j}(t_{n+1}),\mathbf{v}_{f}^{h})_{\Omega_{f}}
-c(\mathbf{u}_{f,j}^{h,n},\mathbf{u}_{f,j}^{h,n+1},\mathbf{v}_{f}^{h})_{\Omega_{f}}
+Pr(e_{p,j}^{n+1},\mathbf{v}_{p}^{h})_{\Omega_{p}}
-\frac{K_{j}}{L^{2}}(\xi_{p,j}^{n+1},\nabla\cdot\mathbf{v}_{p}^{h})_{\Omega_{p}}\\
=&\Bigg(\frac{\mathbf{u}_{f,j}(t_{n+1})-\mathbf{u}_{f,j}(t_{n})}{\Delta t}-\mathbf{u}_{f,j,t}(t_{n+1}),\mathbf{v}_{f}^{h}\Bigg)_{\Omega_{f}}
+\frac{C_{a}}{L^{2}}\Bigg(K_{j}\frac{\mathbf{u}_{p,j}(t_{n+1})-\mathbf{u}_{p,j}(t_{n})}{\Delta t}-\mathbf{u}_{p,j,t}(t_{n+1}),\mathbf{v}_{p}^{h}\Bigg)_{\Omega_{p}}\\
&+PrRa(\mathbf{\xi}\theta_{f,j}(t_{n+1}),\mathbf{v}_{f}^{h})_{\Omega_{f}}
-PrRa(\mathbf{\xi}\theta_{f,j}^{h,n},\mathbf{v}_{f}^{h})_{\Omega_{f}}
+\frac{PrRaK_{j}}{L^{2}}(\mathbf{\xi}\theta_{p,j}(t_{n+1}),\mathbf{v}_{p}^{h})_{\Omega_{p}}\\
&-\frac{PrRaK_{j}}{L^{2}}(\mathbf{\xi}\theta_{p,j}^{h,n},\mathbf{v}_{p}^{h})_{\Omega_{p}},
\end{split}
\end{align}

\begin{align}\label{wucha4}
\begin{split}
&(q^{h},\nabla\cdot e_{f,j}^{n+1})_{\Omega_{f}}=0,
\end{split}
\end{align}

\begin{align}\label{wucha5}
\begin{split}
&\frac{K_{j}}{L^{2}}(\psi^{h},\nabla\cdot e_{p,j}^{n+1})_{\Omega_{p}}=0,
\end{split}
\end{align}

\begin{align}\label{wucha6}
\begin{split}
&\Bigg(\frac{\eta_{f,j}^{n+1}-\eta_{f,j}^{n}}{\Delta t},\varphi^{h}\Bigg)_{\Omega_{f}}
+\Bigg(\frac{\eta_{p,j}^{n+1}-\eta_{p,j}^{n}}{\Delta t},\omega^{h}\Bigg)_{\Omega_{p}}
+k_{f}(\nabla\eta_{f,j}^{n+1},\nabla \varphi^{h})_{\Omega_{f}}
+k_{p}(\nabla\eta_{p,j}^{n+1},\nabla \omega^{h})_{\Omega_{p}}\\
&+t_{f}(\mathbf{u}_{f,j}(t_{n+1}),\theta_{f,j}(t_{n+1}),\varphi^{h})_{\Omega_{f}}
-t_{f}(\mathbf{u}_{f,j}^{h,n},\theta_{f,j}^{h,n+1},\varphi^{h})_{\Omega_{f}}
+t_{p}(\mathbf{u}_{p,j}(t_{n+1}),\theta_{p,j}(t_{n+1}),\omega^{h})_{\Omega_{p}}\\
&-t_{p}(\mathbf{u}_{p,j}^{h,n},\theta_{p,j}^{h,n+1},\omega^{h})_{\Omega_{p}}
-k_{f}\int_{\uppercase\expandafter{\romannumeral1}}\mathbf{n}_{f}\cdot\nabla\theta_{f,j}(t_{n+1})(\varphi^{h}-\omega^{h})dl
+k_{f}\int_{\uppercase\expandafter{\romannumeral1}}\mathbf{n}_{f}\cdot\nabla\theta_{f,j}^{h,n}(\varphi^{h}-\omega^{h})dl\\
&+\frac{k_{f}\gamma}{h}\int_{\uppercase\expandafter{\romannumeral1}}(\theta_{f,j}(t_{n+1})-\theta_{p,j}(t_{n+1}))(\varphi^{h}-\omega^{h})dl
-\frac{k_{f}\gamma}{h}\int_{\uppercase\expandafter{\romannumeral1}}(\theta_{f,j}^{h,n+1}-\theta_{p,j}^{h,n})\varphi^{h}dl
+\frac{k_{f}\gamma}{h}\int_{\uppercase\expandafter{\romannumeral1}}(\theta_{f,j}^{h,n+1}-\theta_{p,j}^{h,n+1})\omega^{h}dl\\
=&\Bigg(\frac{\theta_{f,j}(t_{n+1})-\theta_{f,j}(t_{n})}{\Delta t}-\theta_{f,j,t}(t_{n+1}),\varphi^{h}\Bigg)_{\Omega_{f}}
+\Bigg(\frac{\theta_{p,j}(t_{n+1})-\theta_{p,j}(t_{n})}{\Delta t}-\theta_{p,j,t}(t_{n+1}),\omega^{h}\Bigg)_{\Omega_{p}}.
\end{split}
\end{align}
Choosing $(\mathbf{v}_{f}^{h}, q^{h}, \varphi^{h})=2\Delta t(b_{f,j}^{n+1}, d_{f,j}^{n+1}, \varepsilon_{f,j}^{n+1})$, and  $(\mathbf{v}_{p}^{h}, \psi^{h}, \omega^{h})=2\Delta t(b_{p,j}^{n+1}, d_{p,j}^{n+1}, \varepsilon_{p,j}^{n+1})$ in the equations (\ref{wucha3}) - (\ref{wucha6}), adding these equations together, applying the condition $2(a-b,a)=\parallel a\parallel^{2}-\parallel b\parallel^{2}+\parallel a-b\parallel^{2}$ and by the properties of the project operator (\ref{touyin1}), we obtain the following inequality
\begin{align}\label{wucha7}
\begin{split}
&\Big[\parallel b_{f,j}^{n+1}\parallel_{L^{2}(\Omega_{f})}^{2}-\parallel b_{f,j}^{n}\parallel_{L^{2}(\Omega_{f})}^{2}
+\parallel b_{f,j}^{n+1}-b_{f,j}^{n}\parallel_{L^{2}(\Omega_{f})}^{2}\Big]
+2\Delta tPr\parallel\nabla b_{f,j}^{n+1}\parallel^{2}_{L^{2}(\Omega_{f})}\\
&+\frac{C_{a}k_{min}}{L^{2}}\Big[\parallel b_{p,j}^{n+1}\parallel_{L^{2}(\Omega_{p})}^{2}
-\parallel b_{p,j}^{n}\parallel_{L^{2}(\Omega_{p})}^{2}
+\parallel b_{p,j}^{n+1}-b_{p,j}^{n}\parallel_{L^{2}(\Omega_{p})}^{2}\Big]
+2\Delta tPr\parallel b_{p,j}^{n+1}\parallel^{2}_{L^{2}(\Omega_{p})}\\
&+\Big[\parallel \varepsilon_{f,j}^{n+1}\parallel_{L^{2}(\Omega_{f})}^{2}-\parallel \varepsilon_{f,j}^{n}\parallel_{L^{2}(\Omega_{f})}^{2}
+\parallel \varepsilon_{f,j}^{n+1}-\varepsilon_{f,j}^{n}\parallel_{L^{2}(\Omega_{f})}^{2}\Big]
+2\Delta tk_{f}\parallel\nabla \varepsilon_{f,j}^{n+1}\parallel^{2}_{L^{2}(\Omega_{f})}\\
&+\Big[\parallel \varepsilon_{p,j}^{n+1}\parallel_{L^{2}(\Omega_{p})}^{2}-\parallel \varepsilon_{p,j}^{n}\parallel_{L^{2}(\Omega_{p})}^{2}
+\parallel \varepsilon_{p,j}^{n+1}-\varepsilon_{p,j}^{n}\parallel_{L^{2}(\Omega_{p})}^{2}\Big]
+2\Delta tk_{p}\parallel\nabla \varepsilon_{p,j}^{n+1}\parallel^{2}_{L^{2}(\Omega_{p})}\\
\leq&-2\Delta t\Bigg(\frac{a_{f,j}^{n+1}-a_{f,j}^{n}}{\Delta t},b_{f,j}^{n+1}\Bigg)_{\Omega_{f}}
-2\Delta t\frac{C_{a}k_{max}}{L^{2}}\Bigg(\frac{a_{p,j}^{n+1}-a_{p,j}^{n}}{\Delta t},b_{p,j}^{n+1}\Bigg)_{\Omega_{p}}
-2\Delta t\Bigg(\frac{\tau_{f,j}^{n+1}-\tau_{f,j}^{n}}{\Delta t},\varepsilon_{f,j}^{n+1}\Bigg)_{\Omega_{f}}\\
&-2\Delta t\Bigg(\frac{\tau_{p,j}^{n+1}-\tau_{p,j}^{n}}{\Delta t},\varepsilon_{p,j}^{n+1}\Bigg)_{\Omega_{p}}
+2\Delta t\Big[c(\mathbf{u}_{f,j}(t_{n+1}),\mathbf{u}_{f,j}(t_{n+1}),b_{f,j}^{n+1})_{\Omega_{f}}
-c(\mathbf{u}_{f,j}^{h,n},\mathbf{u}_{f,j}^{h,n+1},b_{f,j}^{n+1})_{\Omega_{f}}\\
&+t_{f}(\mathbf{u}_{f,j}(t_{n+1}),\theta_{f,j}(t_{n+1}),\varepsilon_{f,j}^{n+1})_{\Omega_{f}}
-t_{f}(\mathbf{u}_{f,j}^{h,n},\theta_{f,j}^{h,n+1},\varepsilon_{f,j}^{n+1})_{\Omega_{f}}
+t_{p}(\mathbf{u}_{p,j}(t_{n+1}),\theta_{p,j}(t_{n+1}),\varepsilon_{p,j}^{n+1})_{\Omega_{p}}\\
&-t_{p}(\mathbf{u}_{p,j}^{h,n},\theta_{p,j}^{h,n+1},\varepsilon_{p,j}^{n+1})_{\Omega_{p}}\Big]
+2\Delta t\Bigg(\frac{\mathbf{u}_{f,j}(t_{n+1})-\mathbf{u}_{f,j}(t_{n})}{\Delta t}-\mathbf{u}_{f,j,t}(t_{n+1}),b_{f,j}^{n+1}\Bigg)_{\Omega_{f}}\\
&+\frac{2\Delta tC_{a}k_{max}}{L^{2}}\Bigg(\frac{\mathbf{u}_{p,j}(t_{n+1})-\mathbf{u}_{p,j}(t_{n})}{\Delta t}
-\mathbf{u}_{p,j,t}(t_{n+1}),b_{p,j}^{n+1}\Bigg)_{\Omega_{p}}
+2\Delta tPrRa(\mathbf{\xi}\theta_{f,j}(t_{n+1}),b_{f,j}^{n+1})_{\Omega_{f}}\\
&-2\Delta tPrRa(\mathbf{\xi}\theta_{f,j}^{h,n},b_{f,j}^{n+1})_{\Omega_{f}}
+\frac{2\Delta tPrRa}{L^{2}}(K_{j}\mathbf{\xi}\theta_{p,j}(t_{n+1}),b_{p,j}^{n+1})_{\Omega_{p}}
-\frac{2\Delta tPrRa}{L^{2}}(K_{j}\mathbf{\xi}\theta_{p,j}^{h,n},b_{p,j}^{n+1})_{\Omega_{p}}\\
&+2\Delta t\Bigg(\frac{\theta_{f,j}(t_{n+1})-\theta_{f,j}(t_{n})}{\Delta t}-\theta_{f,j,t}(t_{n+1}),\varepsilon_{f,j}^{n+1}\Bigg)_{\Omega_{f}}
+2\Delta t\Bigg(\frac{\theta_{p,j}(t_{n+1})-\theta_{p,j}(t_{n})}{\Delta t}-\theta_{p,j,t}(t_{n+1}),\varepsilon_{p,j}^{n+1}\Bigg)_{\Omega_{p}}\\
&+2\Delta tk_{f}\int_{\uppercase\expandafter{\romannumeral1}}\mathbf{n}_{f}\cdot\nabla(\theta_{f,j}(t_{n+1})-\theta_{f,j}(t_{n}))(\varepsilon_{f,j}^{n+1}-\varepsilon_{p,j}^{n+1})dl
+2\Delta tk_{f}\int_{\uppercase\expandafter{\romannumeral1}}\mathbf{n}_{f}\cdot\nabla\varepsilon_{f,j}^{n}(\varepsilon_{f,j}^{n+1}-\varepsilon_{p,j}^{n+1})dl\\
&-\frac{2\Delta tk_{f}\gamma}{h}\int_{\uppercase\expandafter{\romannumeral1}}(\theta_{f,j}(t_{n+1})-\theta_{p,j}(t_{n+1}))(\varepsilon_{f,j}^{n+1}-\varepsilon_{p,j}^{n+1})dl
+\frac{2\Delta tk_{f}\gamma}{h}\int_{\uppercase\expandafter{\romannumeral1}}(\theta_{f,j}^{h,n+1}-\theta_{p,j}^{h,n})\varepsilon_{f,j}^{n+1}dl\\
&-\frac{2\Delta tk_{f}\gamma}{h}\int_{\uppercase\expandafter{\romannumeral1}}(\theta_{f,j}^{h,n+1}-\theta_{p,j}^{h,n+1})\varepsilon_{p,j}^{n+1}dl.
\end{split}
\end{align}
Then, we estimate the first two terms on the right hand side of (\ref{wucha7}) by the Cauchy-Schwarz inequality and Young's inequality

\begin{align}\label{wucha8}
\begin{split}
&-2\Delta t\Bigg(\frac{a_{f,j}^{n+1}-a_{f,j}^{n}}{\Delta t},b_{f,j}^{n+1}\Bigg)_{\Omega_{f}}
-2\Delta t\frac{C_{a}k_{max}}{L^{2}}\Bigg(\frac{a_{p,j}^{n+1}-a_{p,j}^{n}}{\Delta t},b_{p,j}^{n+1}\Bigg)_{\Omega_{p}}\\
\leq&\frac{4\Delta tC_{f}^{2}}{Pr}\Big\|\frac{a_{f,j}^{n+1}-a_{f,j}^{n}}{\Delta t}\Big\|_{L^{2}(\Omega_{f})}^{2}
+\frac{1}{4}Pr\Delta t\parallel\nabla b_{f,j}^{n+1}\parallel_{L^{2}(\Omega_{f})}^{2}
+\frac{4\Delta tC_{a}^{2}k_{max}^{2}}{L^{4}Pr}\Big\|\frac{a_{p,j}^{n+1}-a_{p,j}^{n}}{\Delta t}\Big\|_{L^{2}(\Omega_{p})}^{2}\\
&+\frac{1}{4}Pr\Delta t\parallel b_{p,j}^{n+1}\parallel_{L^{2}(\Omega_{p})}^{2}\\
\leq&\frac{4\Delta tC_{f}^{2}}{Pr}\Big\|\frac{1}{\Delta t}\int_{t_{n}}^{t_{n+1}}a_{f,j,t}dt\Big\|_{L^{2}(\Omega_{f})}^{2}
+\frac{1}{4}Pr\Delta t\parallel\nabla b_{f,j}^{n+1}\parallel_{L^{2}(\Omega_{f})}^{2}
+\frac{4\Delta tC_{a}^{2}k_{max}^{2}}{L^{4}Pr}\Big\|\frac{1}{\Delta t}\int_{t_{n}}^{t_{n+1}}a_{p,j,t}dt\Big\|_{L^{2}(\Omega_{p})}^{2}\\
&+\frac{1}{4}Pr\Delta t\parallel b_{p,j}^{n+1}\parallel_{L^{2}(\Omega_{p})}^{2}\\
\leq&\frac{4C_{f}^{2}}{Pr}\int_{t_{n}}^{t_{n+1}}\parallel a_{f,j,t}\parallel_{L^{2}(\Omega_{f})}^{2}dt
+\frac{1}{4}Pr\Delta t\parallel\nabla b_{f,j}^{n+1}\parallel_{L^{2}(\Omega_{f})}^{2}
+\frac{4C_{a}^{2}k_{max}^{2}}{L^{4}Pr}\int_{t_{n}}^{t_{n+1}}\parallel a_{p,j,t}\parallel_{L^{2}(\Omega_{p})}^{2}dt\\
&+\frac{1}{4}Pr\Delta t\parallel b_{p,j}^{n+1}\parallel_{L^{2}(\Omega_{p})}^{2}\\
\leq&Ch^{2k+2}\int_{t_{n}}^{t_{n+1}}\parallel\mathbf{u}_{f,j,t}\parallel_{H^{k+1}(\Omega_{f})}^{2}dt
+Ch^{2k+2}\int_{t_{n}}^{t_{n+1}}\parallel\mathbf{u}_{p,j,t}\parallel_{H^{k+1}(\Omega_{p})}^{2}dt
+\frac{1}{4}Pr\Delta t\parallel\nabla b_{f,j}^{n+1}\parallel_{L^{2}(\Omega_{f})}^{2}\\
&+\frac{1}{4}Pr\Delta t\parallel b_{p,j}^{n+1}\parallel_{L^{2}(\Omega_{p})}^{2},
\end{split}
\end{align}
With the same technique again, we have following result
\begin{align}\label{wucha9}
\begin{split}
&-2\Delta t\Bigg(\frac{\tau_{f,j}^{n+1}-\tau_{f,j}^{n}}{\Delta t},\varepsilon_{f,j}^{n+1}\Bigg)_{\Omega_{f}}
-2\Delta t\Bigg(\frac{\tau_{p,j}^{n+1}-\tau_{p,j}^{n}}{\Delta t},\varepsilon_{p,j}^{n+1}\Bigg)_{\Omega_{p}}\\
\leq&\frac{4C_{t}^{2}}{k_{f}}\int_{t_{n}}^{t_{n+1}}\parallel \tau_{f,j,t}\parallel_{L^{2}(\Omega_{f})}^{2}dt
+\frac{1}{4}k_{f}\Delta t\parallel\nabla \varepsilon_{f,j}^{n+1}\parallel_{L^{2}(\Omega_{f})}^{2}
+\frac{4\tilde{C_{t}}^{2}}{k_{p}}\int_{t_{n}}^{t_{n+1}}\parallel \tau_{p,j,t}\parallel_{L^{2}(\Omega_{p})}^{2}dt\\
&+\frac{1}{4}k_{p}\Delta t\parallel\nabla \varepsilon_{p,j}^{n+1}\parallel_{L^{2}(\Omega_{p})}^{2}\\
\leq&Ch^{2k+2}\int_{t_{n}}^{t_{n+1}}\parallel\theta_{f,j,t}\parallel_{H^{k+1}(\Omega_{f})}^{2}dt
+Ch^{2k+2}\int_{t_{n}}^{t_{n+1}}\parallel\theta_{p,j,t}\parallel_{H^{k+1}(\Omega_{p})}^{2}dt\\
&+\frac{1}{4}k_{f}\Delta t\parallel\nabla \varepsilon_{f,j}^{n+1}\parallel_{L^{2}(\Omega_{f})}^{2}
+\frac{1}{4}k_{p}\Delta t\parallel \nabla\varepsilon_{p,j}^{n+1}\parallel_{L^{2}(\Omega_{p})}^{2}.
\end{split}
\end{align}
Using the Young's inequality and lemma \ref{lemma2.1}, we have
\begin{align}\label{wucha10}
\begin{split}
&2\Delta tc(\mathbf{u}_{f,j}(t_{n+1}),\mathbf{u}_{f,j}(t_{n+1}),b_{f,j}^{n+1})_{\Omega_{f}}
-2\Delta tc(\mathbf{u}_{f,j}^{h,n},\mathbf{u}_{f,j}^{h,n+1},b_{f,j}^{n+1})_{\Omega_{f}}\\
=&2\Delta t\Big[c(\mathbf{u}_{f,j}(t_{n+1})-\mathbf{u}_{f,j}(t_{n}),\mathbf{u}_{f,j}(t_{n+1}),b_{f,j}^{n+1})_{\Omega_{f}}
+c(a_{f,j}^{n},\mathbf{u}_{f,j}(t_{n+1}),b_{f,j}^{n+1})_{\Omega_{f}}\\
&+c(b_{f,j}^{n},\mathbf{u}_{f,j}(t_{n+1}),b_{f,j}^{n+1})_{\Omega_{f}}
+c(\mathbf{u}_{f,j}^{h,n},a_{f,j}^{n+1},b_{f,j}^{n+1})_{\Omega_{f}}\Big]\\
\leq&2\Delta t\parallel\mathbf{u}_{f,j}(t_{n+1})-\mathbf{u}_{f,j}(t_{n})\parallel_{L^{2}(\Omega_{f})}
\parallel\mathbf{u}_{f,j}(t_{n+1})\parallel_{H^{2}(\Omega_{f})}\parallel\nabla b_{f,j}^{n+1}\parallel_{L^{2}(\Omega_{f})}\\
&+2\Delta t\parallel a_{f,j}^{n}\parallel_{L^{2}(\Omega_{f})}
\parallel\mathbf{u}_{f,j}(t_{n+1})\parallel_{H^{2}(\Omega_{f})}\parallel\nabla b_{f,j}^{n+1}\parallel_{L^{2}(\Omega_{f})}\\
&+2\Delta t\parallel b_{f,j}^{n}\parallel_{L^{2}(\Omega_{f})}
\parallel\mathbf{u}_{f,j}(t_{n+1})\parallel_{H^{2}(\Omega_{f})}\parallel\nabla b_{f,j}^{n+1}\parallel_{L^{2}(\Omega_{f})}\\
&+2\Delta t\parallel \nabla\mathbf{u}_{f,j}^{h,n}\parallel_{L^{2}(\Omega_{f})}
\parallel\nabla a_{f,j}^{n+1}\parallel_{L^{2}(\Omega_{f})}\parallel\nabla b_{f,j}^{n+1}\parallel_{L^{2}(\Omega_{f})}\\
\leq&\frac{8\Delta t^{3}}{Pr}\parallel\mathbf{u}_{f,j,t}(t_{n+1})\parallel_{L^{2}(\Omega_{f})}^{2}\parallel\mathbf{u}_{f,j}(t_{n+1})\parallel_{H^{2}(\Omega_{f})}^{2}
+\frac{8\Delta t^{4}}{3Pr}\int_{t_{n}}^{t_{n+1}}\parallel\mathbf{u}_{f,j,tt}\parallel_{L^{2}(\Omega_{f})}^{2}dt
\parallel \mathbf{u}_{f,j}(t_{n+1})\parallel_{H^{2}(\Omega_{f})}^{2}\\
&+\frac{8\Delta t}{Pr}Ch^{2k+2}\parallel\mathbf{u}_{f,j}(t_{n})\parallel_{H^{k+1}(\Omega_{f})}^{2}
\parallel \mathbf{u}_{f,j}(t_{n+1})\parallel_{H^{2}(\Omega_{f})}^{2}
+\frac{8\Delta t}{Pr}\parallel b_{f,j}^{n}\parallel_{L^{2}(\Omega_{f})}^{2}
\parallel \mathbf{u}_{f,j}(t_{n+1})\parallel_{H^{2}(\Omega_{f})}^{2}\\
&+\frac{8\Delta t}{Pr}Ch^{2k}\parallel\mathbf{u}_{f,j}(t_{n+1})\parallel_{H^{k+1}(\Omega_{f})}^{2}
\parallel\nabla\mathbf{u}_{f,j}^{h,n}\parallel_{L^{2}(\Omega_{f})}^{2}
+\frac{1}{2}Pr\Delta t\parallel\nabla b_{f,j}^{n+1}\parallel_{L^{2}(\Omega_{f})}^{2}.
\end{split}
\end{align}
Applying a similar proof method as (\ref{wucha10}), we get
\begin{align}\label{wucha11}
\begin{split}
&2\Delta tt_{f}(\mathbf{u}_{f,j}(t_{n+1}),\theta_{f,j}(t_{n+1}),\varepsilon_{f,j}^{n+1})_{\Omega_{f}}
-2\Delta tt_{f}(\mathbf{u}_{f,j}^{h,n},\theta_{f,j}^{h,n+1},\varepsilon_{f,j}^{n+1})_{\Omega_{f}}\\
\leq&C(\Delta t^{2}+h^{2k})
+\frac{8\Delta t}{k_{f}}\parallel b_{f,j}^{n}\parallel_{L^{2}(\Omega_{f})}^{2}
\parallel \theta_{f,j}(t_{n+1})\parallel_{H^{2}(\Omega_{f})}^{2}
+\frac{1}{2}k_{f}\Delta t\parallel\nabla \varepsilon_{f,j}^{n+1}\parallel_{L^{2}(\Omega_{f})}^{2},
\end{split}
\end{align}

\begin{align}\label{wucha12}
\begin{split}
&2\Delta tt_{p}(\mathbf{u}_{p,j}(t_{n+1}),\theta_{p,j}(t_{n+1}),\varepsilon_{p,j}^{n+1})_{\Omega_{p}}
-2\Delta tt_{p}(\mathbf{u}_{p,j}^{h,n},\theta_{p,j}^{h,n+1},\varepsilon_{p,j}^{n+1})_{\Omega_{p}}\\
=&2\Delta tt_{p}(\mathbf{u}_{p,j}(t_{n+1})-\mathbf{u}_{p,j}(t_{n}),\theta_{p,j}(t_{n+1}),\varepsilon_{p,j}^{n+1})_{\Omega_{p}}
+2\Delta tt_{p}(a_{p,j}^{n},\theta_{p,j}(t_{n+1}),\varepsilon_{p,j}^{n+1})_{\Omega_{p}}\\
&+2\Delta tt_{p}(b_{p,j}^{n},\theta_{p,j}(t_{n+1}),\varepsilon_{p,j}^{n+1})_{\Omega_{p}}
+2\Delta tt_{p}(\mathbf{u}_{p,j}(t_{n}),\tau_{p,j}^{n+1},\varepsilon_{p,j}^{n+1})_{\Omega_{p}}\\
&+2\Delta tt_{p}(a_{p,j}^{n},\tau_{p,j}^{n+1},\varepsilon_{p,j}^{n+1})_{\Omega_{p}}
+2\Delta tt_{p}(b_{p,j}^{n},\tau_{p,j}^{n+1},\varepsilon_{p,j}^{n+1})_{\Omega_{p}}\\
\leq&2\Delta t\parallel\mathbf{u}_{p,j}(t_{n+1})-\mathbf{u}_{p,j}(t_{n})\parallel_{L^{2}(\Omega_{p})}
\parallel\theta_{p,j}(t_{n+1})\parallel_{H^{2}(\Omega_{p})}\parallel\nabla\varepsilon_{p,j}^{n+1}\parallel_{L^{2}(\Omega_{p})}\\
&+2\Delta t\parallel a_{p,j}^{n}\parallel_{L^{2}(\Omega_{p})}
\parallel\theta_{p,j}(t_{n+1})\parallel_{H^{2}(\Omega_{p})}\parallel\nabla\varepsilon_{p,j}^{n+1}\parallel_{L^{2}(\Omega_{p})}\\
&+2\Delta t\parallel b_{p,j}^{n}\parallel_{L^{2}(\Omega_{p})}
\parallel\theta_{p,j}(t_{n+1})\parallel_{H^{2}(\Omega_{p})}\parallel\nabla\varepsilon_{p,j}^{n+1}\parallel_{L^{2}(\Omega_{p})}\\
&+2\Delta t\parallel \mathbf{u}_{p,j}(t_{n})\parallel_{H^{2}(\Omega_{p})}
\parallel\nabla\tau_{p,j}^{n+1}\parallel_{L^{2}(\Omega_{p})}\parallel\nabla\varepsilon_{p,j}^{n+1}\parallel_{L^{2}(\Omega_{p})}\\
&+2\Delta t\parallel a_{p,j}^{n}\parallel_{L^{2}(\Omega_{p})}
\parallel\tau_{p,j}^{n+1}\parallel_{H^{2}(\Omega_{p})}\parallel\nabla\varepsilon_{p,j}^{n+1}\parallel_{L^{2}(\Omega_{p})}\\
&+2\Delta t\parallel b_{p,j}^{n}\parallel_{L^{2}(\Omega_{p})}
\parallel\tau_{p,j}^{n+1}\parallel_{H^{2}(\Omega_{p})}\parallel\nabla\varepsilon_{p,j}^{n+1}\parallel_{L^{2}(\Omega_{p})}\\
\leq&\frac{12 \Delta t^{3}}{k_{p}}\parallel\mathbf{u}_{p,j,t}(t_{n+1})\parallel_{L^{2}(\Omega_{p})}^{2}
\parallel\theta_{p,j}(t_{n+1})\parallel_{H^{2}(\Omega_{p})}^{2}
+\frac{4 \Delta t^{4}}{k_{p}}\int_{t_{n}}^{t_{n+1}}\parallel\mathbf{u}_{p,j,tt}\parallel_{L^{2}(\Omega_{p})}^{2}dt
\parallel\theta_{p,j}(t_{n+1})\parallel_{H^{2}(\Omega_{p})}^{2}\\
&+\frac{12 \Delta t}{k_{p}}\parallel b_{p,j}^{n}\parallel_{L^{2}(\Omega_{p})}^{2}
\parallel\theta_{p,j}(t_{n+1})\parallel_{H^{2}(\Omega_{p})}^{2}
+\frac{12 \Delta t}{k_{p}}Ch^{2k+2}\parallel \mathbf{u}_{p,j}(t_{n})\parallel_{H^{k+1}(\Omega_{p})}^{2}
\parallel\theta_{p,j}(t_{n+1})\parallel_{H^{2}(\Omega_{p})}^{2}\\
&+\frac{12 \Delta t}{k_{p}}Ch^{2k+2}\parallel \mathbf{u}_{p,j}(t_{n})\parallel_{H^{k+1}(\Omega_{p})}^{2}
\parallel\tau_{p,j}^{n+1}\parallel_{H^{2}(\Omega_{p})}^{2}
+\frac{12 \Delta t}{k_{p}}\parallel b_{p,j}^{n}\parallel_{L^{2}(\Omega_{p})}^{2}
\parallel \tau_{p,j}^{n+1}\parallel_{H^{2}(\Omega_{p})}^{2}\\
&+\frac{12 \Delta t}{k_{p}}Ch^{2k}\parallel \mathbf{u}_{p,j}(t_{n})\parallel_{H^{2}(\Omega_{p})}^{2}
\parallel\theta_{p,j}(t_{n+1})\parallel_{H^{k+1}(\Omega_{p})}^{2}
+\frac{1}{2}k_{p}\Delta t\parallel\nabla\varepsilon_{p,j}^{n+1}\parallel_{L^{2}(\Omega_{p})}^{2}.
\end{split}
\end{align}
By the Cauchy-Schwarz inequality and Young's inequality, we have the following estimates
\begin{align}\label{wucha13}
\begin{split}
&2\Delta tPrRa(\mathbf{\xi}\theta_{f,j}(t_{n+1}),b_{f,j}^{n+1})_{\Omega_{f}}
-2\Delta tPrRa(\mathbf{\xi}\theta_{f,j}^{h,n},b_{f,j}^{n+1})_{\Omega_{f}}\\
=&2\Delta tPrRa\mathbf{\xi}(\theta_{f,j}(t_{n+1})-\theta_{f,j}(t_{n}),b_{f,j}^{n+1})_{\Omega_{f}}
+2\Delta tPrRa(\mathbf{\xi}\eta_{f,j}^{n},b_{f,j}^{n+1})_{\Omega_{f}}\\
\leq&4\Delta t^{3}PrRa^{2}C_{f}^{2}\parallel \theta_{f,j,t}(t_{n+1})\parallel_{L^{2}(\Omega_{f})}^{2}
+\frac{4\Delta t^{4}PrRa^{2}C_{f}^{2}}{3}\int_{t_{n}}^{t_{n+1}}\parallel \theta_{f,j,tt}\parallel_{L^{2}(\Omega_{f})}^{2}dt\\
&+4\Delta tPrRa^{2}C_{f}^{2}Ch^{2k+2}\parallel \theta_{f,j}(t_{n})\parallel_{H^{k+1}(\Omega_{f})}^{2}
+4\Delta tPrRa^{2}C_{f}^{2}\parallel \varepsilon_{f,j}^{n}\parallel_{L^{2}(\Omega_{f})}^{2}\\
&+\frac{1}{2}\Delta tPr\parallel \nabla b_{f,j}^{n+1}\parallel_{L^{2}(\Omega_{f})}^{2},
\end{split}
\end{align}

\begin{align}\label{wucha14}
\begin{split}
&\frac{2\Delta tPrRa}{L^{2}}(K_{j}\mathbf{\xi}\theta_{p,j}(t_{n+1}),b_{p,j}^{n+1})_{\Omega_{p}}
-\frac{2\Delta tPrRa}{L^{2}}(K_{j}\mathbf{\xi}\theta_{p,j}^{h,n},b_{p,j}^{n+1})_{\Omega_{p}}\\
\leq&\frac{4\Delta t^{3}PrRa^{2}k_{max}^{2}}{L^{4}}\parallel \theta_{p,j,t}(t_{n+1})\parallel_{L^{2}(\Omega_{p})}^{2}
+\frac{4\Delta t^{4}PrRa^{2}k_{max}^{2}}{3L^{4}}\int_{t_{n}}^{t_{n+1}}\parallel \theta_{p,j,tt}\parallel_{L^{2}(\Omega_{p})}^{2}dt\\
&+\frac{4\Delta tPrRa^{2}k_{max}^{2}Ch^{2k+2}}{L^{4}}\parallel \theta_{p,j}(t_{n})\parallel_{H^{k+1}(\Omega_{p})}^{2}
+\frac{4\Delta tPrRa^{2}k_{max}^{2}}{L^{4}}\parallel \varepsilon_{p,j}^{n}\parallel_{L^{2}(\Omega_{p})}^{2}\\
&+\frac{1}{2}\Delta tPr\parallel b_{p,j}^{n+1}\parallel_{L^{2}(\Omega_{p})}^{2}.
\end{split}
\end{align}
Thanks to Taylor expansion with the integral remainder, we obtain
\begin{align*}
\begin{split}
&2\Delta t\Bigg(\frac{\mathbf{u}_{f,j}(t_{n+1})-\mathbf{u}_{f,j}(t_{n})}{\Delta t}-\mathbf{u}_{f,j,t}(t_{n+1}),b_{f,j}^{n+1}\Bigg)_{\Omega_{f}}\\
\end{split}
\end{align*}
\begin{align}\label{wucha15}
\begin{split}
=&2\bigg(\int_{t_{n}}^{t_{n+1}}(t-t_{n})\mathbf{u}_{f,j,tt}dt,b_{f,j}^{n+1}\bigg)_{\Omega_{f}}\\
\leq&\frac{4C_{f}^{2}\Delta t^{2}}{3Pr}\int_{t_{n}}^{t_{n+1}}\parallel\mathbf{u}_{f,j,tt}\parallel_{L^{2}(\Omega_{f})}^{2}dt
+\frac{1}{4}Pr\Delta t\parallel\nabla b_{f,j}^{n+1}\parallel_{L^{2}(\Omega_{f})}^{2}.
\end{split}
\end{align}
Using the same ways as (\ref{wucha15}), yield the following inequalities
\begin{align}\label{wucha16}
\begin{split}
&\frac{2\Delta tC_{a}k_{max}}{L^{2}}\Bigg(\frac{\mathbf{u}_{p,j}(t_{n+1})-\mathbf{u}_{p,j}(t_{n})}{\Delta t}
-\mathbf{u}_{p,j,t}(t_{n+1}),b_{p,j}^{n+1}\Bigg)_{\Omega_{p}}\\
\leq&\frac{4C_{a}^{2}\Delta t^{2}k_{max}^{2}}{3PrL^{4}}\int_{t_{n}}^{t_{n+1}}\parallel\mathbf{u}_{p,j,tt}\parallel_{L^{2}(\Omega_{p})}^{2}dt
+\frac{1}{4}Pr\Delta t\parallel b_{p,j}^{n+1}\parallel_{L^{2}(\Omega_{p})}^{2},
\end{split}
\end{align}

\begin{align}\label{wucha17}
\begin{split}
&2\Delta t\Bigg(\frac{\theta_{f,j}(t_{n+1})-\theta_{f,j}(t_{n})}{\Delta t}-\theta_{f,j,t}(t_{n+1}),\varepsilon_{f,j}^{n+1}\Bigg)_{\Omega_{f}}\\
\leq&\frac{4C_{t}^{2}\Delta t^{2}}{3k_{f}}\int_{t_{n}}^{t_{n+1}}\parallel\theta_{f,j,tt}\parallel_{L^{2}(\Omega_{f})}^{2}dt
+\frac{1}{4}k_{f}\Delta t\parallel\nabla \varepsilon_{f,j}^{n+1}\parallel_{L^{2}(\Omega_{f})}^{2},
\end{split}
\end{align}

\begin{align}\label{wucha18}
\begin{split}
&2\Delta t\Bigg(\frac{\theta_{p,j}(t_{n+1})-\theta_{p,j}(t_{n})}{\Delta t}-\theta_{p,j,t}(t_{n+1}),\varepsilon_{p,j}^{n+1}\Bigg)_{\Omega_{p}}\\
\leq&\frac{4\tilde{C_{t}}^{2}\Delta t^{2}}{3k_{p}}\int_{t_{n}}^{t_{n+1}}\parallel\theta_{p,j,tt}\parallel_{L^{2}(\Omega_{p})}^{2}dt
+\frac{1}{4}k_{p}\Delta t\parallel\nabla \varepsilon_{p,j}^{n+1}\parallel_{L^{2}(\Omega_{p})}^{2}.
\end{split}
\end{align}
Due to trace and inverse inequalities, we have
\begin{align}\label{wucha19}
\begin{split}
&2\Delta tk_{f}\int_{\uppercase\expandafter{\romannumeral1}}\mathbf{n}_{f}\cdot\nabla(\theta_{f,j}(t_{n+1})-\theta_{f,j}(t_{n}))(\varepsilon_{f,j}^{n+1}-\varepsilon_{p,j}^{n+1})dl
+2\Delta tk_{f}\int_{\uppercase\expandafter{\romannumeral1}}\mathbf{n}_{f}\cdot\nabla\varepsilon_{f,j}^{n}(\varepsilon_{f,j}^{n+1}-\varepsilon_{p,j}^{n+1})dl\\
\leq&\frac{4k_{f}C_{T}C_{inv}\Delta t}{\gamma}\parallel\nabla(\theta_{f,j}(t_{n+1})-\theta_{f,j}(t_{n}))\parallel_{L^{2}(\Omega_{f})}^{2}
+\frac{4k_{f}C_{T}C_{inv}\Delta t}{\gamma}\parallel\nabla\varepsilon_{f,j}^{n}\parallel_{L^{2}(\Omega_{f})}^{2}\\
&+\frac{k_{f}\gamma\Delta t}{2h}\parallel\varepsilon_{f,j}^{n+1}-\varepsilon_{p,j}^{n+1}\parallel_{L^{2}(\uppercase\expandafter{\romannumeral1})}^{2}\\
\leq&\frac{4k_{f}C_{T}C_{inv}\Delta t^{3}}{\gamma}\parallel\nabla\theta_{f,j,t}(t_{n+1})\parallel_{L^{2}(\Omega_{f})}^{2}
+\frac{4k_{f}C_{T}C_{inv}\Delta t^{4}}{3\gamma}\int_{t_{n}}^{t_{n+1}}\parallel\nabla\theta_{f,j,tt}\parallel_{L^{2}(\Omega_{f})}^{2}dt\\
&+\frac{4k_{f}C_{T}C_{inv}\Delta t}{\gamma}\parallel\nabla\varepsilon_{f,j}^{n}\parallel_{L^{2}(\Omega_{f})}^{2}
+\frac{k_{f}\gamma\Delta t}{2h}\parallel\varepsilon_{f,j}^{n+1}-\varepsilon_{p,j}^{n+1}\parallel_{L^{2}(\uppercase\expandafter{\romannumeral1})}^{2},
\end{split}
\end{align}
and
\begin{align}\label{wucha20}
\begin{split}
&-\frac{2\Delta tk_{f}\gamma}{h}\int_{\uppercase\expandafter{\romannumeral1}}(\theta_{f,j}(t_{n+1})-\theta_{p,j}(t_{n+1}))(\varepsilon_{f,j}^{n+1}-\varepsilon_{p,j}^{n+1})dl
+\frac{2\Delta tk_{f}\gamma}{h}\int_{\uppercase\expandafter{\romannumeral1}}(\theta_{f,j}^{h,n+1}-\theta_{p,j}^{h,n})\varepsilon_{f,j}^{n+1}dl\\
&-\frac{2\Delta tk_{f}\gamma}{h}\int_{\uppercase\expandafter{\romannumeral1}}(\theta_{f,j}^{h,n+1}-\theta_{p,j}^{h,n+1})\varepsilon_{p,j}^{n+1}dl\\
=&-\frac{2\Delta tk_{f}\gamma}{h}\parallel\varepsilon_{f,j}^{n+1}-\varepsilon_{p,j}^{n+1}\parallel_{L^{2}(\uppercase\expandafter{\romannumeral1})}^{2}
-\frac{2\Delta tk_{f}\gamma}{h}\int_{\uppercase\expandafter{\romannumeral1}}\tau_{f,j}^{n+1}(\varepsilon_{f,j}^{n+1}-\varepsilon_{p,j}^{n+1})dl\\
&+\frac{2\Delta tk_{f}\gamma}{h}\int_{\uppercase\expandafter{\romannumeral1}}\tau_{p,j}^{n+1}(\varepsilon_{f,j}^{n+1}-\varepsilon_{p,j}^{n+1})dl
+\frac{2\Delta tk_{f}\gamma}{h}\int_{\uppercase\expandafter{\romannumeral1}}(\theta_{p,j}^{h,n+1}-\theta_{p,j}^{h,n})\varepsilon_{f,j}^{n+1}dl.
\end{split}
\end{align}
For the second and third terms at the right end of (\ref{wucha20}), using the trace and inverse inequalities, we can obtain the following estimates
\begin{align}\label{wucha21}
\begin{split}
&\frac{2\Delta tk_{f}\gamma}{h}\int_{\uppercase\expandafter{\romannumeral1}}\tau_{f,j}^{n+1}(\varepsilon_{f,j}^{n+1}-\varepsilon_{p,j}^{n+1})dl\\
\leq &4\Delta tk_{f}\gamma C_{T}C_{inv}Ch^{2k}\parallel\theta_{f,j}(t_{n+1})\parallel_{H^{k+1}(\Omega_{f})}^{2}
+\frac{\Delta tk_{f}\gamma}{4h}\parallel\varepsilon_{f,j}^{n+1}-\varepsilon_{p,j}^{n+1}\parallel_{L^{2}(\uppercase\expandafter{\romannumeral1})}^{2},
\end{split}
\end{align}

\begin{align}\label{wucha22}
\begin{split}
&\frac{2\Delta tk_{f}\gamma}{h}\int_{\uppercase\expandafter{\romannumeral1}}\tau_{p,j}^{n+1}(\varepsilon_{f,j}^{n+1}-\varepsilon_{p,j}^{n+1})dl\\
\leq& 4\Delta tk_{f}\gamma \tilde{C}_{T}\tilde{C}_{inv}Ch^{2k}\parallel\theta_{p,j}(t_{n+1})\parallel_{H^{k+1}(\Omega_{f})}^{2}
+\frac{\Delta tk_{f}\gamma}{4h}\parallel\varepsilon_{f,j}^{n+1}-\varepsilon_{p,j}^{n+1}\parallel_{L^{2}(\uppercase\expandafter{\romannumeral1})}^{2}.
\end{split}
\end{align}
Similarly, for the fourth term on the right end of (\ref{wucha20}), we deduce
\begin{align}\label{wucha23}
\begin{split}
&\frac{2\Delta tk_{f}\gamma}{h}\int_{\uppercase\expandafter{\romannumeral1}}(\theta_{p,j}^{h,n+1}-\theta_{p,j}^{h,n})\varepsilon_{f,j}^{n+1}dl\\
=&\frac{2\Delta tk_{f}\gamma}{h}\int_{\uppercase\expandafter{\romannumeral1}}(\theta_{p,j}^{h,n+1}-\theta_{p,j}(t_{n+1}))\varepsilon_{f,j}^{n+1}dl
+\frac{2\Delta tk_{f}\gamma}{h}\int_{\uppercase\expandafter{\romannumeral1}}(\theta_{p,j}(t_{n})-\theta_{p,j}^{h,n})\varepsilon_{f,j}^{n+1}dl\\
&+\frac{2\Delta tk_{f}\gamma}{h}\int_{\uppercase\expandafter{\romannumeral1}}(\theta_{p,j}(t_{n+1})-\theta_{p,j}(t_{n}))\varepsilon_{f,j}^{n+1}dl\\
=&-\frac{2\Delta tk_{f}\gamma}{h}\int_{\uppercase\expandafter{\romannumeral1}}(\varepsilon_{p,j}^{n+1}-\varepsilon_{p,j}^{n})\varepsilon_{f,j}^{n+1}dl
-\frac{2\Delta tk_{f}\gamma}{h}\int_{\uppercase\expandafter{\romannumeral1}}(\tau_{p,j}^{n+1}-\tau_{p,j}^{n})\varepsilon_{f,j}^{n+1}dl\\
&+\frac{2\Delta tk_{f}\gamma}{h}\int_{\uppercase\expandafter{\romannumeral1}}(\theta_{p,j}(t_{n+1})-\theta_{p,j}(t_{n}))\varepsilon_{f,j}^{n+1}dl\\
\leq&-\frac{\Delta tk_{f}\gamma}{h}\big[\parallel\varepsilon_{p,j}^{n+1}\parallel_{L^{2}(\uppercase\expandafter{\romannumeral1})}^{2}
-\parallel\varepsilon_{p,j}^{n}\parallel_{L^{2}(\uppercase\expandafter{\romannumeral1})}^{2}
-\parallel\varepsilon_{f,j}^{n+1}-\varepsilon_{p,j}^{n+1}\parallel_{L^{2}(\uppercase\expandafter{\romannumeral1})}^{2}\big]
+\frac{2\Delta tk_{f}\gamma}{h}\int_{\uppercase\expandafter{\romannumeral1}}\tau_{p,j}^{n+1}\varepsilon_{f,j}^{n+1}dl\\
&+\frac{2\Delta tk_{f}\gamma}{h}\int_{\uppercase\expandafter{\romannumeral1}}\tau_{p,j}^{n}\varepsilon_{f,j}^{n+1}dl
+\frac{2\Delta tk_{f}\gamma}{h}\int_{\uppercase\expandafter{\romannumeral1}}(\theta_{p,j}(t_{n+1})-\theta_{p,j}(t_{n}))\varepsilon_{f,j}^{n+1}dl\\
\leq&-\frac{\Delta tk_{f}\gamma}{h}\big[\parallel\varepsilon_{p,j}^{n+1}\parallel_{L^{2}(\uppercase\expandafter{\romannumeral1})}^{2}
-\parallel\varepsilon_{p,j}^{n}\parallel_{L^{2}(\uppercase\expandafter{\romannumeral1})}^{2}
-\parallel\varepsilon_{f,j}^{n+1}-\varepsilon_{p,j}^{n+1}\parallel_{L^{2}(\uppercase\expandafter{\romannumeral1})}^{2}\big]\\
&+\frac{3\gamma^{3}k_{f}\tilde{C}_{T}C_{t}\tilde{C}_{t}\tilde{C}_{inv}Ch^{2k}}{4C_{inv}}\big(\parallel\theta_{p,j}(t_{n+1})\parallel_{H^{k+1}(\Omega_{p})}^{2}
+\parallel\theta_{p,j}(t_{n})\parallel_{H^{k+1}(\Omega_{p})}^{2}\big)\\
&+\frac{\Delta t^{4}\gamma^{3}k_{f}\tilde{C}_{T}C_{t}\tilde{C}_{t}}{4h^{2}C_{inv}}\int_{t_{n}}^{t_{n+1}}\parallel\nabla\theta_{p,j,tt}\parallel_{L^{2}(\Omega_{p})}^{2}dt
+\frac{3\Delta t^{3}\gamma^{3}k_{f}\tilde{C}_{T}C_{t}\tilde{C}_{t}}{4h^{2}C_{inv}}\parallel\theta_{p,j,t}(t_{n+1})\parallel_{H^{1}(\Omega_{p})}^{2}\\
&+\frac{4\Delta t k_{f}C_{T}C_{inv}}{\gamma}\parallel\nabla\varepsilon_{f,j}^{n+1}\parallel_{L^{2}(\Omega_{f})}^{2}.
\end{split}
\end{align}
Inserting (\ref{wucha8}) - (\ref{wucha23}) into the equation (\ref{wucha7}), summing up from $n=0$ to $N-1$ and using Gronwall's lemma yields
\begin{align}\label{wucha24}
\begin{split}
&\parallel b_{f,j}^{N}\parallel_{L^{2}(\Omega_{f})}^{2}
+\sum\limits_{n=0}^{N-1}\parallel b_{f,j}^{n+1}-b_{f,j}^{n}\parallel_{L^{2}(\Omega_{f})}^{2}
+\frac{C_{a}k_{min}}{L^{2}}\big[\parallel b_{p,j}^{N}\parallel_{L^{2}(\Omega_{p})}^{2}
+\sum\limits_{n=0}^{N-1}\parallel b_{p,j}^{n+1}-b_{p,j}^{n}\parallel_{L^{2}(\Omega_{p})}^{2}\big]\\
&+\parallel \varepsilon_{f,j}^{N}\parallel_{L^{2}(\Omega_{f})}^{2}
+\sum\limits_{n=0}^{N-1}\parallel \varepsilon_{f,j}^{n+1}-\varepsilon_{f,j}^{n}\parallel_{L^{2}(\Omega_{f})}^{2}
+\parallel \varepsilon_{p,j}^{N}\parallel_{L^{2}(\Omega_{p})}^{2}
+\sum\limits_{n=0}^{N-1}\parallel \varepsilon_{p,j}^{n+1}-\varepsilon_{p,j}^{n}\parallel_{L^{2}(\Omega_{p})}^{2}\\
&+\frac{1}{2}Pr\Delta t\sum\limits_{n=0}^{N-1}\parallel\nabla b_{f,j}^{n+1}\parallel_{L^{2}(\Omega_{f})}^{2}
+Pr\Delta t\sum\limits_{n=0}^{N-1}\parallel b_{p,j}^{n+1}\parallel_{L^{2}(\Omega_{p})}^{2}
+k_{f}\Delta t\bigg(1-\frac{8C_{T}C_{inv}}{\gamma}\bigg)\sum\limits_{n=0}^{N-1}\parallel\nabla \varepsilon_{f,j}^{n+1}\parallel_{L^{2}(\Omega_{f})}^{2}\\
&+k_{p}\Delta t\sum\limits_{n=0}^{N-1}\parallel\nabla \varepsilon_{p,j}^{n+1}\parallel_{L^{2}(\Omega_{p})}^{2}
+\frac{\Delta tk_{f}\gamma}{h}\parallel \varepsilon_{p,j}^{N}\parallel_{L^{2}(\uppercase\expandafter{\romannumeral1})}^{2}\\
\leq&C(\Delta t^{2}+h^{2k}),
\end{split}
\end{align}
where $\gamma>8C_{T}C_{inv}$. Finally, applying the expectation on (\ref{wucha24}), we have the error estimate (\ref{wucha}).  $\qedsymbol$
\section{Numerical experiments}
~\par In this section, we present numerical simulation of the fully discrete FEM-MCM method for the 2D/3D stochastic closed-loop geothermal system. The first example with an analytic solution is provided to show the convergence of finite element solutions when the hydraulic conduction tensor takes different values. The second numerical example uses the FEM-MCM method to solve the 2D stochastic closed-loop geothermal system, and the optimal convergence order is obtained. In addition, the influence of penalty parameter $\gamma$ on convergence is demonstrated. The third numerical test is conducted to investigate a test on 3D natural convection in a cubical cavity with left hand side heating. The applicability of the proposed model and numerical method in 3D stochastic closed-loop geothermal systems is demonstrated by a simplified device in the fourth and fifth examples respectively.\par
In the first two experiments, we apply the MINI elements $(P1b-P1)$ for the Navier-Stokes region, piecewise constant elements $P0$ for Darcy pressure, and Brezzi-Douglas-Marini $(BDM1)$ finite elements for Darcy velocity \cite{W. J. Layton-2002,F. Brezzi-1985}. For the temperatures in the two regions, we apply the linear Lagrangian elements $P1$. In the following table, the error is defined as $\bar{e}^{h,n}_{\Psi_{\kappa}}= \Psi^{h,n}_{\kappa}-\Psi_{\kappa}(t_{n})$, where $\Psi = \mathbf{u}$ or $\theta$ and $\kappa = f$ or $p$.
For the latter three experiments, we use the MINI elements $(P1b3d-P13d)$ for the Navier-Stokes region, piecewise constant elements $P03d$ for Darcy pressure, and Raviart-Thomas $(RT03d)$ finite elements for Darcy velocity. For the temperatures in the two regions, we use the linear Lagrangian elements $P13d$.\par
\subsection{Two-dimensional regimes}
\subsubsection{Convergence test when K takes different values}
~\par In this section, in order to test the convergence, we consider the problem on $\Omega=(0,1)\times(0,2)$, where $\Omega_{f}=(0,1)\times(1,2)$, and $\Omega_{p}=(0,1)\times(0,1)$. We take $a=1.0,\ Pr=1.0,\ Ra=1.0,\ C_{a}=1.0,\ k_{f}=k_{p}=1.0,\  \gamma=10^{5},\ L=1.0,\ T=0.5,\ \Delta t=0.001$, and
\begin{align}
\begin{split}
&K=K_{j}=
\left[ \begin{matrix}
	k_{11}^{j}&		0         \\
	0         &		k_{22}^{j}\\
\end{matrix} \right] ,
\qquad j=1,...,J,\nonumber
\end{split}
\end{align}
where $K$ is a random hydraulic conductivity tensor and $K_{j}$ is one of the samples of $K$.

Moreover, the exact solution is

\begin{align*}
\begin{split}
\begin{cases}
&\mathbf{u}_{f}=\Big(10x^{2}(x-1)^{2}y(y-1)(2y-1)cos(t)k^{j}_{11},-10x(x-1)(2x-1)y^{2}(y-1)^{2}cos(t)k^{j}_{22}\Big),\\
&p_{f}=10(2x-1)(2y-1)cos(t),\\
&\mathbf{u}_{p}=\Big([2\pi sin^{2}(\pi x)sin(\pi y)cos(\pi y)]cos(t),[-2\pi sin(\pi x)sin^{2}(\pi y)cos(\pi x)]cos(t)\Big),\\
&\phi_{p}=cos(\pi x)cos(\pi y)cos(t),\\
&\theta_{f}=ax(1-x)(1-y)e^{-t},\\
&\theta_{f}=ax(1-x)(y-y^{2})e^{-t}.
\end{cases}
\end{split}
\end{align*}

The physical parameters hydraulic conductivity tensor was chosen as $k_{11}^{1}=k_{22}^{1}=2.21,\ k_{11}^{2}=k_{22}^{2}=4.21$, and $\ k_{11}^{3}=k_{22}^{3}=6.21$. We have reported in Tables 1-3 the errors and convergence rates of the decoupled method for the stochastic closed-loop geothermal system. We can achieve the optimal convergence order for both the $L^{2}$-norm and $H^{1}$-norm, which is consistent with the theoretical results.

\begin{table}[h]
\centering
\caption{Errors and convergence rates of the algorithm: $k_{11}=k_{22}=2.21$}
\begin{tabular}{cccccccccc}
\toprule
$h$ & $\parallel \bar{e}_{\mathbf{u}_{f}}^{h,n}\parallel_{0}$ & $rate$
& $ \parallel \bar{e}_{\theta_{f}}^{h,n}\parallel_{0}$ & $rate$
& $\parallel \bar{e}_{\mathbf{u}_{p}}^{h,n}\parallel_{0}$ & $rate$
& $\parallel \bar{e}_{\theta_{p}}^{h,n}\parallel_{0}$ & $rate$ &\\
\midrule
$\frac{1}{4}$  & $0.968279$    & $-$           & $0.00463045$    & $-$           & $0.261987$   & $-$         & $0.00362939$   & $-$       & \\
$\frac{1}{8}$  & $0.236474$    & $2.03374$     & $0.00100742$    & $2.20048$     & $0.0520901$  & $2.33041$   & $0.000679279$  & $2.41765$ & \\
$\frac{1}{16}$ & $0.0580309$   & $2.02679$     & $0.000219494$   & $2.19842$     & $0.0120885$  & $2.10738$   & $0.000147805$  & $2.20031$ & \\
$\frac{1}{32}$ & $0.0146095$   & $1.98992$     & $5.42306e-05$   & $2.017$       & $0.00313974$ & $1.94492$   & $3.49639e-05$  & $2.07975$ & \\
\bottomrule
\toprule
$h$ & $\parallel \bar{e}_{\mathbf{u}_{f}}^{h,n}\parallel_{1}$ & $rate$
& $ \parallel \bar{e}_{\theta_{f}}^{h,n}\parallel_{1}$ & $rate$
& $\parallel \bar{e}_{\mathbf{u}_{p}}^{h,n}\parallel_{H(div)}$ & $rate$
& $\parallel \bar{e}_{\theta_{p}}^{h,n}\parallel_{1}$ & $rate$ & \\
\midrule
$\frac{1}{4}$  & $14.7899$  & $-$           & $0.062879$   & $-$           & $3.35363$  & $-$          & $0.0367264$  & $-$       & \\
$\frac{1}{8}$  & $6.44033$  & $1.19941$     & $0.03087$    & $1.02637$     & $1.87647$  & $1.11391$    & $0.0155855$  & $1.23662$ & \\
$\frac{1}{16}$ & $2.59281$  & $1.31262$     & $0.0148033$  & $1.06029$     & $0.881757$ & $1.04109$    & $0.00730445$ & $1.09336$ & \\
$\frac{1}{32}$ & $1.22146$  & $1.08591$     & $0.00751796$ & $0.977505$    & $0.444293$ & $0.975316$   & $0.0035737$  & $1.03136$ & \\
\bottomrule
\end{tabular}
\end{table}

\begin{table}[h]
\centering
\caption{Errors and convergence rates of the algorithm: $k_{11}=k_{22}=4.21$}
\begin{tabular}{cccccccccc}
\toprule
$h$ & $\parallel \bar{e}_{\mathbf{u}_{f}}^{h,n}\parallel_{0}$ & $rate$
& $ \parallel \bar{e}_{\theta_{f}}^{h,n}\parallel_{0}$ & $rate$
& $\parallel \bar{e}_{\mathbf{u}_{p}}^{h,n}\parallel_{0}$ & $rate$
& $\parallel \bar{e}_{\theta_{p}}^{h,n}\parallel_{0}$ & $rate$ &\\
\midrule
$\frac{1}{4}$  & $1.84732$     & $-$           & $0.00470008$     & $-$           & $0.264952$   & $-$         & $0.0036294$   & $-$       &\\
$\frac{1}{8}$  & $0.451223$    & $2.03352$     & $0.000969245$    & $2.27775$     & $0.0524375$  & $2.33706$   & $0.000679279$ & $2.41765$  &\\
$\frac{1}{16}$ & $0.110871$    & $2.02496$     & $0.000207871$    & $2.22117$     & $0.0121527$  & $2.10932$   & $0.000147805$ & $2.20031$ &\\
$\frac{1}{32}$ & $0.0279204$   & $1.98949$     & $5.20351e-05$    & $1.99813$     & $0.00315756$ & $1.94439$   & $3.49639e-05$ & $2.07975$ &\\
\bottomrule
\toprule
$h$ & $\parallel \bar{e}_{\mathbf{u}_{f}}^{h,n}\parallel_{1}$ & $rate$
& $ \parallel \bar{e}_{\theta_{f}}^{h,n}\parallel_{1}$ & $rate$
& $\parallel \bar{e}_{\mathbf{u}_{p}}^{h,n}\parallel_{H(div)}$ & $rate$
& $\parallel \bar{e}_{\theta_{p}}^{h,n}\parallel_{1}$ & $rate$ & \\
\midrule
$\frac{1}{4}$  & $28.2506$  & $-$           & $0.0633716$   & $-$           & $3.38017$  & $-$          & $0.0367264$  & $-$       & \\
$\frac{1}{8}$  & $12.2873$  & $1.20111$     & $0.0309526$   & $1.03378$     & $1.88396$  & $1.11856$    & $0.0155855$  & $1.23662$ & \\
$\frac{1}{16}$ & $4.94214$  & $1.31396$     & $0.0148113$   & $1.06336$     & $0.884263$ & $1.04183$    & $0.00730445$ & $1.09336$ & \\
$\frac{1}{32}$ & $2.32732$  & $1.08647$     & $0.00751879$  & $0.978129$    & $0.445271$ & $0.975695$   & $0.0035737$  & $1.03136$ & \\
\bottomrule
\end{tabular}
\end{table}
\begin{table}[h]
\centering
\caption{Errors and convergence rates of the algorithm: $k_{11}=k_{22}=6.21$}
\begin{tabular}{cccccccccc}
\toprule
$h$ & $\parallel \bar{e}_{\mathbf{u}_{f}}^{h,n}\parallel_{0}$ & $rate$
& $ \parallel \bar{e}_{\theta_{f}}^{h,n}\parallel_{0}$ & $rate$
& $\parallel \bar{e}_{\mathbf{u}_{p}}^{h,n}\parallel_{0}$ & $rate$
& $\parallel \bar{e}_{\theta_{p}}^{h,n}\parallel_{0}$ & $rate$ &\\
\midrule
$\frac{1}{4}$  & $2.73063$     & $-$           & $0.00481102$     & $-$           & $0.266161$   & $-$         & $0.0036294$     & $-$      &\\
$\frac{1}{8}$  & $0.667075$    & $2.03331$     & $0.000944391$    & $2.34889$     & $0.0525797$  & $2.33973$   & $0.000679279$   & $2.41766$&\\
$\frac{1}{16}$ & $0.164148$    & $2.02285$     & $0.000199411$    & $2.24364$     & $0.012179$   & $2.11011$   & $0.000147805$   & $2.20031$&\\
$\frac{1}{32}$ & $0.0413512$   & $1.98899$     & $5.06762e-05$    & $1.97636$     & $0.00316478$ & $1.94421$   & $3.49639e-05$   & $2.07975$&\\
\bottomrule
\toprule
$h$ & $\parallel \bar{e}_{\mathbf{u}_{f}}^{h,n}\parallel_{1}$ & $rate$
& $ \parallel \bar{e}_{\theta_{f}}^{h,n}\parallel_{1}$ & $rate$
& $\parallel \bar{e}_{\mathbf{u}_{p}}^{h,n}\parallel_{H(div)}$ & $rate$
& $\parallel \bar{e}_{\theta_{p}}^{h,n}\parallel_{1}$ & $rate$ & \\
\midrule
$\frac{1}{4}$  & $41.8501$  & $-$           & $0.0643057$   & $-$           & $3.39062$  & $-$          & $0.0367264$  & $-$       & \\
$\frac{1}{8}$  & $18.1688$  & $1.20377$     & $0.0310468$   & $1.0505$      & $1.8869$   & $1.12043$    & $0.0155855$  & $1.23662$ & \\
$\frac{1}{16}$ & $7.2968$   & $1.31613$     & $0.0148211$   & $1.06679$     & $0.885262$ & $1.04213$    & $0.00730445$ & $1.09336$ & \\
$\frac{1}{32}$ & $3.43403$  & $1.08736$     & $0.00751983$  & $0.978885$    & $0.445665$ & $0.975846$   & $0.0035737$  & $1.03136$ & \\
\bottomrule
\end{tabular}
\end{table}
\subsubsection{Optimal convergence order with random hydraulic conductivity tensor}
~\par In this example, we apply the FEM-MCM method for solving a stochastic closed-loop geothermal system with a random hydraulic conduction tensor. Here we choose the following random hydraulic conductivity tensor
\begin{align}
\begin{split}
&K(x,y,\mathbf{\lambda})=K_{j}(x,y,\mathbf{\lambda})
\left[ \begin{matrix}
	k_{11}^{j}(x,y,\mathbf{\lambda})&		0                              \\
	0                               &		k_{22}^{j}(x,y,\mathbf{\lambda})\\
\end{matrix} \right] ,
\qquad j=1,...,J,\nonumber
\end{split}
\end{align}
and
\begin{align*}
\begin{split}
&k_{11}^{j}(x,y,\mathbf{\lambda})=k_{22}^{j}(x,y,\mathbf{\lambda})=3+\sigma(\lambda_{1}+\lambda_{2}).
\end{split}
\end{align*}

The exact solution is set by
\begin{align*}
\begin{split}
\begin{cases}
&\mathbf{u}_{f}=\Big(\frac{1}{3+\lambda_{1}+\lambda_{2}}10x^{2}(x-1)^{2}y(y-1)(2y-1)cos(t),-\frac{1}{3+\lambda_{1}+\lambda_{2}}10x(x-1)(2x-1)y^{2}(y-1)^{2}cos(t)\Big),\\
&\mathbf{u}_{p}=\Big([2\pi sin^{2}(\pi x)sin(\pi y)cos(\pi y)]cos(t),[-2\pi sin(\pi x)sin^{2}(\pi y)cos(\pi x)]cos(t)\Big),\\
&p_{f}=10(2x-1)(2y-1)cos(t),\ \phi_{p}=cos(\pi x)cos(\pi y)cos(t),\\
&\theta_{f}=ax(1-x)(1-y)e^{-t},\ \theta_{p}=ax(1-x)(y-y^{2})e^{-t},
\end{cases}
\end{split}
\end{align*}
where $(x, y)\in \Omega=(0, 1)\times(0,2),\ \Omega_{f}=(0,1)\times(1,2),\ \Omega_{p}=(0,1)\times(0,1),\ \sigma=0.1$, and $\lambda_{1},\ \lambda_{2}$ are independent and identically uniform distributed in the interval $[-1,1]$. The initial conditions, boundary conditions, and forcing terms are chosen to match the exact solution. We set $T=0.5$ and the number of samples $J=64$. In addition, we take $a=1.0,\ Pr=1.0,\ Ra=1.0,\ C_{a}=1.0,\ k_{f}=k_{p}=1.0,\  \gamma=10^{5}$, and $L=1.0$.

The MCM with $J=64$ realizations is used to compute the expectation. For any $1\leq n \leq N$, we use the following numerical integration formulas to represent the approximate errors
\begin{align*}
\begin{split}
&\parallel\bar{e}^{h,n}_{\Psi_{\kappa}}\parallel_{a}
:=\bigg(E\big[\parallel\Psi^{h,n}_{\kappa}-\Psi_{\kappa}(t_{n})\parallel_{a}^{2}\big]\bigg)^{\frac{1}{2}}
\approx\bigg(\frac{1}{J}\sum\limits_{j=1}^{J}
\parallel\Psi^{h,n}_{\kappa,j}-\Psi_{\kappa,j}(t_{n})\parallel_{a}^{2}\bigg)^{\frac{1}{2}},
\end{split}
\end{align*}
where $a=0$ or 1.

\begin{table}[!ht]
\centering
\caption{Errors and convergence rates with random hydraulic conductivity tensor}
\begin{tabular}{cccccccccc}
\toprule
$h$ & $\parallel \bar{e}_{\mathbf{u}_{f}}^{h,n}\parallel_{0}$ & $\rho_{\mathbf{u}_{f}}$
& $ \parallel \bar{e}_{\theta_{f}}^{h,n}\parallel_{0}$ & $\rho_{\theta_{f}}$
& $\parallel \bar{e}_{\mathbf{u}_{p}}^{h,n}\parallel_{0}$ & $\rho_{\mathbf{u}_{p}}$
& $\parallel \bar{e}_{\theta_{p}}^{h,n}\parallel_{0}$ & $\rho_{\theta_{p}}$ & \\
\midrule
$\frac{1}{4}$  & $0.152341403$    & $-$          & $0.00460378$    & $-$         & $0.26377$     & $-$         & $0.0036294$     & $-$      &\\
$\frac{1}{8}$  & $0.037243306$    & $2.03226$    & $0.00103291$    & $2.15611$   & $0.0522988$   & $2.33443$   & $0.000679279$   & $2.41765$&\\
$\frac{1}{16}$ & $0.009136531$    & $2.02726$    & $0.000228093$   & $2.17902$   & $0.012127$    & $2.10855$   & $0.000147805$   & $2.20031$&\\
$\frac{1}{32}$ & $0.002298631$    & $1.99087$    & $5.61951e-05$   & $2.02111$   & $0.00315048$  & $1.94459$   & $3.49639e-05$   & $2.07975$&\\
\bottomrule
\toprule
$h$ & $\parallel \bar{e}_{\mathbf{u}_{f}}^{h,n}\parallel_{1}$ & $\rho_{\mathbf{u}_{f}}$
& $ \parallel \bar{e}_{\theta_{f}}^{h,n}\parallel_{1}$ & $\rho_{\theta_{f}}$
& $\parallel \bar{e}_{\mathbf{u}_{p}}^{h,n}\parallel_{H(div)}$ & $\rho_{\mathbf{u}_{p}}$
& $\parallel \bar{e}_{\theta_{p}}^{h,n}\parallel_{1}$ & $\rho_{\theta_{p}}$ & \\
\midrule
$\frac{1}{4}$  & $2.33270375$   & $-$           & $0.062712$    & $-$          & $3.36976$   & $-$          & $0.0367264$   & $-$       & \\
$\frac{1}{8}$  & $1.015652766$  & $1.19959$     & $0.0308302$   & $1.0244$     & $1.88102$   & $1.11672$    & $0.0155855$   & $1.23662$ & \\
$\frac{1}{16}$ & $0.408966672$  & $1.31235$     & $0.0147995$   & $1.05879$    & $0.883275$  & $1.04153$    & $0.00730445$  & $1.09336$ & \\
$\frac{1}{32}$ & $0.19268$      & $1.08578$     & $0.00751758$  & $0.97721$    & $0.444883$  & $0.975545$   & $0.0035737$   & $1.03136$ & \\
\bottomrule
\end{tabular}
\end{table}

In order to calculate the convergence order of mesh size $h$, we define
\begin{align*}
\begin{split}
&\rho_{\Psi_{\kappa}}=\frac{log(e_{\Psi_{\kappa}}^{h_{1},n}/e_{\Psi_{\kappa}}^{h_{2},n})}{log(h_{1}/h_{2})}.
\end{split}
\end{align*}

We fix $\Delta t=0.001$, and use different spatial mesh size $h=\frac{1}{4},\frac{1}{8},\frac{1}{16}$, and $\frac{1}{32}$. We can observe from Table 5 that the decoupled method of the stochastic closed-loop geothermal system achieves optimal convergence rates in both the $L^{2}$-norm and the $H^{1}$-norm. Obviously, we get the same optimal convergence order as the deterministic closed-loop geothermal system.

\begin{table}[!ht]
\centering
\caption{Errors and convergence rates with random hydraulic conductivity tensor}
\begin{tabular}{cccccccccc}
\toprule
$\Delta t$ & $\parallel \bar{e}_{\mathbf{u}_{f}}^{h,n}\parallel_{0}$ & $\beta_{\mathbf{u}_{f}}$
& $ \parallel \bar{e}_{\theta_{f}}^{h,n}\parallel_{0}$ & $\beta_{\theta_{f}}$
& $\parallel \bar{e}_{\mathbf{u}_{p}}^{h,n}\parallel_{0}$ & $\beta_{\mathbf{u}_{p}}$
& $\parallel \bar{e}_{\theta_{p}}^{h,n}\parallel_{0}$ & $\beta_{\theta_{p}}$ &\\
\midrule
$\frac{1}{20}$  & $4.31087e-05$   & $-$        & $7.05199e-05$   & $-$         & $1.64515e-05$   & $-$         & $2.8029e-05$    & $-$      &\\
$\frac{1}{40}$  & $2.27774e-05$   & $1.89261$  & $3.46326e-05$   & $2.03623$   & $8.34363e-06$   & $1.97174$   & $1.37453e-05$   & $2.03916$&\\
$\frac{1}{80}$  & $1.16927e-05$   & $1.94801$  & $1.71185e-05$   & $2.02311$   & $4.20186e-06$   & $1.9857$    & $6.7897e-06$    & $2.02444$&\\
$\frac{1}{160}$ & $5.9222e-06$    & $1.97438$  & $8.50615e-06$   & $2.01249$   & $2.10851e-06$   & $1.99281$   & $3.37284e-06$   & $2.01305$&\\
\bottomrule
\toprule
$\Delta t$ & $\parallel \bar{e}_{\mathbf{u}_{f}}^{h,n}\parallel_{1}$ & $\beta_{\mathbf{u}_{f}}$
& $ \parallel \bar{e}_{\theta_{f}}^{h,n}\parallel_{1}$ & $\beta_{\theta_{f}}$
& $\parallel \bar{e}_{\mathbf{u}_{p}}^{h,n}\parallel_{H(div)}$ & $\beta_{\mathbf{u}_{p}}$
& $\parallel \bar{e}_{\theta_{p}}^{h,n}\parallel_{1}$ & $\beta_{\theta_{p}}$ & \\
\midrule
$\frac{1}{20}$  & $0.0005228583$ & $-$         & $0.000326623$   & $-$           & $0.00363762$   & $-$          & $0.000124672$  & $-$       & \\
$\frac{1}{40}$  & $0.0002773511$ & $1.88512$   & $0.000160381$   & $2.03654$     & $0.00184399$   & $1.97269$    & $6.1139e-05$   & $2.03917$ & \\
$\frac{1}{80}$  & $0.0001426276$ & $1.94458$   & $7.92846e-05$   & $2.02285$     & $0.000928418$  & $1.98616$    & $3.02005e-05$  & $2.02444$ & \\
$\frac{1}{160}$ & $7.22978e-05$  & $1.97278$   & $3.94005e-05$   & $2.01228$     & $0.000465831$  & $1.99304$    & $1.50023e-05$  & $2.01305$ & \\
\bottomrule
\end{tabular}
\end{table}

Here we use the same method as in \cite{M. Mu-2009} to calculate the convergence order with respect to the time step $\Delta t$. For this we define the parameter
\begin{align}
\begin{split}
&\beta_{v}=\frac{\parallel v^{h,\Delta t}-v^{h,\frac{1}{2}\Delta t}\parallel_{0}}{\parallel v^{h,\frac{1}{2}\Delta t}-v^{h,\frac{1}{4}\Delta t}\parallel_{0}},\nonumber
\end{split}
\end{align}
where $v=\mathbf{u}_{f},\mathbf{u}_{p},\theta_{f}$, and $\theta_{p}$.

In Table 5, we fix $h=\frac{1}{32}$ and then set the time step $\Delta t=\frac{1}{20}, \frac{1}{40}, \frac{1}{80}$, and $\frac{1}{160}$ for numerical calculation. Table 5 shows $L^{2}$-norm and $H^{1}$-norm approximate errors of velocity and temperature of the decoupled scheme. In addition, the data in Table 5 illustrates that the decoupled method of the stochastic closed-loop geothermal system achieves optimal convergence order in the $L^{2}$-norm and the $H^{1}$-norm. These results completely agree with the theoretical results.

Below we examine the effect of penalty parameters $\gamma$ on convergence. We fix $h=\frac{1}{8}$, $\Delta t=0.001$, and show the errors of velocity and temperature under penalty parameters $\gamma=0, 10^{-3}, 1, 10^{3}$, and $10^{5}$ in Figure 2 and Figure 3. We can observe that there is almost no difference in velocity error when taking different values of the penalty parameters. However, as Figure 3 shows, the temperature error is smaller with the increase of parameters $\gamma$. That is, a larger penalty parameter provides better convergence performance for $\theta_{f}$ and $\theta_{p}$.
\begin{figure}[!ht]
\centering
\begin{minipage}[t]{0.48\linewidth}
\centering
\includegraphics[height=5.3cm,width=5.3cm]{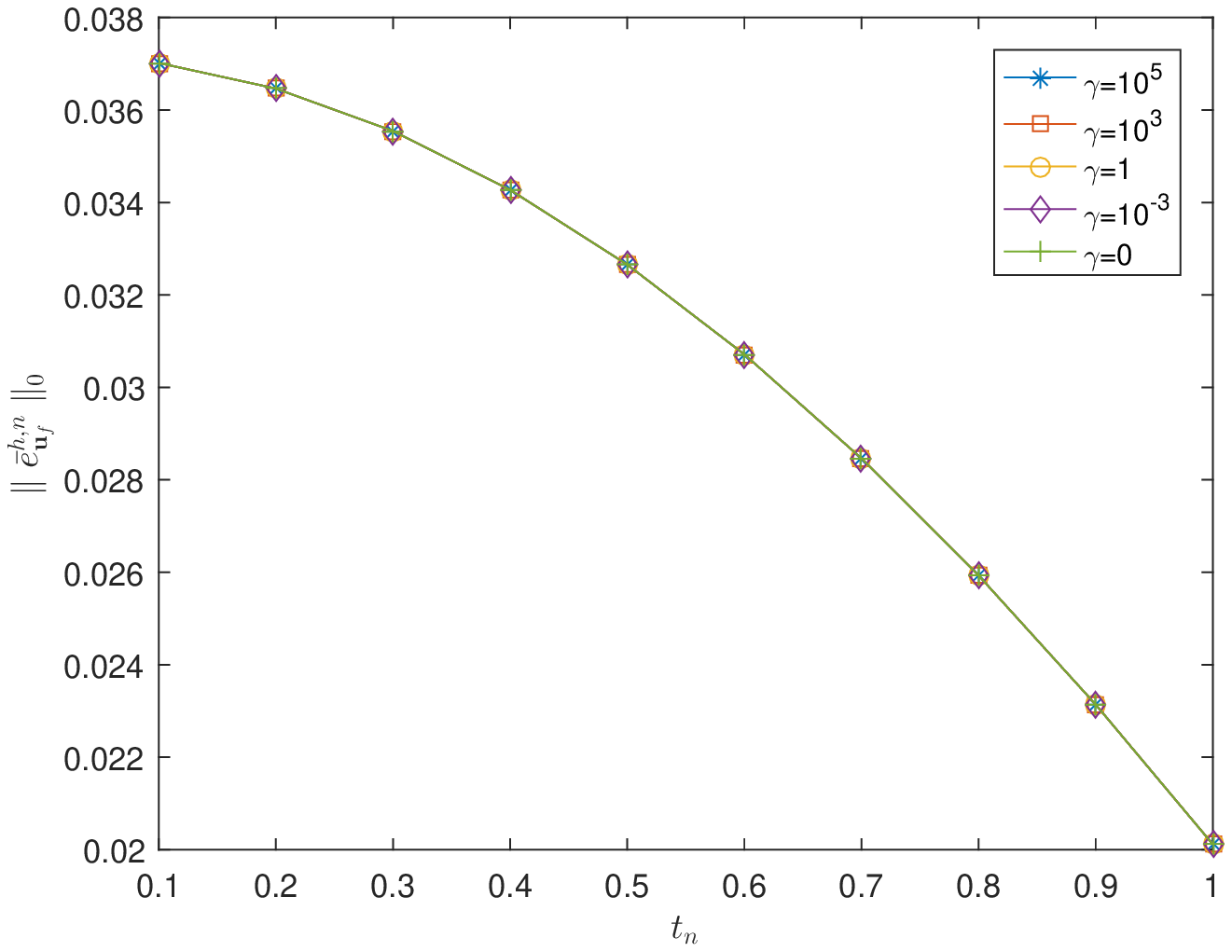}
\end{minipage}
\hfill
\begin{minipage}[t]{0.48\linewidth}
\centering
\includegraphics[height=5.3cm,width=5.3cm]{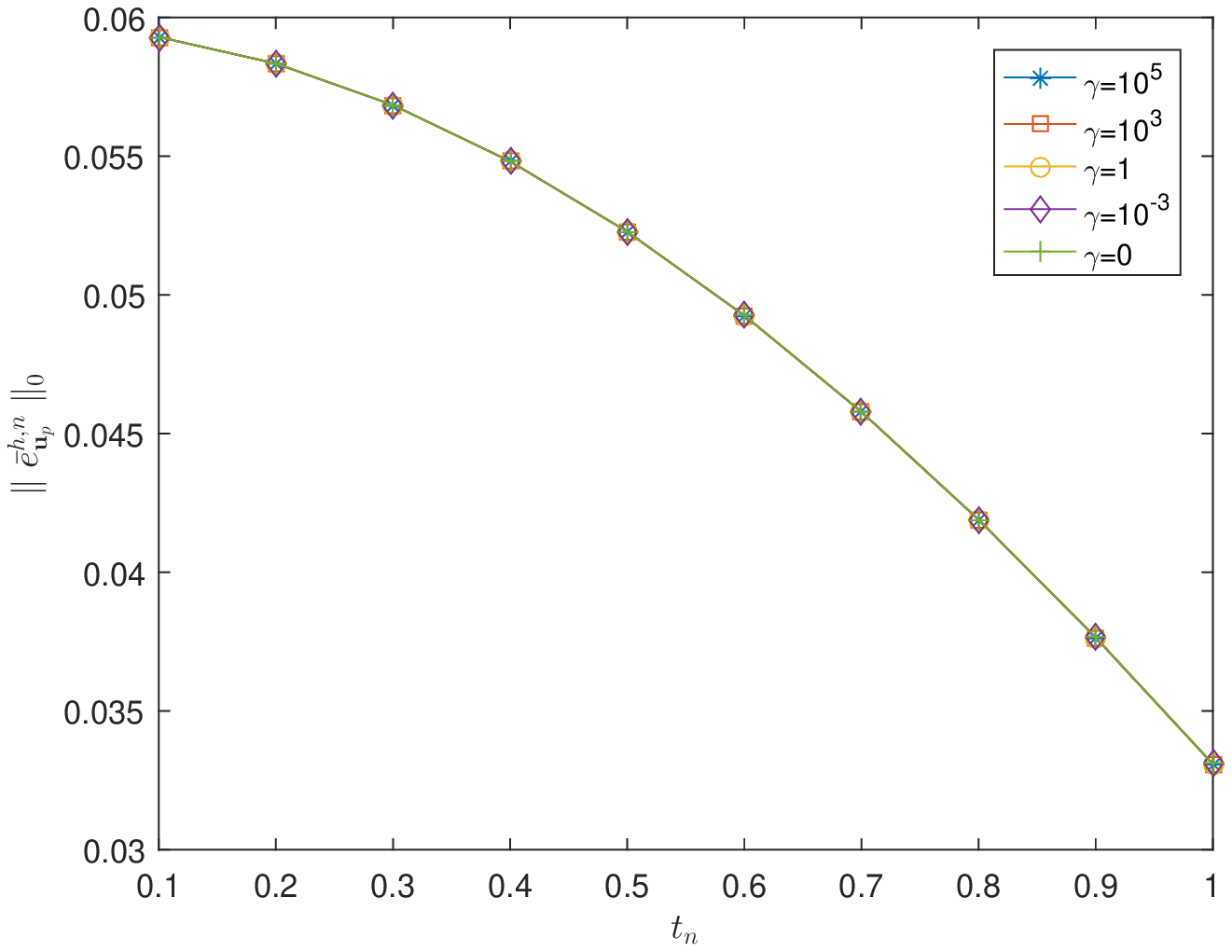}
\end{minipage}
\caption{$\parallel \bar{e}_{\mathbf{u}_{f}}^{h,n}\parallel_{0}$ (left) and $\parallel \bar{e}_{\mathbf{u}_{p}}^{h,n}\parallel_{0}$ (right) with different penalty parameter $\gamma$.}
\end{figure}

\begin{figure}[!ht]
\centering
\begin{minipage}[t]{0.48\linewidth}
\centering
\includegraphics[height=5.3cm,width=5.3cm]{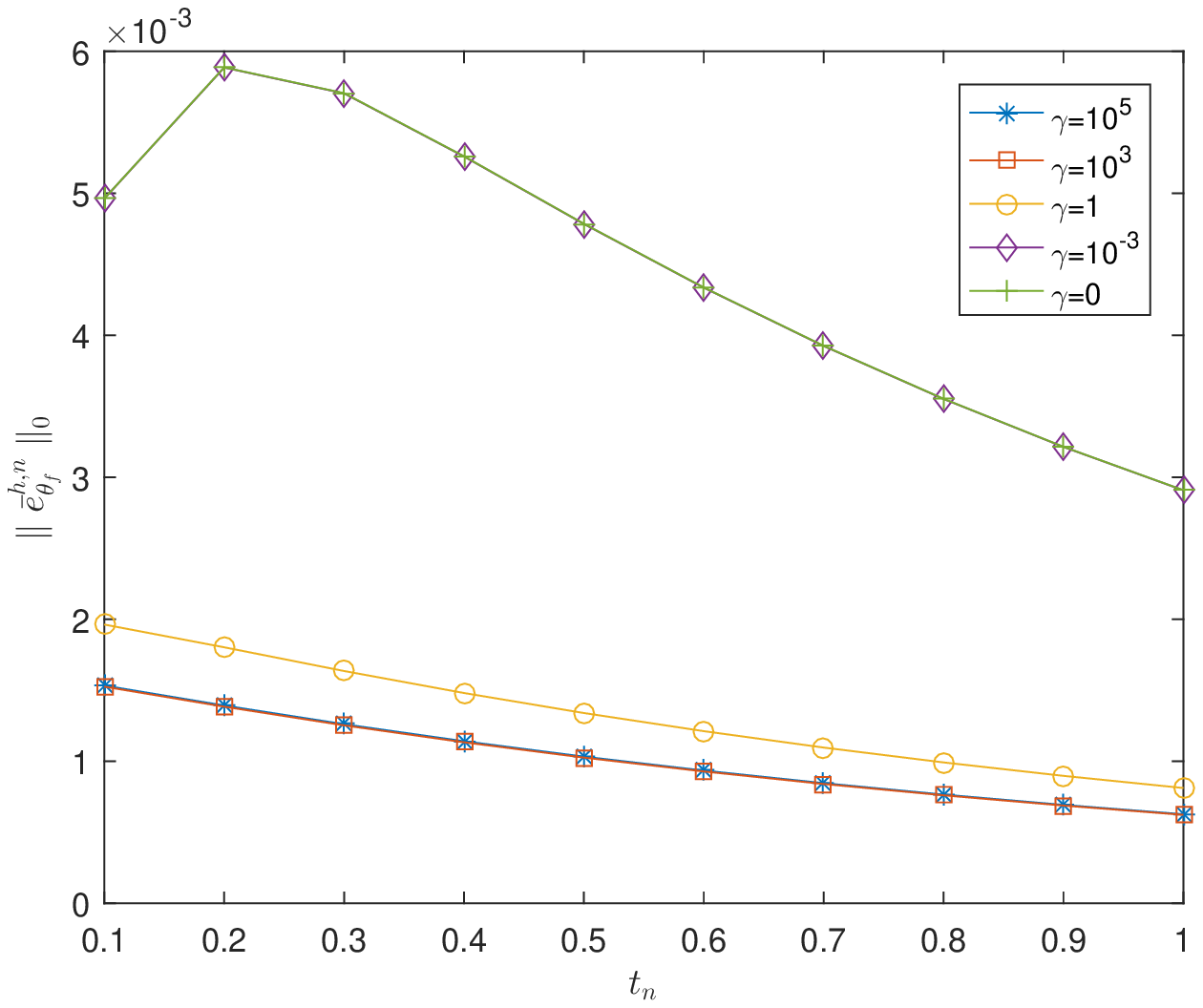}
\end{minipage}
\hfill
\begin{minipage}[t]{0.48\linewidth}
\centering
\includegraphics[height=5.3cm,width=5.3cm]{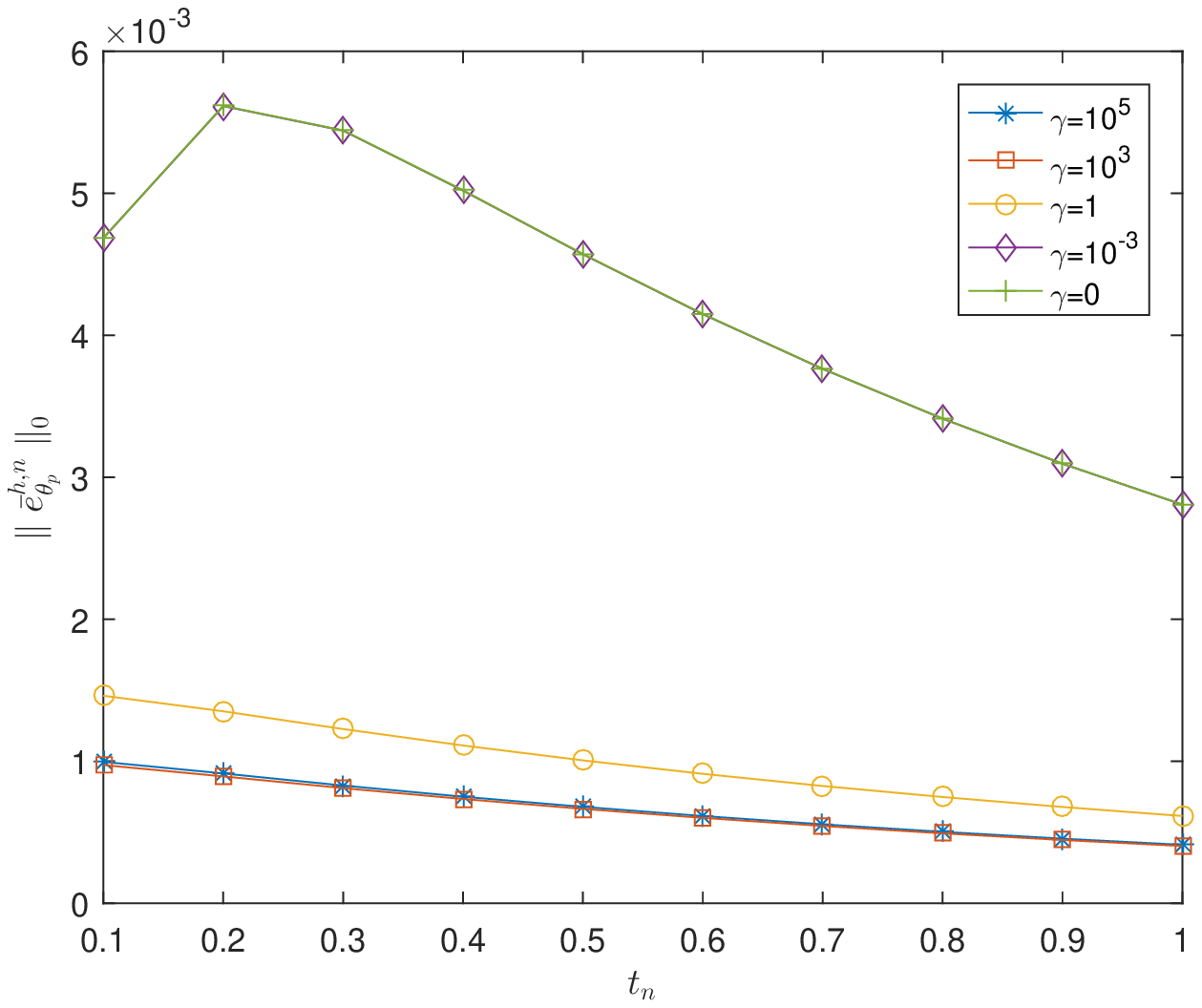}
\end{minipage}
\caption{$\parallel \bar{e}_{\theta_{f}}^{h,n}\parallel_{0}$ (left) and $\parallel \bar{e}_{\theta_{p}}^{h,n}\parallel_{0}$ (right) with different penalty parameter $\gamma$.}
\end{figure}
\subsection{Three-dimensional regimes}
\subsubsection{3D natural convection in a cubical cavity with the left wall heating}
~\par In order to illustrate the effectiveness of the proposed method, we used the benchmark tests proposed in \cite{D C Wan-2001,Z Y Zhang-2013,Z Y Zhang-2014}. We present a 3D natural convection in a cubical cavity with the left wall heating. The domain $\Omega_{f}=(0,1)\times(1,2)\times(0,1)$ and $\Omega_{p}=(0,1)^{3}$ are given. We set the boundary conditions of the computational domain as follows: the left boundary has $\theta_{f}=\theta_{p}=1$ and the right boundary has $\theta_{f}=\theta_{p}=0$. In addition, except for the interface, all the other boundaries of the pipeline region satisfy the no-slip boundary condition, and all the other boundaries of the porous media region satisfy the no-flow boundary condition:
\begin{align*}
\begin{split}
&\mathbf{u}_{f}=0\quad on\ \partial\Omega_{f}\backslash\uppercase\expandafter{\romannumeral1},\qquad
\mathbf{u}_{p}\cdot\mathbf{n}_{p}=0\quad on\ \partial\Omega_{p}\backslash\uppercase\expandafter{\romannumeral1}.
\end{split}
\end{align*}

We choose the parameters $Pr=0.71,\ k_{f}=k_{p}=C_{a}=1.0,\ L=10.0,\ \gamma=1.0,\ T=3.0$, the time step $\Delta t=0.01$, and the mesh size $h=\frac{1}{9}$. In addition, the random hydraulic conductivity tensor is chosen as follows
\begin{align*}
\begin{split}
&K(x,y,z,\mathbf{\lambda})=
\left[ \begin{matrix}
	k_{11}(x,y,z,\mathbf{\lambda})&		0                              &		0                              \\
	0                             &		k_{22}(x,y,z,\mathbf{\lambda}) &		0                              \\
	0                             &		0                              &		k_{33}(x,y,z,\mathbf{\lambda}) \\
\end{matrix} \right] ,
\nonumber
\end{split}
\end{align*}
and
\begin{align*}
\begin{split}
&k_{11}(x,y,z,\mathbf{\lambda})=k_{22}(x,y,z,\mathbf{\lambda})=k_{33}(x,y,z,\mathbf{\lambda})=a_{0}+\sigma\sqrt\omega_{0}Y_{0}(\lambda)+\sum\limits_{i=1}^{n_{f}}\sigma\sqrt\omega_{i}
[Y_{i}(\lambda)cos(i\pi y)+Y_{n_{f}+i}(\lambda)sin(i\pi y)],
\end{split}
\end{align*}
where $\omega_{0}=\frac{\sqrt{\pi L_{c}}}{2},\omega_{i}=\sqrt{\pi}L_{c}e^{-\frac{(i\pi L_{c})^{2}}{4}}$ for $i=1,...,n_{f}$ and $Y_{0},...,Y_{2n_{f}}$ are uncorrelated random variables with zero mean and unit variance. We take the desired physical correlation length $L_{c} = 0.25$ for the random field and $a_{0}= 1, \sigma = 0.15, n_{f} = 3$. We assume the random variables $Y_{0},...,Y_{2n_{f}}$ are independent and uniformly distributed in the interval $[-\sqrt{3}, \sqrt{3}]$.

\begin{figure}[!ht]
\centering
\begin{minipage}[t]{0.48\linewidth}
\centering
\includegraphics[height=6.5cm,width=6.5cm]{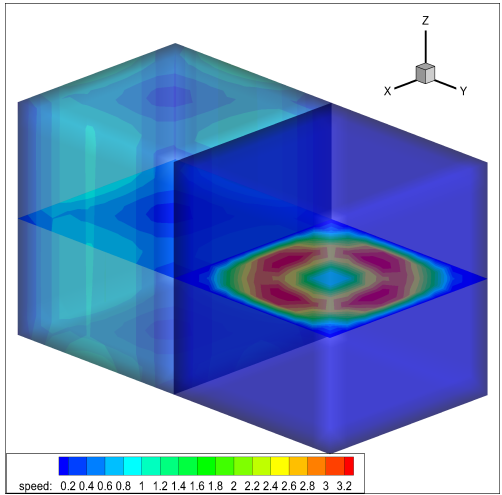}
\end{minipage}
\hfill
\begin{minipage}[t]{0.48\linewidth}
\centering
\includegraphics[height=6.5cm,width=6.5cm]{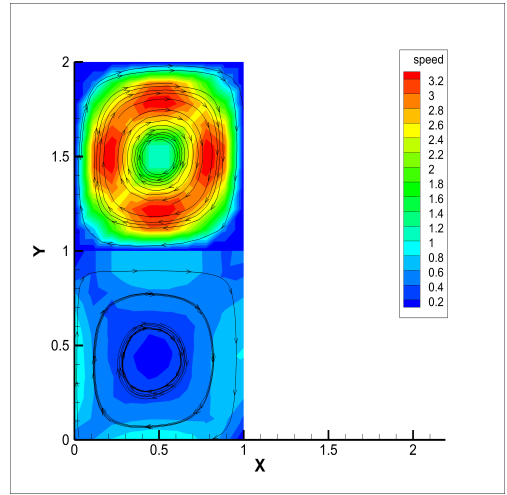}
\end{minipage}
\caption{The streamlines and magnitudes of velocity in a cubical cavity (left) and the cross-section view at $z=0.5$ (right) when $Ra=1.0\times 10^3$.}
\end{figure}

\begin{figure}[!ht]
\centering
\begin{minipage}[t]{0.48\linewidth}
\centering
\includegraphics[height=6.5cm,width=6.5cm]{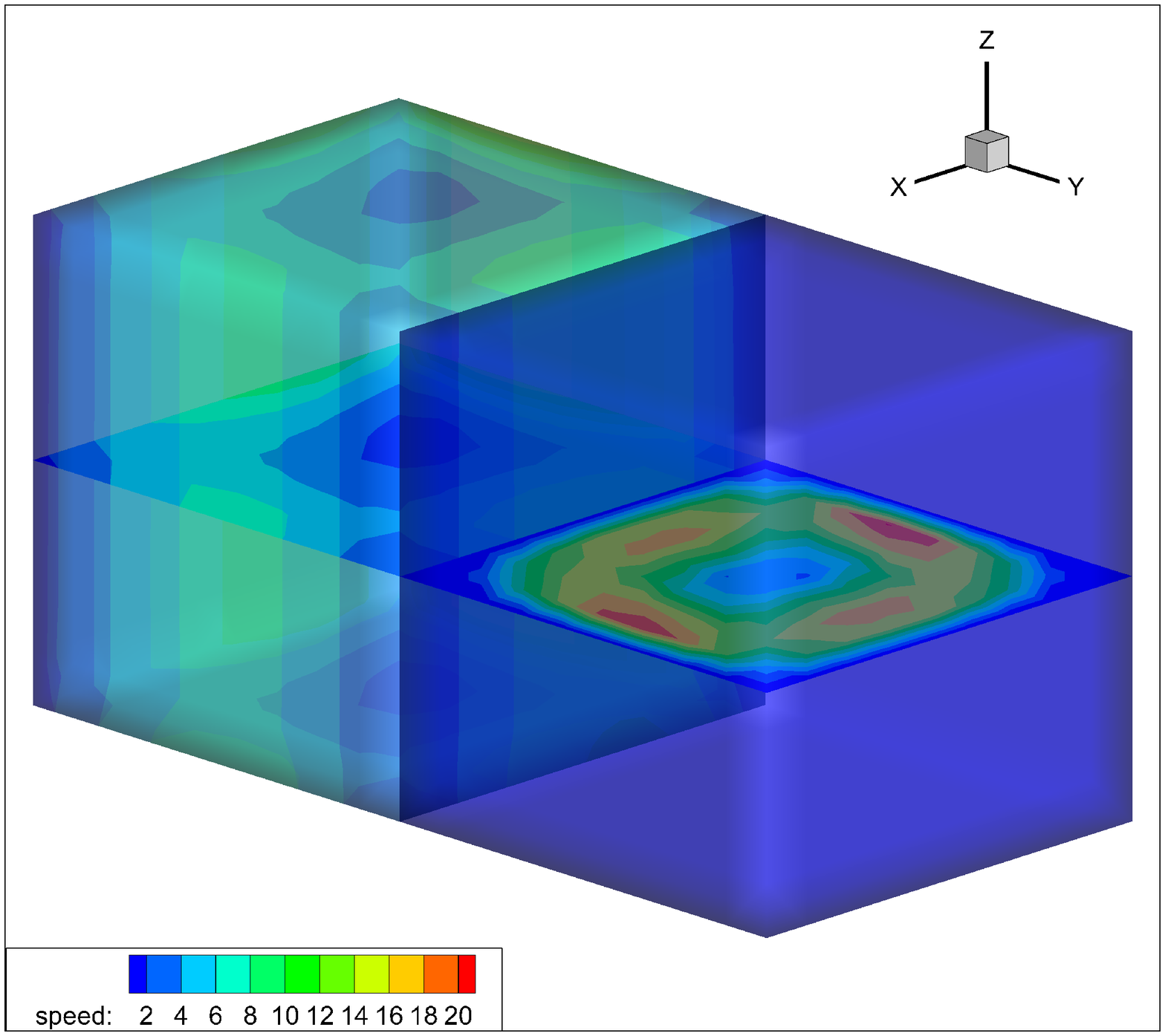}
\end{minipage}
\hfill
\begin{minipage}[t]{0.48\linewidth}
\centering
\includegraphics[height=6.5cm,width=6.5cm]{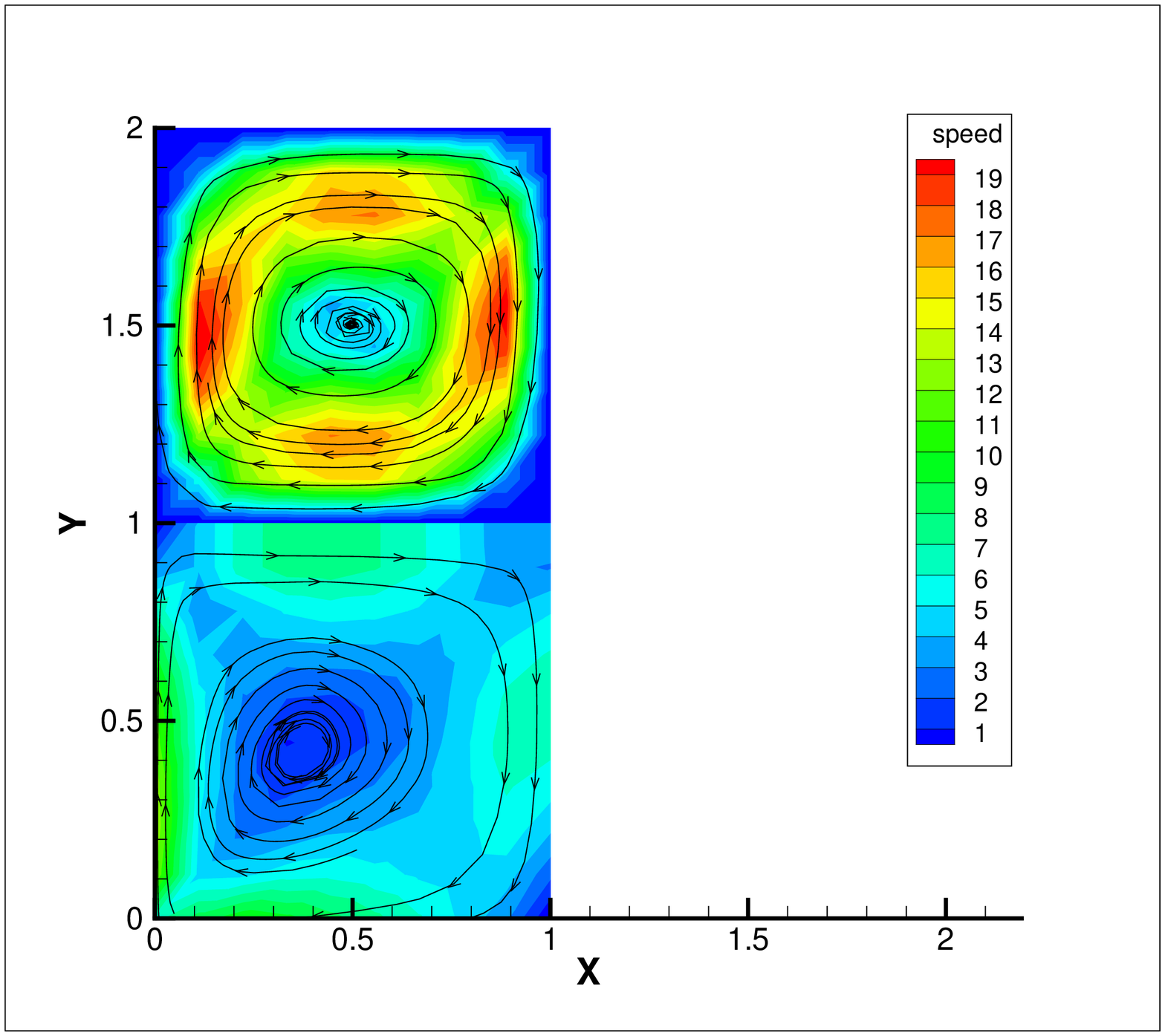}
\end{minipage}
\caption{The streamlines and magnitudes of velocity in a cubical cavity (left) and the cross-section view at $z=0.5$ (right) when $Ra=1.0\times 10^4$.}
\end{figure}

\begin{figure}[!ht]
\centering
\begin{minipage}[t]{0.48\linewidth}
\centering
\includegraphics[height=6.5cm,width=6.5cm]{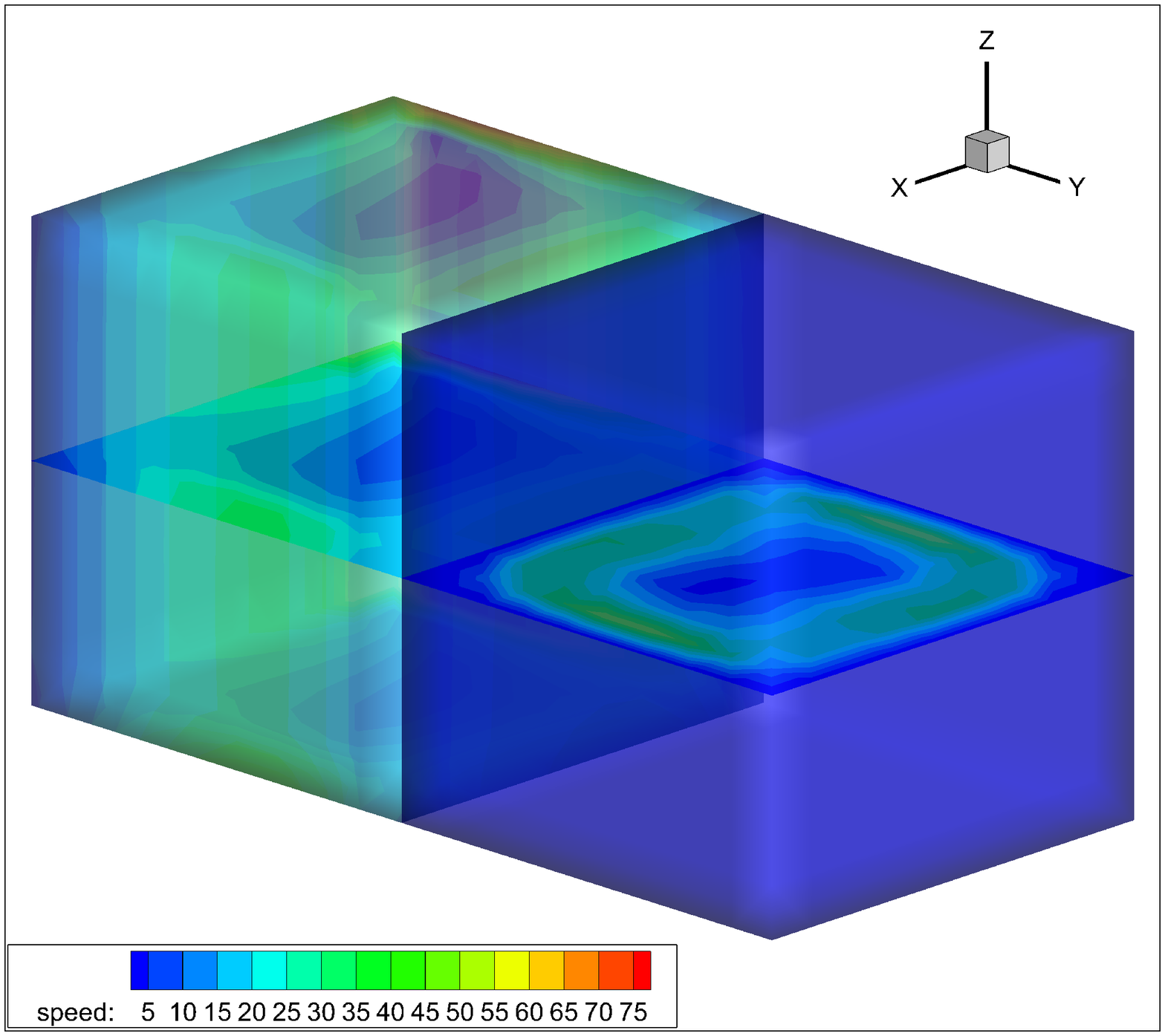}
\end{minipage}
\hfill
\begin{minipage}[t]{0.48\linewidth}
\centering
\includegraphics[height=6.5cm,width=6.5cm]{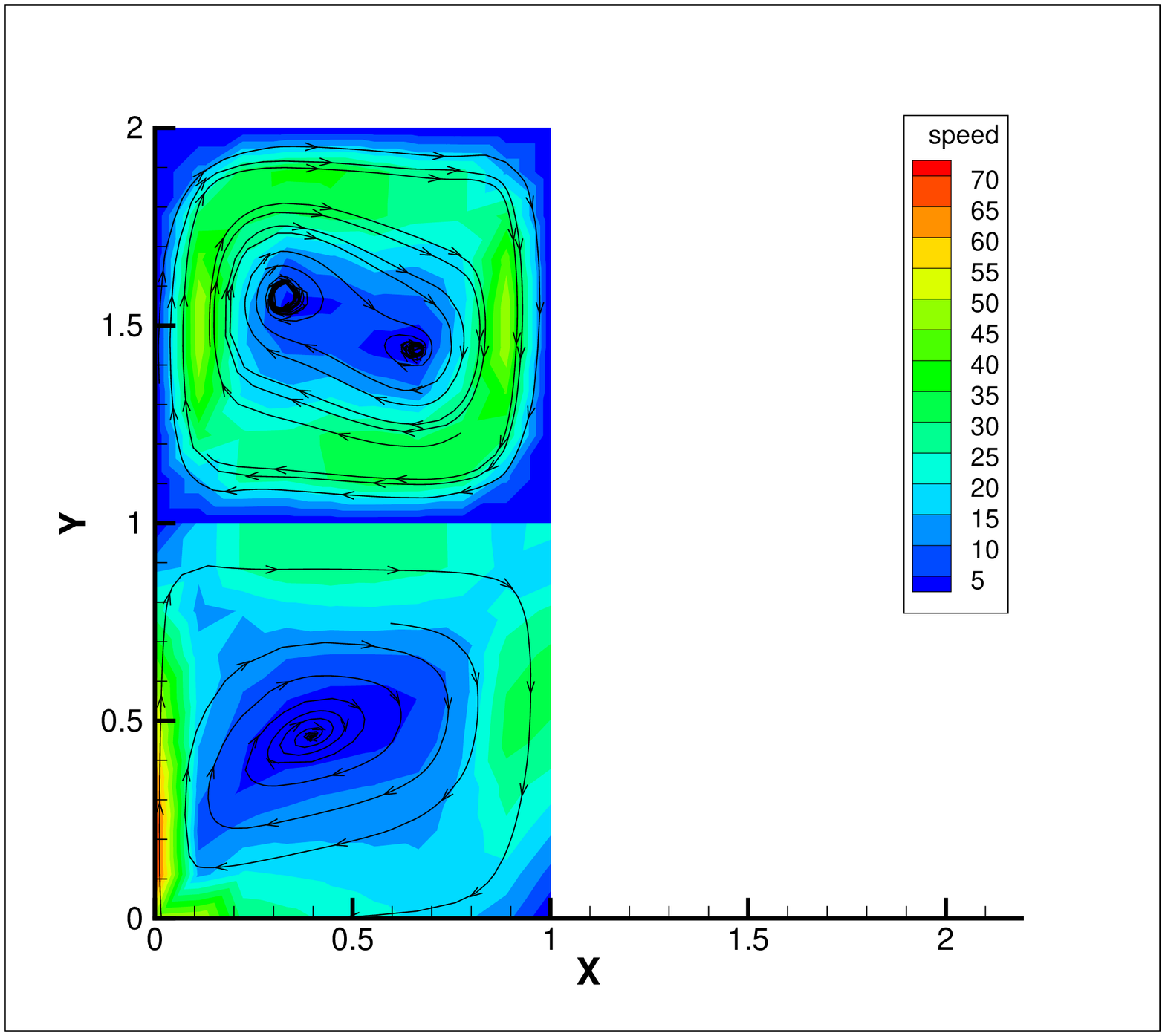}
\end{minipage}
\caption{The streamlines and magnitudes of velocity in a cubical cavity (left) and the cross-section view at $z=0.5$ (right) when $Ra=5.0\times10^4$.}
\end{figure}

\begin{figure}[!ht]
\centering
\begin{minipage}[t]{0.48\linewidth}
\centering
\includegraphics[height=6.5cm,width=6.5cm]{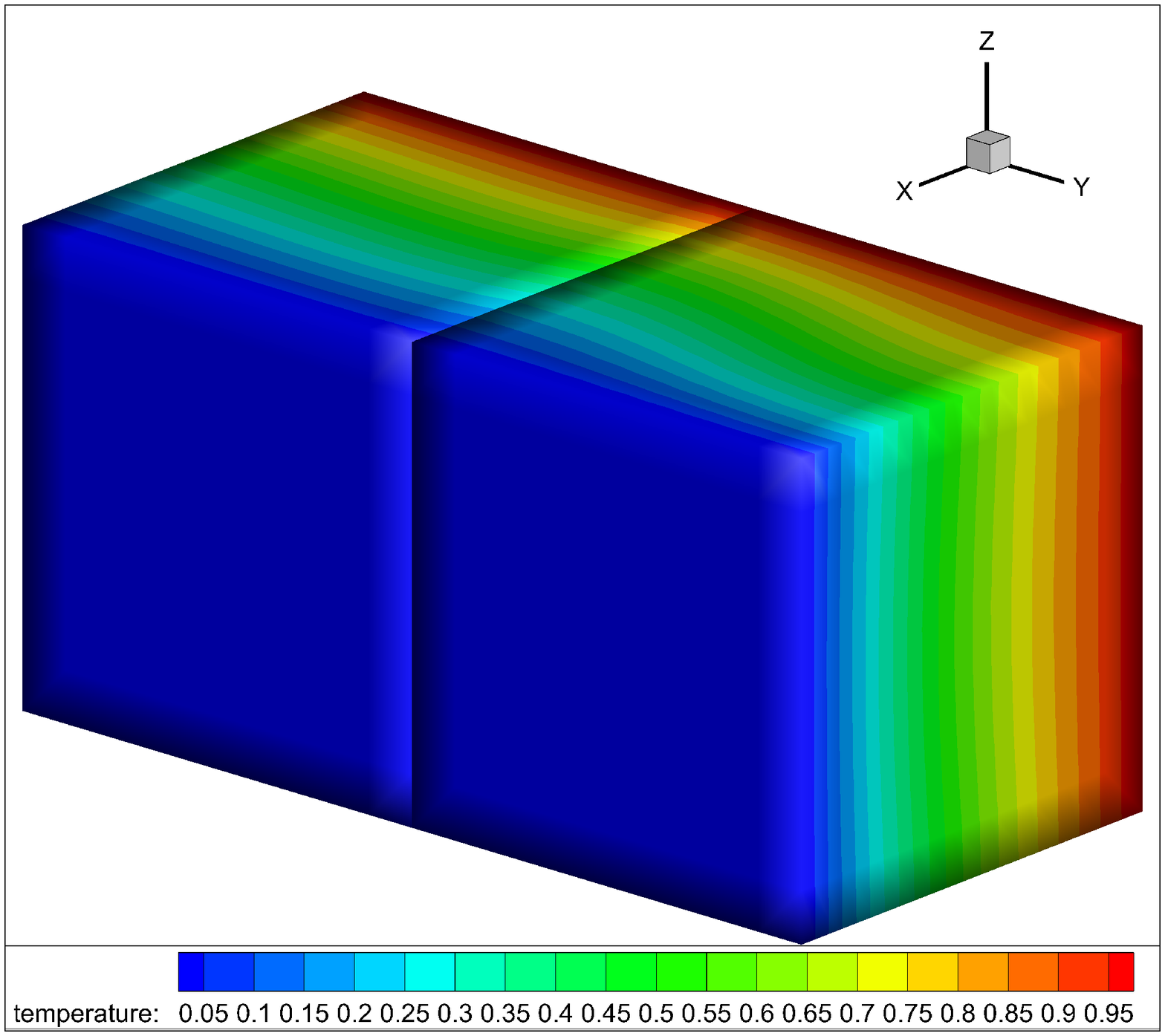}
\end{minipage}
\hfill
\begin{minipage}[t]{0.48\linewidth}
\centering
\includegraphics[height=6.5cm,width=6.5cm]{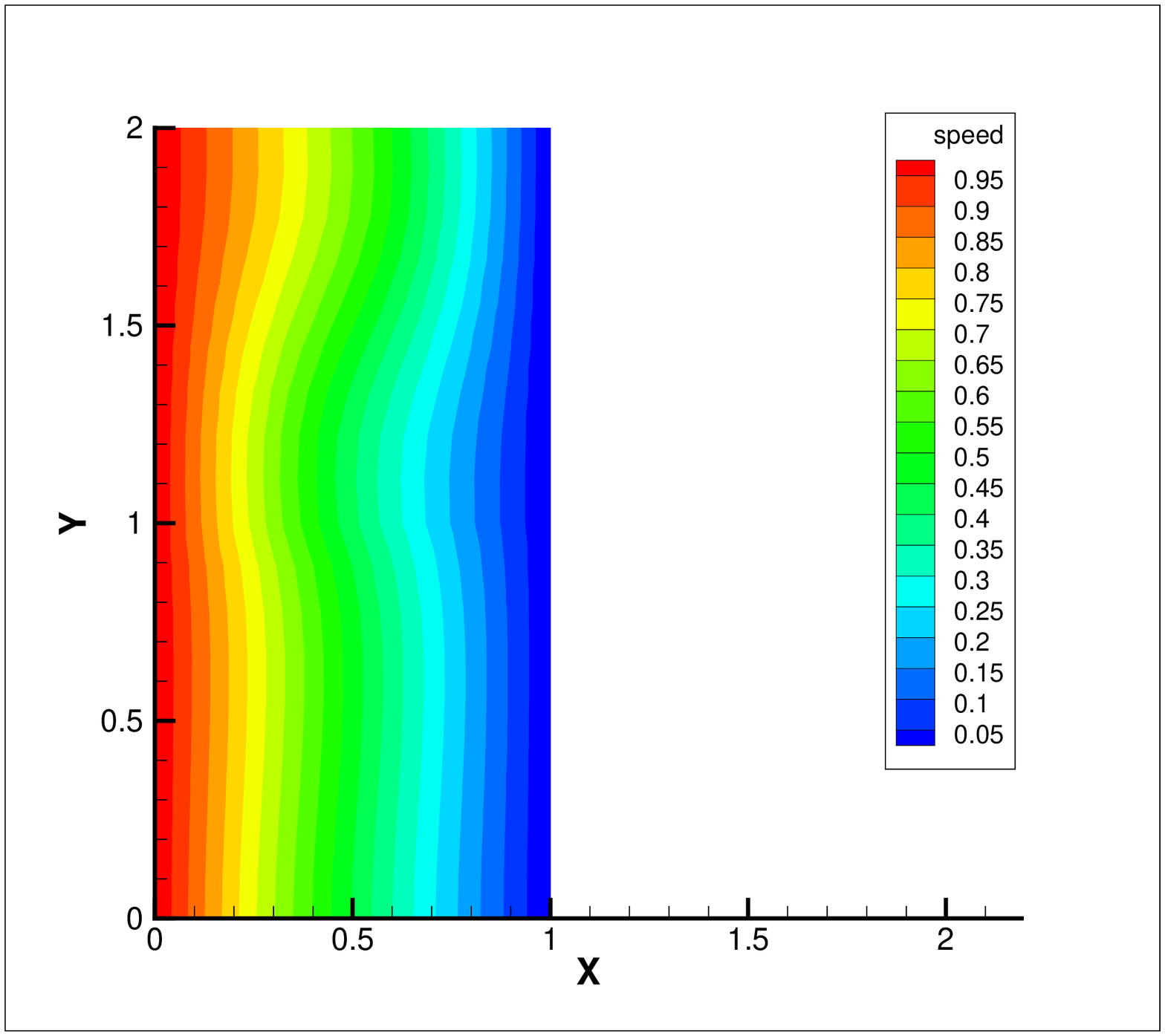}
\end{minipage}
\caption{The temperature distribution in a cubical cavity (left) and the cross-section view at $z=0.5$ (right) when $Ra=1.0\times10^3$.}
\end{figure}

\begin{figure}[!ht]
\centering
\begin{minipage}[t]{0.48\linewidth}
\centering
\includegraphics[height=6.5cm,width=6.5cm]{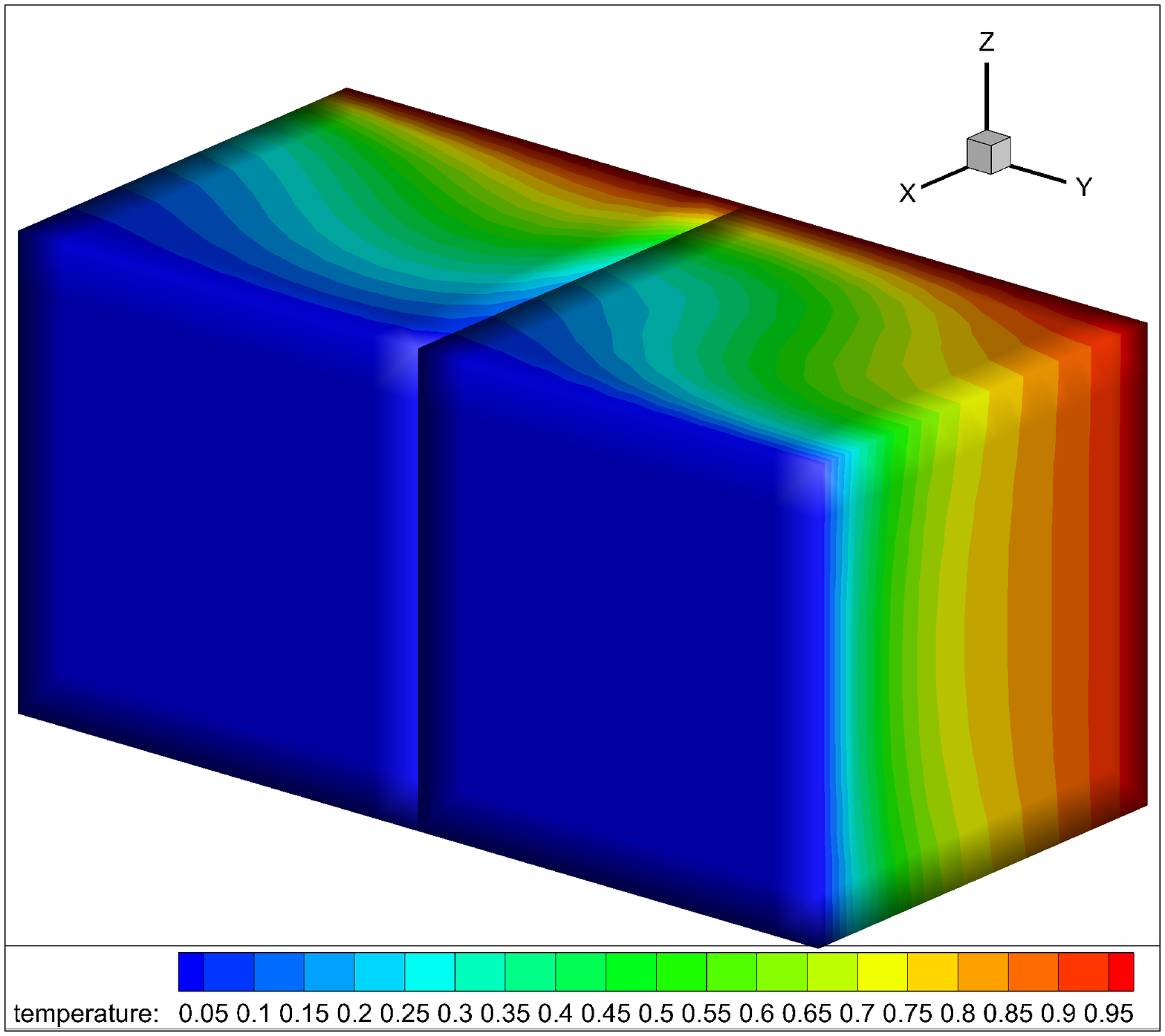}
\end{minipage}
\hfill
\begin{minipage}[t]{0.48\linewidth}
\centering
\includegraphics[height=6.5cm,width=6.5cm]{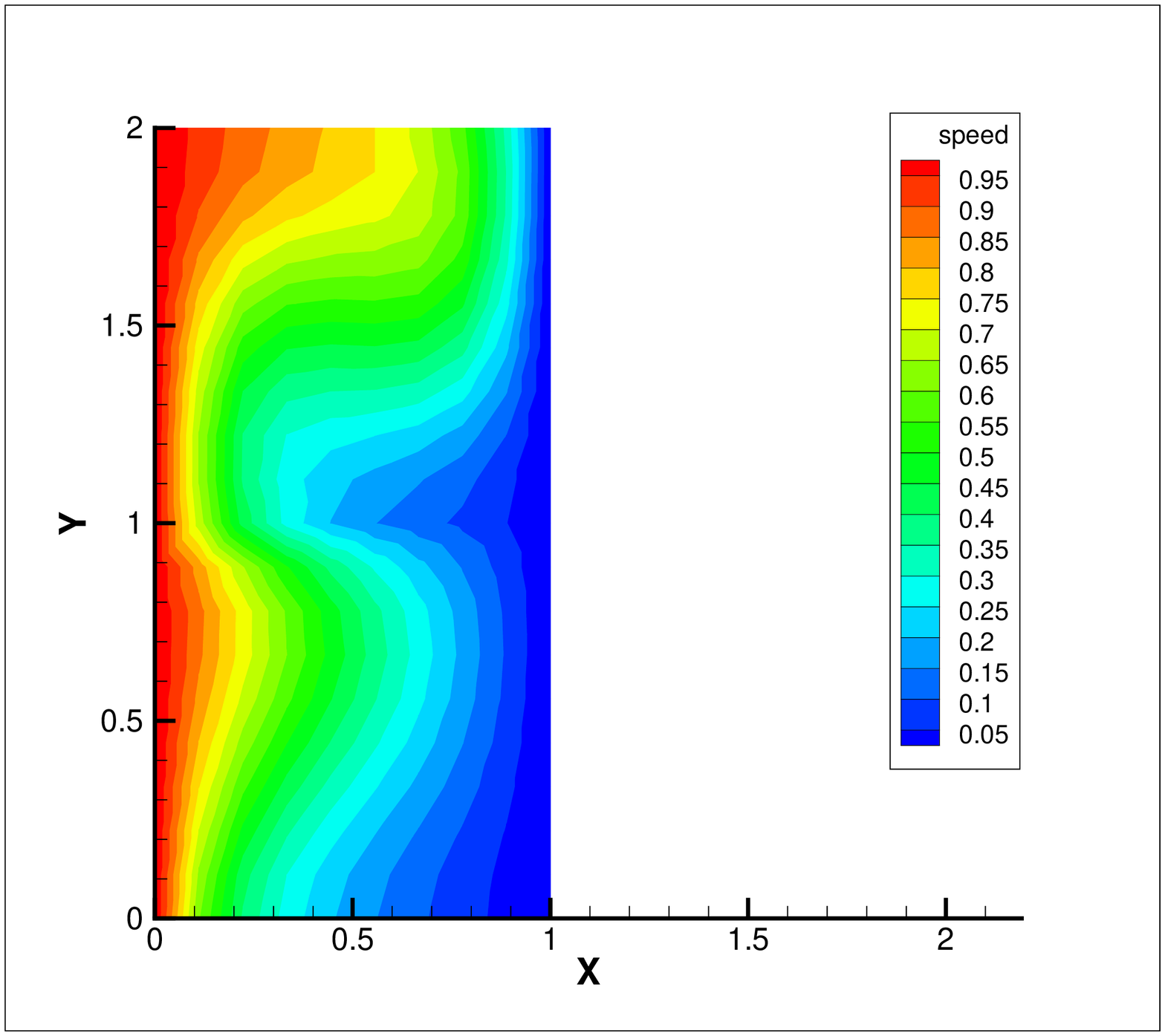}
\end{minipage}
\caption{The temperature distribution in a cubical cavity (left) and the cross-section view at $z=0.5$ (right) when $Ra=1.0\times10^4$.}
\end{figure}

\begin{figure}[!ht]
\centering
\begin{minipage}[t]{0.48\linewidth}
\centering
\includegraphics[height=6.5cm,width=6.5cm]{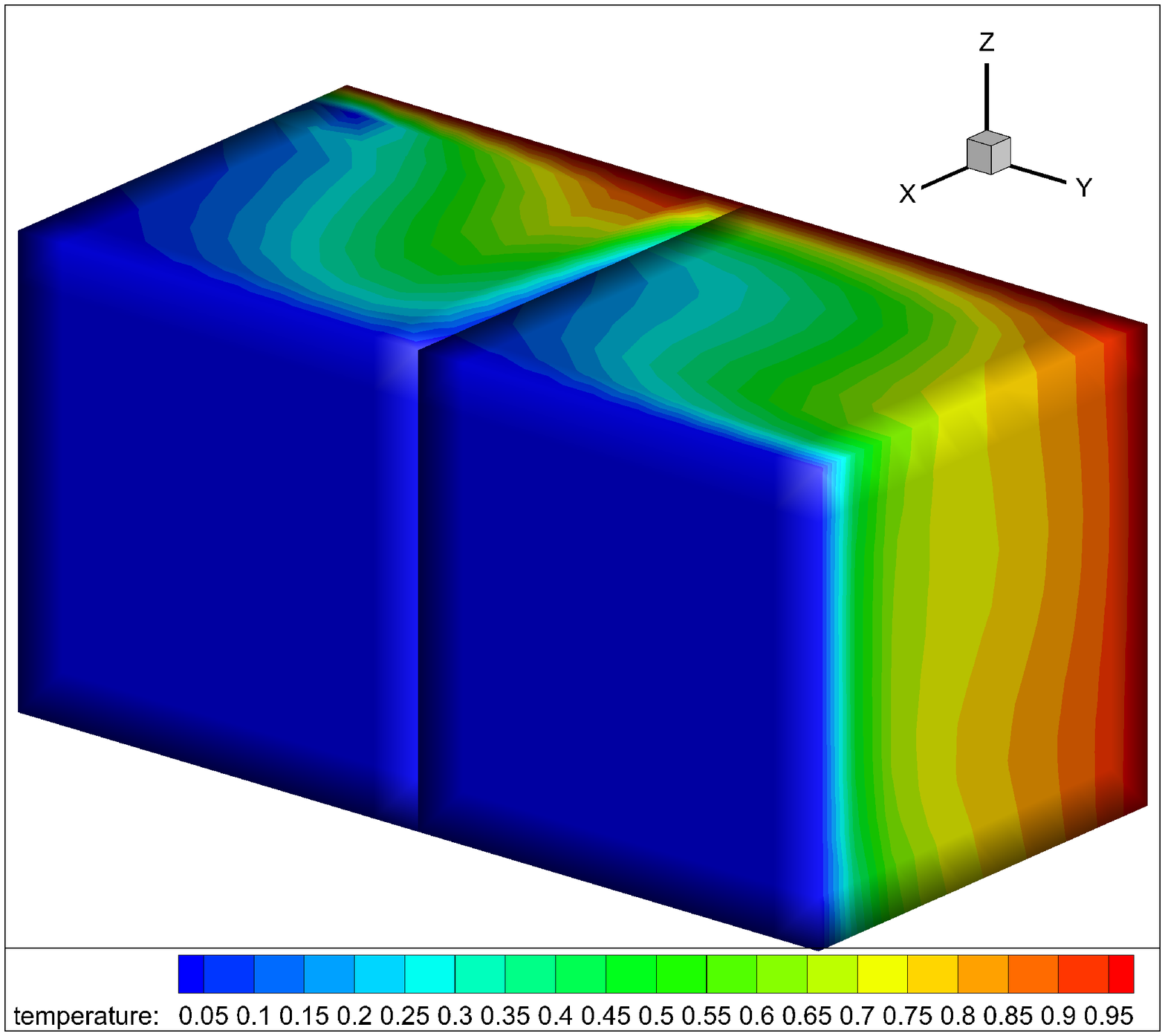}
\end{minipage}
\hfill
\begin{minipage}[t]{0.48\linewidth}
\centering
\includegraphics[height=6.5cm,width=6.5cm]{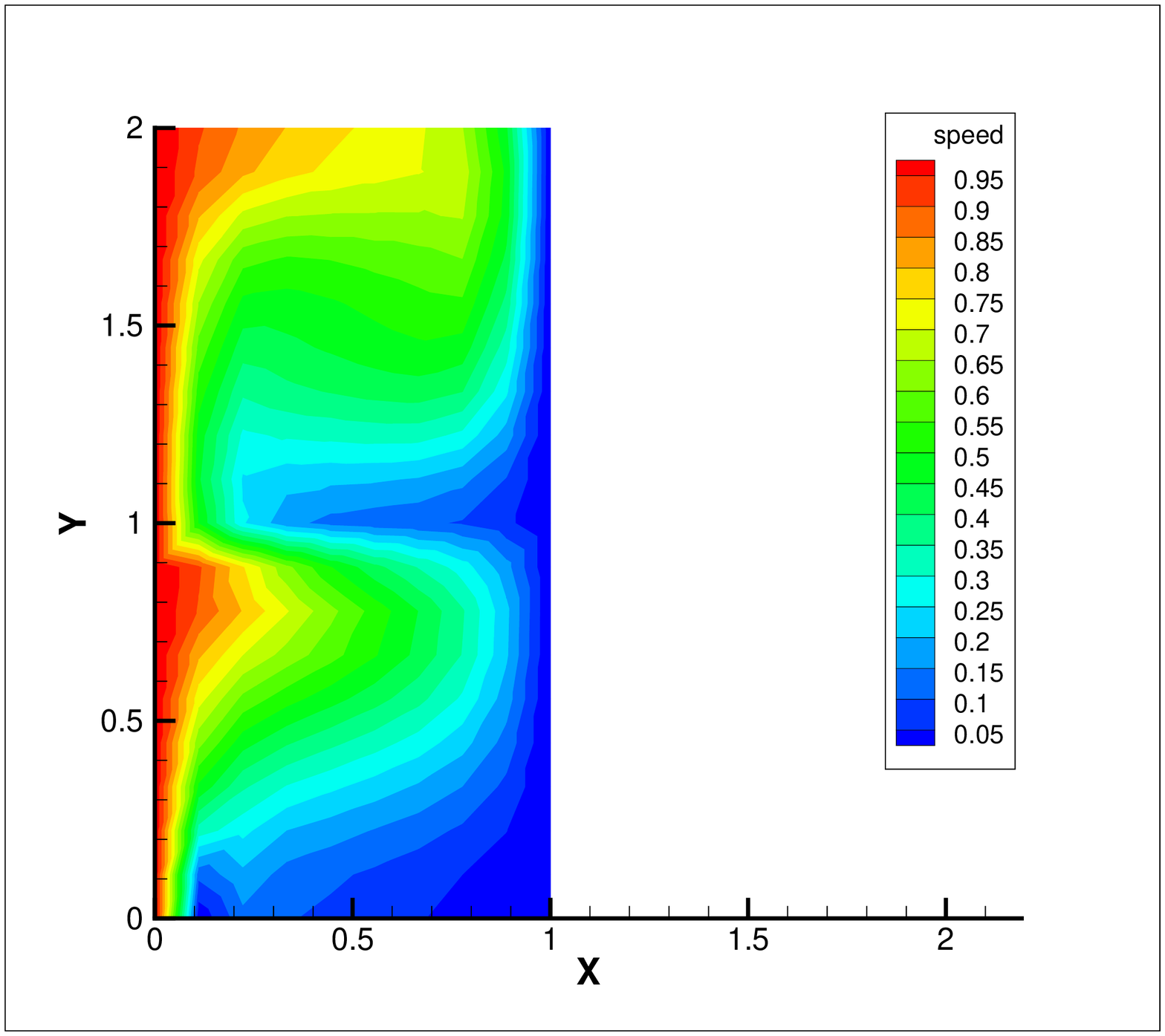}
\end{minipage}
\caption{The temperature distribution in a cubical cavity (left) and the cross-section view at $z=0.5$ (right) when $Ra=5.0\times10^4$.}
\end{figure}

For any $1\leq n\leq N$, we calculate the finite element solution using the classical MCM with $J=100$ sample simulations
\begin{align*}
\begin{split}
&(E[\|\mathbf{u}_{f}\|^2_{L^2(\Omega_{f})}])^{1/2}\approx \bigg(\frac{1}{J} \sum_{l=1}^{J}\|\mathbf{u}_{f}(w_{l})\|^2_{L^2(\Omega_{f})}\bigg)^{1/2},\qquad
(E[\|\mathbf{u}_{p}\|^2_{L^2(\Omega_{p})}])^{1/2}\approx \bigg(\frac{1}{J} \sum_{l=1}^{J}\|\mathbf{u}_{p}(w_{l})\|^2_{L^2(\Omega_{p})}\bigg)^{1/2},\\
&(E[\|\theta_{f}\|^2_{L^2(\Omega_{f})}])^{1/2}\approx \bigg(\frac{1}{J} \sum_{l=1}^{J}\|\theta_{f}(w_{l})\|^2_{L^2(\Omega_{f})}\bigg)^{1/2},\qquad
(E[\|\theta_{p}\|^2_{L^2(\Omega_{p})}])^{1/2}\approx \bigg(\frac{1}{J} \sum_{l=1}^{J}\|\theta_{p}(w_{l})\|^2_{L^2(\Omega_{p})}\bigg)^{1/2}.\\
\end{split}
\end{align*}

Figures 4-9 are numerical results of this example. Observing Figures 4-6, it can be seen from the flow lines that the elliptical vortex at the cavity center breaks up into two vortices as the Rayleigh number increases, as studied in \cite{D C Wan-2001,Z Y Zhang-2013,Z Y Zhang-2014}. As can be seen from Figures 7-9, along with the increase of Rayleigh numbers, the temperature convection becomes increasingly prominent. The proposed model and numerical method are validated by these physically valid simulation results.
\subsubsection{Simulation for a 3D stochastic closed-loop geothermal system by a closed U-tube}
\begin{figure}[!ht]
\centering
\begin{minipage}[t]{0.48\linewidth}
\centering
\includegraphics[height=6.5cm,width=6.5cm]{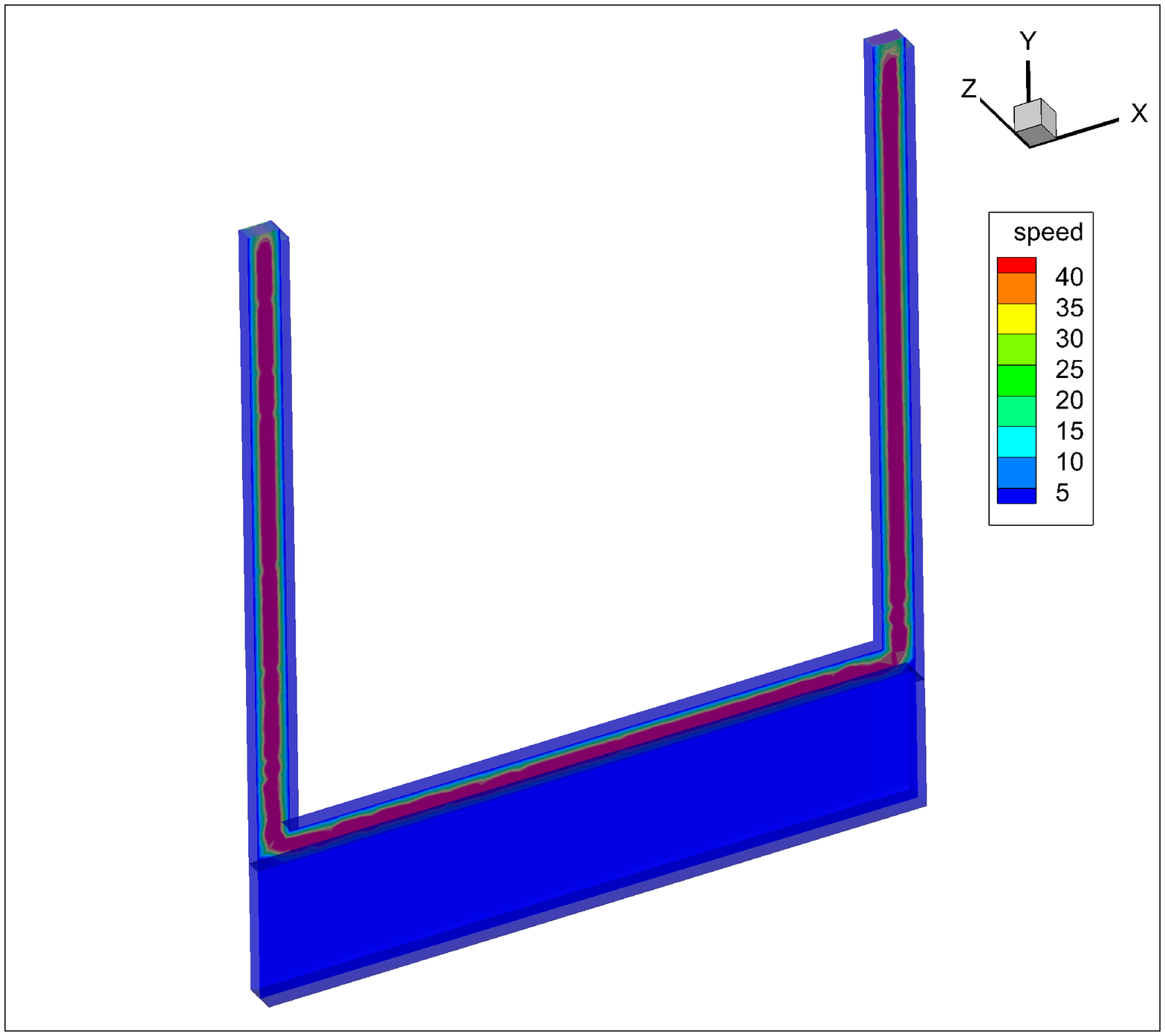}
\end{minipage}
\hfill
\begin{minipage}[t]{0.48\linewidth}
\centering
\includegraphics[height=6.5cm,width=6.5cm]{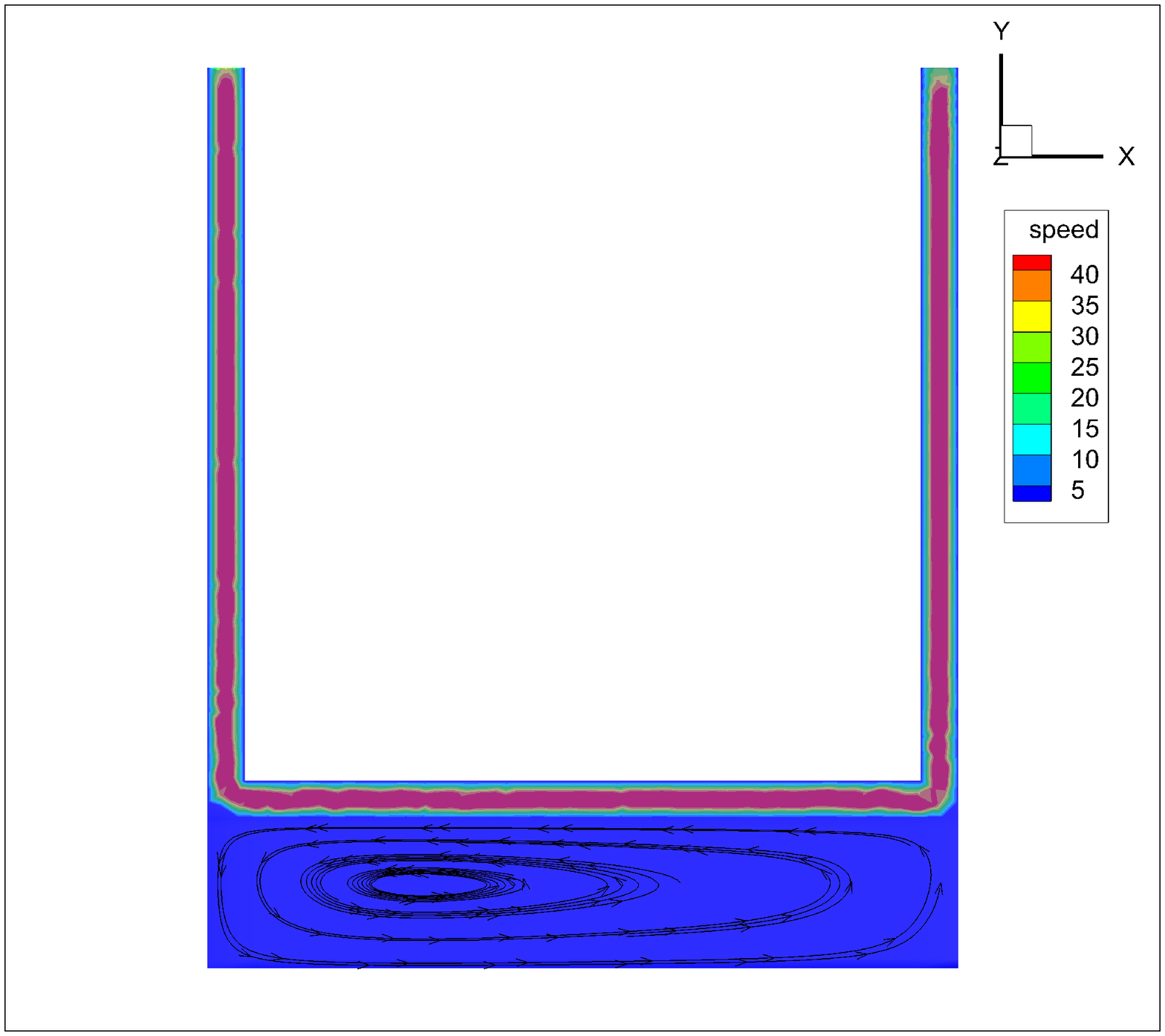}
\end{minipage}
\caption{The streamlines and magnitudes of velocity in a 3D U-shape pipeline (left) and the cross-section view at $z = 0.125$ (right) with $J=100$.}
\end{figure}

\begin{figure}[!ht]
\centering
\begin{minipage}[t]{0.48\linewidth}
\centering
\includegraphics[height=6.5cm,width=6.5cm]{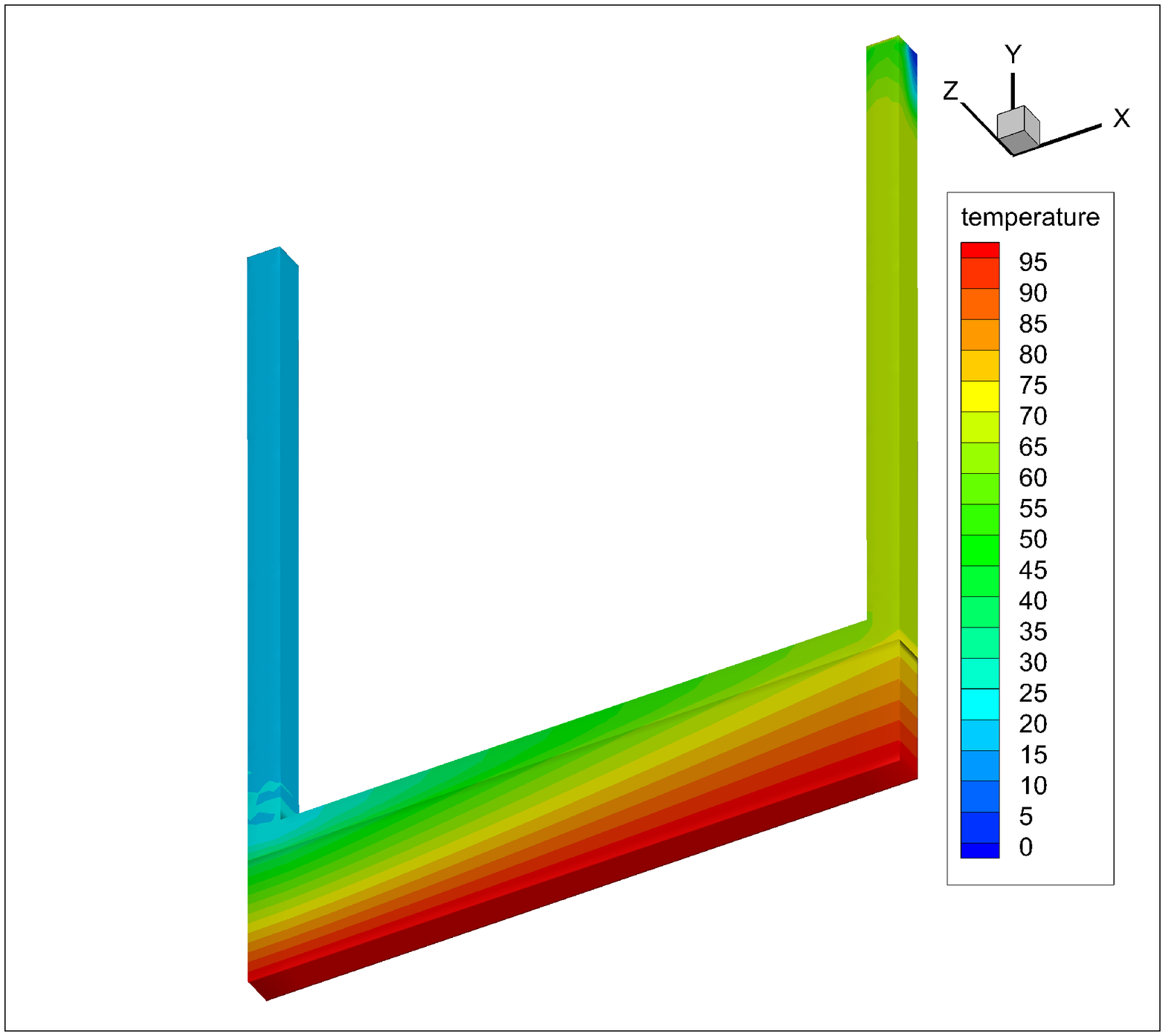}
\end{minipage}
\hfill
\begin{minipage}[t]{0.48\linewidth}
\centering
\includegraphics[height=6.5cm,width=6.5cm]{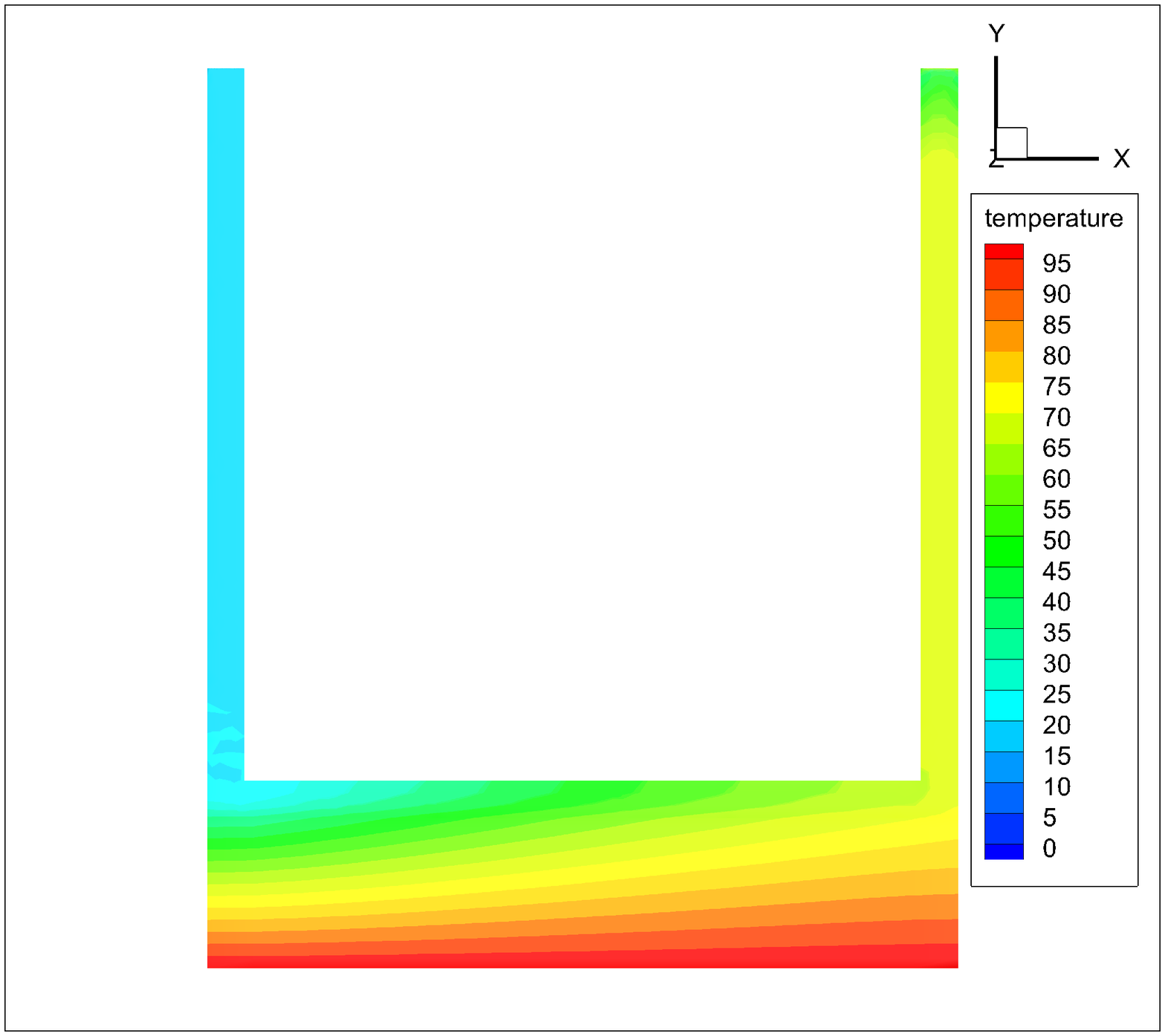}
\end{minipage}
\caption{The temperature distribution in a 3D U-shape pipeline (left) and the cross-section view at $z=0.18$ (right) with $J=100$.}
\end{figure}

In this numerical experiment, we simulate the velocity distribution and thermal convection of a 3D stochastic closed-loop geothermal system.
The pipeline region consists of a 3D U-tube closed-loop heat exchange pipeline that injects cold fluid from a vertical injection well on the left, and after passing through the geothermal reservoir, the hot fluid is pumped out from a vertical production well on the right \cite{W. Choi-2016, C. M. Oldenburg-2016}.

For the boundaries $\partial\Omega_{fin}=\{(x,y,z):y=6,0\leq x\leq0.25,0\leq z\leq0.25\}$, $\partial\Omega_{fout}=\{(x,y,x):y=6,4.75\leq x\leq5,0\leq z\leq0.25\}$, $\{(x,y,z):x=0, 1\leq y\leq6,0\leq z\leq0.25\}$, $\{(x,y,z):x=0.25, 1.25\leq y\leq6,0\leq z\leq0.25\}$, $\{(x,y,z):y=1.25, 0.25\leq x\leq4.75,0\leq z\leq0.25\}$, $\{(x,y,z):x=4.75, 1.25\leq y\leq6,0\leq z\leq0.25\}$, $\{(x,y,z):x=5, 1\leq y\leq6,0\leq z\leq0.25\}$ and the area of the U-shaped pipe when $z=0$ and $z=0.25$ of the closed-loop pipe $\Omega_{f}$, we impose the velocity and temperature boundary conditions as
\begin{align*}
\begin{split}
\begin{cases}
&U_{x}=0,\ U_{y}= -2048.0x(0.25-x),\ U_{z}=0,\ \theta_{f}=20.0\qquad on\  \partial\Omega_{fin},\\
&(-p_{f}\mathbf{I}+Pr\nabla\mathbf{u}_{f})\cdot\mathbf{n}_{f}=0,\ \mathbf{n}_{f}\cdot k_{f}\nabla\theta_{f}=0\qquad on\  \partial\Omega_{fout},\\
&\mathbf{u}_{f}=\mathbf{0},\ \mathbf{n}_{f}\cdot k_{f}\nabla\theta_{f}=0\qquad
on\ \partial\Omega_{f}\setminus\uppercase\expandafter{\romannumeral1}\setminus\partial\Omega_{fin}\setminus\partial\Omega_{fout}.
\end{cases}
\end{split}
\end{align*}

On the interface $\uppercase\expandafter{\romannumeral1}=\{(x,y,z):y=1,0\leq x\leq5,0\leq z\leq0.25\}$, use the interface conditions (\ref{i1}) - (\ref{i4}) mentioned earlier. The geothermal reservoir domain $\Omega_{p}=(0,5)\times(0,1)\times(0,0.25)$. The no-flow boundary condition $\mathbf{u}_{p}\cdot\mathbf{n}_{p}=0$ is imposed on $\partial\Omega_{p}\setminus\uppercase\expandafter{\romannumeral1}$. The homogeneous Neumann boundary condition is considered for the temperature on the left and right walls of $\Omega_{p}$. We impose the boundary condition $\theta_{p}=100.0$ at the bottom of the $\Omega_{p}$.

Assume that the initial conditions of temperature and velocity in the closed-loop pipe are $\mathbf{u}_{f}=(0,0,0)$ and $\theta_{f}=20.0$, respectively. Consider $\mathbf{u}_{p}=(0,0,0)$ and $\theta_{p}=100.0$ as the initial conditions for the velocity and temperature of $\Omega_{p}$, respectively. Moreover, we take $Pr=3.0,\ Ra=3.0\times 10^3,\ k_{f}=0.6,\ k_{p}=C_{a}=1.0, \ L=1.0\times 10^3,\ \gamma=1.0,\ T=3.0$ and $\Delta t=0.01$. We choose the grid size $h_{f}=\frac{1}{49}$ and $\ h_{p}=\frac{1}{16}$. The selection of the random hydraulic conductivity tensor is consistent with experiment 5.2.1.

We consider a group of simulations with $J = 100$ using the MCM for sampling. We can observe from Figure 10 that the no-fluid communication interface conditions are reflected in the velocity streamlines. Figure 11 shows that the temperature distribution is reasonable, the cold fluid injected from the vertical injection well, heated through the geothermal reservoir, and hot fluid is pumped from the production well.
\subsubsection{Simulation for a 3D stochastic closed-loop geothermal system by a closed Co-axial}
~\par Like the U-tube design, the Co-axial configuration transfers heat from the geothermal reservoir to the well by conduction only. And, as the central 'Up' pipe is encased by the 'Down' pipe, greater insulation as the heated water travels to the surface \cite{Y.T. He-2021}. In this numerical experiment, we simulate the velocity distribution and thermal convection of a 3D stochastic closed-loop geothermal system by a closed Co-axial.

For the boundaries $\partial\Omega_{fin1}=\{(x,y,z):y=1,0\leq x\leq0.2,0\leq z\leq0.25\}$, $\partial\Omega_{fin2}=\{(x,y,z):y=1,0.4\leq x\leq0.6,0\leq z\leq0.25\}$, $\partial\Omega_{fout}=\{(x,y,x):y=1,0.2\leq x\leq0.4,0\leq z\leq0.25\}$, $\{(x,y,z):x=0.6, 0.1\leq y\leq1,0\leq z\leq0.25\}$, $\{(x,y,z):x=0, 1\leq y\leq0.1,0\leq z\leq0.25\}$, $\{(x,y,z):x=0.4, 0.3\leq x\leq1,0\leq z\leq0.25\}$, $\{(x,y,z):x=0.2, 0.3\leq y\leq1,0\leq z\leq0.25\}$, and the area of the closed Co-Axial pipe when $z=0$ and $z=0.25$ of the closed-loop pipe $\Omega_{f}$, we impose the velocity and temperature boundary conditions as
\begin{align*}
\begin{split}
\begin{cases}
&U_{x}=0,\ U_{y}= 1000.0(x-0.1)^{2}-10,\ U_{z}=0,\ \theta_{f}=20.0\qquad on\  \partial\Omega_{fin1},\\
&U_{x}=0,\ U_{y}= 1000.0(x-0.5)^{2}-10,\ U_{z}=0,\ \theta_{f}=20.0\qquad on\  \partial\Omega_{fin2},\\
&(-p_{f}\mathbf{I}+Pr\nabla\mathbf{u}_{f})\cdot\mathbf{n}_{f}=0,\ \mathbf{n}_{f}\cdot k_{f}\nabla\theta_{f}=0\qquad on\  \partial\Omega_{fout},\\
&\mathbf{u}_{f}=\mathbf{0},\ \mathbf{n}_{f}\cdot k_{f}\nabla\theta_{f}=0\qquad
on\ \partial\Omega_{f}\setminus\uppercase\expandafter{\romannumeral1}\setminus\partial\Omega_{fin}\setminus\partial\Omega_{fout}.
\end{cases}
\end{split}
\end{align*}

On the interface $\uppercase\expandafter{\romannumeral1}=\{(x,y,z):y=0.1,0\leq x\leq0.6,0\leq z\leq0.25\}$, use the interface conditions (\ref{i1}) - (\ref{i4}) mentioned earlier. The geothermal reservoir domain $\Omega_{p}=(0,0.6)\times(0,0.1)\times(0,0.25)$. The no-flow boundary condition $\mathbf{u}_{p}\cdot\mathbf{n}_{p}=0$ is imposed on $\partial\Omega_{p}\setminus\uppercase\expandafter{\romannumeral1}$. The homogeneous Neumann boundary condition is considered for the temperature on the left and right walls of $\Omega_{p}$. We impose the boundary condition $\theta_{p}=120.0$ at the bottom of the $\Omega_{p}$.

Assume that the initial conditions of temperature and velocity in the closed-loop pipe are $\mathbf{u}_{f}=(0,0,0)$ and $\theta_{f}=20.0$, respectively. Consider $\mathbf{u}_{p}=(0,0,0)$ and $\theta_{p}=120.0$ as the initial conditions for the velocity and temperature of $\Omega_{p}$, respectively. The selection of all parameters and the random hydraulic conductivity tensor is consistent with experiment 5.2.2.
The mesh size $h_{max} = 1/16$ and $Ra=1.0\times 10^3$.

Similarly, we consider a group of simulations with $J = 100$ using the MCM for sampling. Figure 12 shows that the no-fluid communication interface conditions are reflected in the velocity streamlines. We can observe from Figure 13 that the cold flow is injected from the injection Wells on both sides, and after absorbing heat from the geothermal reservoir, the heat flow is pumped out along the inner pipe.

\begin{figure}[!ht]
\centering
\begin{minipage}[t]{0.48\linewidth}
\centering
\includegraphics[height=6.5cm,width=6.5cm]{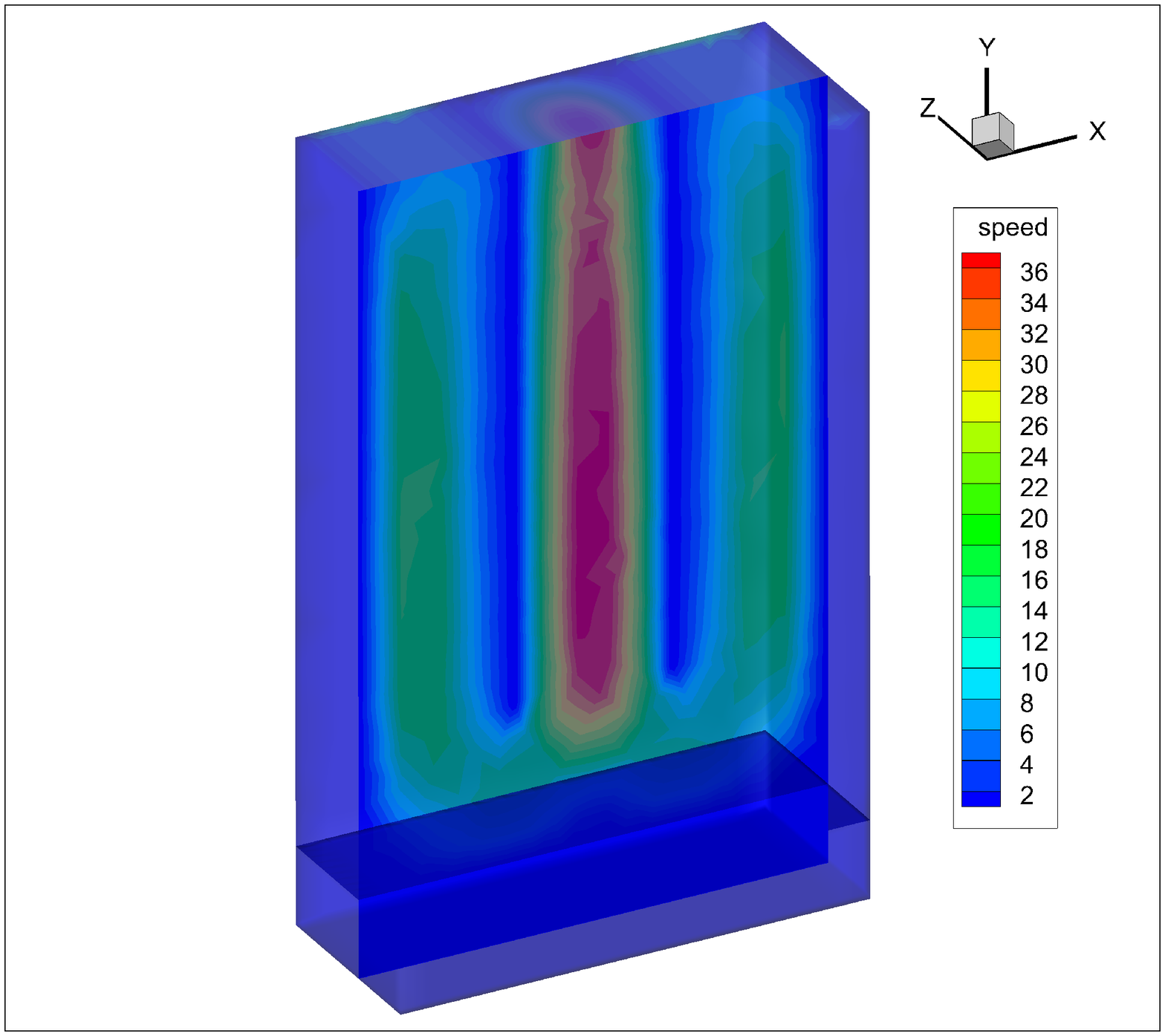}
\end{minipage}
\hfill
\begin{minipage}[t]{0.48\linewidth}
\centering
\includegraphics[height=6.5cm,width=6.5cm]{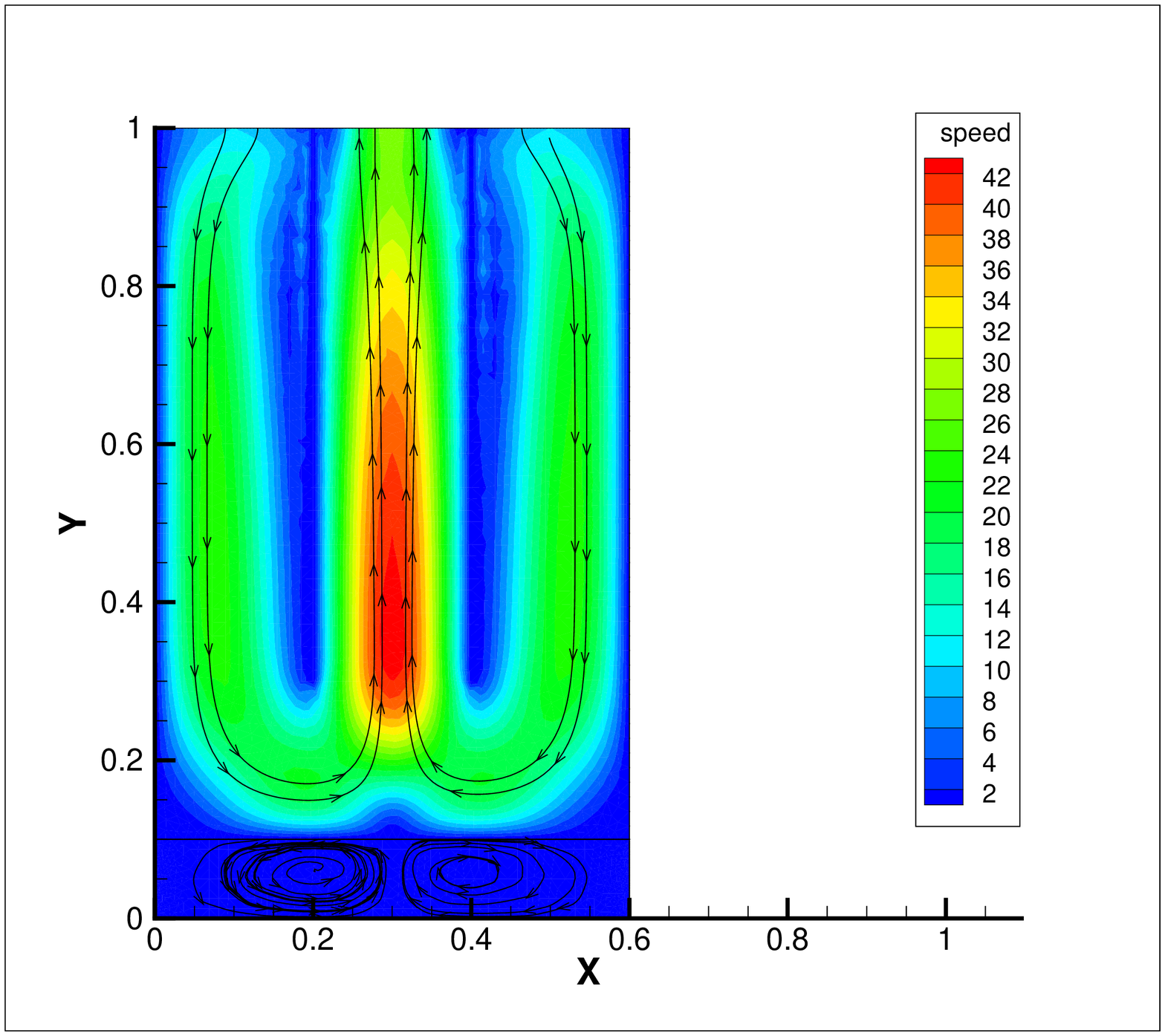}
\end{minipage}
\caption{The streamlines and magnitudes of velocity in a 3D closed Co-axial pipeline (left) and the cross-section view at $z = 0.125$ (right) with $J=100$.}
\end{figure}

\begin{figure}[!ht]
\centering
\begin{minipage}[t]{0.48\linewidth}
\centering
\includegraphics[height=6.5cm,width=6.5cm]{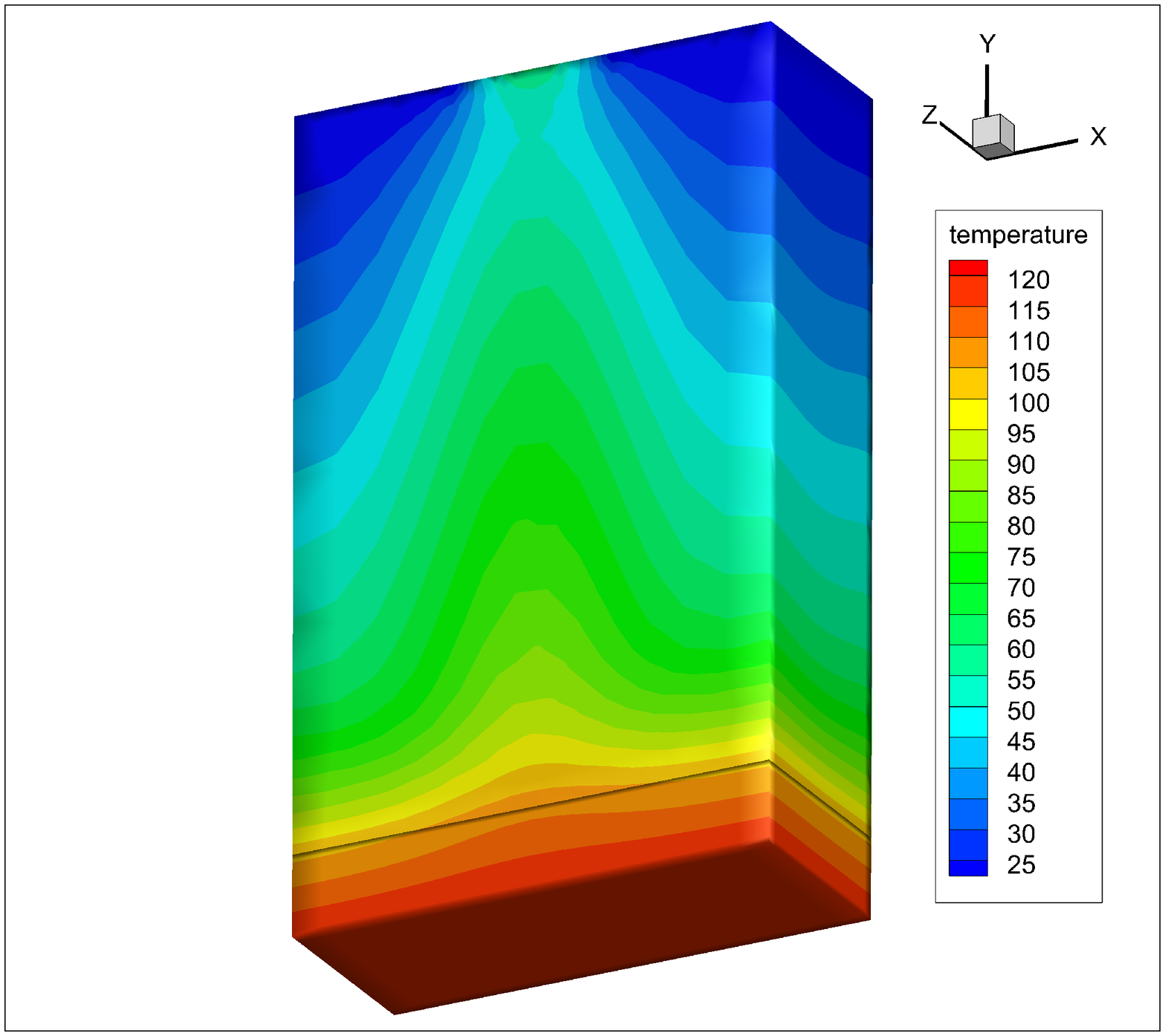}
\end{minipage}
\hfill
\begin{minipage}[t]{0.48\linewidth}
\centering
\includegraphics[height=6.5cm,width=6.5cm]{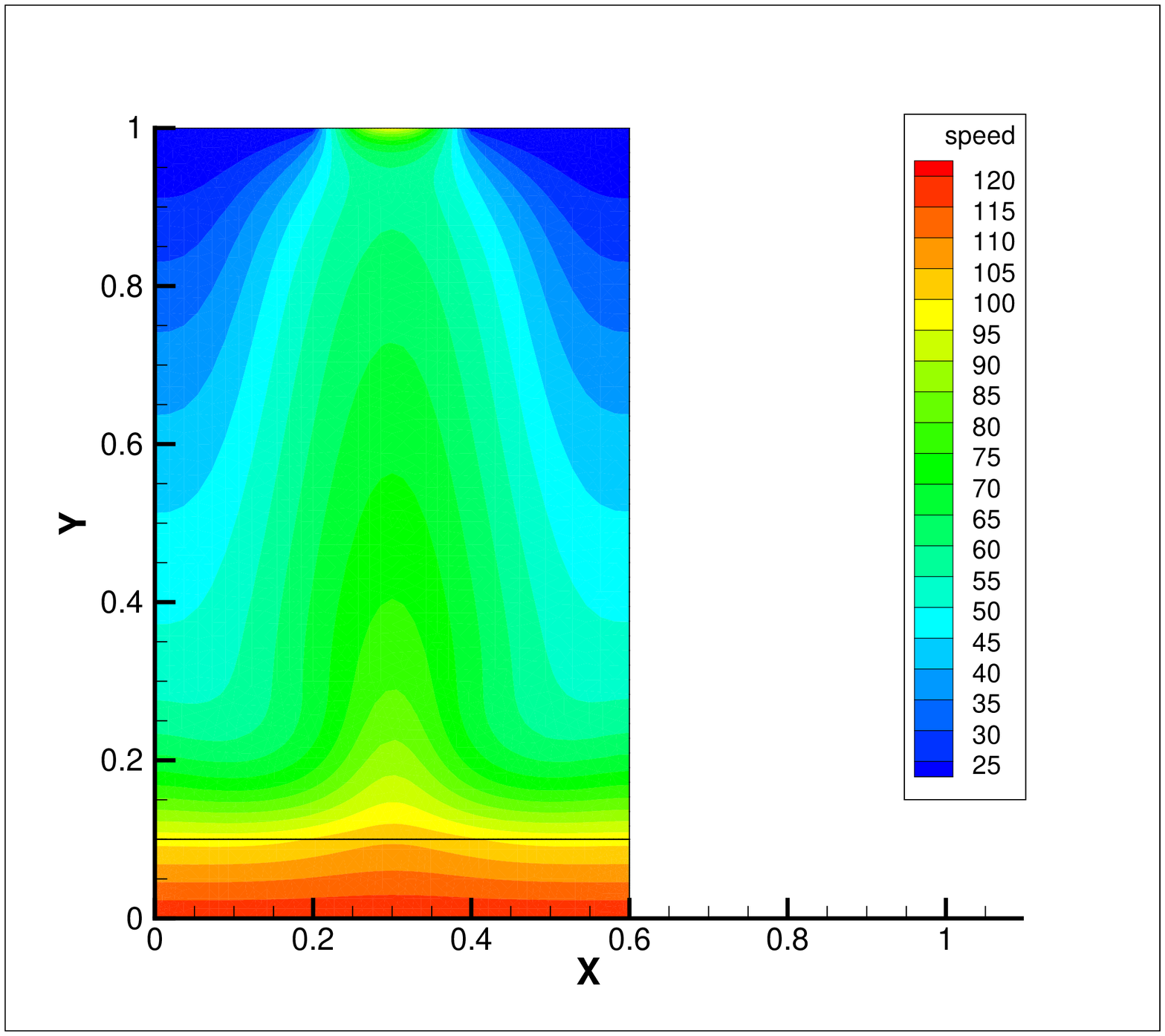}
\end{minipage}
\caption{The temperature distribution in a 3D closed Co-axial pipeline (left) and the cross-section view at $z=0.18$ (right) with $J=100$.}
\end{figure}
\section*{Declaration of competing interest}
The authors declare that they have no known competing financial interests or personal relationships.
\section*{Declaration of data availability}
The datasets generated during and analysed during the current study are available from the corresponding author on reasonable request.

\end{document}